\documentclass[12pt]{article}

\usepackage{amssymb}
\usepackage{latexsym}
\usepackage{amsfonts}
\newtheorem{thm}{Theorem}[section]
\newtheorem{cor}[thm]{Corollary}
\newtheorem{lem}[thm]{Lemma}
\newtheorem{pro}[thm]{Proposition}
\newtheorem{defn}[thm]{Definition}

\title{Circuits, coNP-completeness, and the groups of Richard Thompson }

\author{ Jean-Camille Birget
  \thanks{Research supported in part by NSF grant DMS-9970471, 
          and in part by NSERC grant 216872-1999}
       }
\date{}
\begin{document}
\maketitle

\begin{abstract}
We construct a finitely presented group with coNP-complete word problem, 
and a finitely generated simple group with coNP-complete word problem. 
These groups are represented as Thompson groups, hence as partial 
transformation groups of strings.  
The proof provides a simulation of combinational circuits by elements of
the Thompson-Higman group $G_{3,1}$. 
\end{abstract}

%%%%%%%%%%%%%%%%%%%%%%%%%%%%%%%%%%%%%%%%%%%%%%%%%%%%%%%%
% Section 1
%%%%%%%%%%%%%%%%%%%%%%%%%%%%%%%%%%%%%%%%%%%%%%%%%%%%%%%%

\section{Introduction}

There are many open problems in computational complexity, e.g., the famous
questions ``P $\neq$ NP ?'', and ``NP $\neq$ coNP ?'', that are believed 
to be very difficult. One way to approach very difficult problems is to 
relate them to other disciplines. For computational complexity 
there are interesting relations with combinatorial group theory.
An early connection was Max Dehn's formulation of the word problem of a
group (1910). It took 45 years until it was shown that there is a finitely
presented group whose word problem is undecidable, and that certain finitely 
presented groups can simulate universal Turing machines (Novikov 1955, Boone 
1954-57). Soon after, Higman's embedding theorem (1961) gave an algebraic 
characterization of recursive enumerability of the word problem of a 
group $G$ (namely, $G$ has a recursively enumerable word problem iff $G$ is
isomorphic to a subgroup of some finitely presented group). Boone and Higman
(1976) gave an algebraic characterization of decidability of the word problem
of a group $G$ (namely, $G$ has a decidable word problem iff $G$ is 
isomorphic to a subgroup of some simple group, which itself is a subgroup 
of some finitely presented group). It was also proved that some finitely
presented groups have a primitive recursive word problem; in fact, Madlener 
and Otto \cite{MO} gave a version of the Higman embedding theorem that 
preserves the Grzegorczyk hierarchy from level 3 upward. Madlener and
Otto also introduced what was later called the isoperimetric function
of a group.

It has long been folklore knowledge that (un)decidability, recursive
enumerability, primitive recursiveness, and the Grzegorczyk level, of the
word problem of a finitely presented group $G$ is an algebraic property of
the group, i.e., if one changes over to a different finite set of generators
of the same group, the property is preserved. 
Madlener and Otto showed that the isoperimetric function of a group changes 
only linearly when one changes the finite presentation of the group. 
A similar argument shows that the computational complexity (time or space,
deterministic, nondeterministic, or co-nondeterministic) of a group 
changes only linearly under change of finite generating set. 
So, combinatorial 
group theory gives us the following advantage over the ordinary formal
language formulation of computational complexity: algebraic invariance.
Complexity is a property of the group, no matter how the group arizes
(as words in a presentation, or as transformations of a space, or a set
with a composition operation). Note, however, that this invariance only
holds as long as we stick to finite generating sets.
 
It was also shown \cite{BiRedFunct} that every decision problem $L$  can 
be reduced (by a one-to-one linear-time reduction) to the word problem of 
some finitely generated group $G_L$, with the property that this word problem
has the same time complexity (up to a linear factor $n$) as the problem $L$.
(This was proved for deterministic and nondeterministic time complexity,
but the proof for the deterministic case also works for co-nondeterministic
time complexity.)  So, the word problem for finitely generated groups is as
general as decision problems overall, as far as time complexity is concerned.
Also, as a consequence, there exist finitely generated groups whose word 
problem is NP-complete, or coNP-complete.

The word problem of finitely presented groups is naturally related to 
nondeterministic time complexity; indeed, for a finite presentation 
$\langle A,R \rangle$, a word $w$ over $A^{\pm 1}$ is equivalent to 
$\varepsilon$ (the empty word) iff there exists a rewrite sequence consisting 
of applying relators in $R$; this rewrite process can be ``guessed'' and 
carried out by 
a nondeterministic Turing machine. More precisely, there is a close 
connection between the isoperimetric function and nondeterministic time 
complexity. In \cite{SBR} and \cite{BORS} it was shown that the word problem 
of a finitely generated group $G$ is in NP iff $G$ is embeddable in a 
finitely presented group whose isoperimetric function is polynomially bounded.
This implies that there exist finitely presented groups with
NP-complete word problem. 
The theorem extends to other nondeterministic time complexity classes.
A semigroup version of this result had been proved earlier \cite{Bi}.
It was also shown in \cite{SBR} that a function $f$ (with $f(n) \geq n^4$) 
is an isoperimetric function if $f(n)^4$ is the time complexity of a 
nondeterministic Turing machine and if $f$ is superadditive (i.e, 
$f(x+y) \geq f(x) + f(y)$). In particular, all functions 
$n^{\alpha}$ with $\alpha \geq 4$ ($\alpha \in \mathbb{Q}$) are 
isoperimetric functions.  Later, Brady and Bridson \cite{BradyBridson} also
proved that $n^{\alpha}$ is an isoperimetric function for all $\alpha$ 
ranging over a countable dense set of real numbers $\geq 2$. 
See also Section 3 of \cite{Bridson}.
On the other hand, there is no isoperimetric function between $n$ and $n^2$. 
More precisely, if an isoperimetric function $f$ satisfies $f(n) = o(n^2)$
then $f(n) = O(n)$; this is known as the Gromov gap. Groups with
linear isoperimetric function are called ``word hyperbolic'' (see 
\cite{Gromov}); they have shown up in many situations, and they have many 
special properties (e.g., their word problem can be decided in linear time 
by a deterministic  Turing machine). In summary, the study of connections 
between combinatorial group theory and nondeterministic time complexity 
has been successful, especially for combinatorial group theory, regarding 
isoperimetric functions.

\medskip

In this paper we look at connections between {\it co-nondeterministic} time 
complexity and combinatorial group theory. By definition, a problem 
$L$ (represented by a formal language) is in coNTime($T$) iff $L$ is accepted
by a co-nondeterministic Turing machine in time $T$. A co-nondeterministic
Turing machine is a Turing machine $M$ which is allowed to make choices
(just like a nondeterministic Turing machine), but which uses the following
acceptance rule: a word $w$ is accepted by $M$ iff {\it all} computation
paths of $M$ with input $w$ lead to an accept state. So, ``for-all'' is used
instead of nondeterminism's ``there-exists''. An equivalent definition is 
that coNTime($T$) consists of the languages whose complement is in
NTime($T$). Most books on computational complexity discuss coNTime($T$) and 
coNP; see e.g., \cite{Handb}. 
There are many famous coNP-complete decision problems,
that are as significant as the well-known NP-compete problems (although NP 
is much more popular than coNP in Computer Science). Here is a sampling: \\  
$\bullet$ The {\it tautology problem}: Given a boolean formula, is it a 
tautology, i.e., is it true for all truth-value assignments? Informal 
versions of this problem goes back to antiquity; the tautology problem is 
``the decision problem of boolean logic''.  \\  
$\bullet$ The {\it circuit equivalence problem}: Given two acyclic boolean
circuits (also called ``combinational circuits''), do they have the same
input-output function?  \\   
$\bullet$ {\it Integer linear programming equivalence problem}: Given two
instances of integer linear programming, do they have the same set of 
feasible solutions? \\   
$\bullet$ The {\it 4-coloring problem}: Given a planar graph, do we need four
colors to vertex-color it? (Note that every planar graph is 4-colorable, and 
the question whether a planar graph is 3-colorable is NP-complete.) \\  
$\bullet$ {\it Connectivity lower-bound}: Given a graph and an integer $k$, 
is the connectivity of the graph greater than $k$? Equivalently, does the 
graph remain connected when any $k$ edges are removed?    
  
Since there is a close connection between nondeterminism and finitely
presented groups, as we saw, and since NP is believed to be different from
coNP, one might expect at first that there is no natural connection between 
co-nondeterminism and combinatorial group theory.
However, if we take transformation groups as our starting point we see a 
hint at a connection: In a transformation group two elements $g_1$ and $g_2$ 
(permutations) are equal iff $g_1(x) = g_2(x)$ {\it for all} $x$ in the 
action space. Here again the for-all quantifier shows up, which
corresponds to co-nondeterminism. In order to investigate the complexity of 
problems about transformation groups, it is convenient to consider groups    
of transformations of words (i.e., strings over a finite alphabet). 
The groups introduced by Richard Thompson \cite{Th0} in the 1960s turn out to
be appropriate for this, not only based on their nice definition, but also
based on their history: they were used for constructing finitely presented 
groups with undecidable word problem \cite{McKTh}, and for proving a stronger 
form of the Boone-Higman theorem \cite{Th}. Below we give some background on
these groups. Note that here we do not view the Thompson groups as a special 
class of groups (as is usually done in the literature), but as a general 
formalism for describing all countable groups; in fact, all subgroups 
of ${\mathfrak S}_{\mathbb N}$ can be represented as Thompson groups 
(${\mathfrak S}_{\mathbb N}$ denotes the group of all permutations of the 
natural integers). 

In order to achieve coNP-hardness we show that every acyclic circuit can be 
``simulated'' by an element of a particular Thompson group (namely the 
finitely presented Thompson-Higman group $G_{3,1}$, defined below). 
So, we simulate a circuit by a permutation of strings over the 3-letter 
alphabet $\{0,1,\#\}$. The simulation is such 
that two circuits are equivalent iff their simulating permutations are  
equal when restricted to all strings that start with 0. 
Technically, the Thompson group elements are partial permutations of 
$\{0,1,\#\}^*$ that map certain maximal prefix codes bijectively to maximal 
prefix codes (see the background on Thompson groups below). 
This simulation is a polynomial-time many-to-one reduction from the circuit
equivalence problem (which is coNP-complete) to a problem about the
Thompson-Higman group $G_{3,1}$. In a succession of steps (see the more
detailed outline of the paper below), we reduce the latter problem to the 
word problem of another finitely presented Thompson group. We also reduce 
this problem to the word problem of a finitely generated {\it simple}
group (and we conjecture that this simple group is actually finitely 
presented). Moreover, we show that all the groups above have their word 
problem in coNP.

Our simulation of acyclic circuits by group elements is similar to the 
construction of a reversible circuit. This connects this paper with the 
classical topic of {\it reversible computation} (see \cite{Lec}, 
\cite{Ben73}, \cite{Ben89} for reversible Turing machines, and 
\cite{FredToff} for reversible acyclic circuits).  In our case the result
is stronger, since we do not just get reversibility but a finitely presented
group. 

\smallskip

Another motivation for this paper is a conjecture attributed to Higman
about a stronger form of the Boone-Higman theorem. The conjecture is that 
a finitely generated group $G$ has decidable word problem iff $G$ is 
embeddable into a finitely presented simple group.  
(It is well known that every finitely presented simple group has a decidable
word problem.)

A consequence of this conjecture would be that the word problem of finitely 
presented simple groups can have arbitrarily large time complexity. This 
means that for every function $T$ which is the time complexity of a 
deterministic Turing machine, there is a finitely presented simple group 
whose word problem cannot be decided in time $\leq T$.   
(Indeed, by \cite{BiRedFunct} finitely generated groups $G$ have arbitrarily
high complexity; moreover, a finitely generated subgroup $G$ of a group $S$
cannot have higher complexity than $S$, up to linear changes in the 
complexity function.)

On the other hand, all known finitely presented simple groups have word
problems with rather low complexity (in the cases where the complexity has 
been analyzed in detail it always turned out to be in the complexity class 
P). In that connection, see \cite{RoeverGR} and also \cite{GZ}, 
\cite{LiptZalc}. 
So, one might ask the opposite question: 
{\it Is there some cap on the computational complexity of 
the word problem of finitely presented simple groups? }
At the moment, neither Higman's conjecture nor the opposite question have 
much evidence in their favor (and, a priori, they could both be wrong).     
A contribution of this paper, in the direction of Higman's conjecture, is
the construction of a finitely generated simple group with coNP-complete
word problem; we conjecture that this group is also finitely presented.

\bigskip

\noindent {\bf Some background and notations on the Thompson groups}

\medskip

The Thompson groups, introduced by Richard Thompson in the 1960s 
\cite{Th0, Th}, provided the first known examples of simple finitely 
presented infinite groups.
Although Thompson defined his groups as permutation groups of certain sets 
of infinite words over the alphabet $\{ 0, 1\}$, we prefer the 
approach of E.~Scott \cite{ESc} and G.~Higman \cite{Hig74}, which enables 
us to define the Thompson groups as partial actions on the words over a 
finite alphabet. The advantage of finite words is that algorithmic 
problems and their complexity can be defined in a direct way. 

Let us introduce some terminology; we have made an effort to stay close to 
classical or widely used concepts. More details (and proofs) can be found 
in \cite{BiThomps}, and often also in \cite{ESc}, \cite{Hig74}, and 
\cite{Th}.  For a finite alphabet $A$, the set of all words over $A$ 
(including the empty word $\varepsilon$) is denoted by $A^*$. We will 
assume from now on that $A$ has a least two letters. Concatenation of two 
words $u, v \in A^*$ is denoted by $u \cdot v$ or $uv$; $A^*$ is a monoid 
under concatenation. For $X_1, X_2 \subseteq A^*$ the concatenation is  
$X_1 \cdot X_2 = X_1X_2 = \{x_1x_2 \in A^* : x_1 \in X_1, x_2 \in X_2\}$. 
A {\it right ideal} of $A^*$ is defined to be a subset $R \subseteq A^*$ such
that \,  $R \cdot A^* \subseteq R$ \, (i.e., $R$ is closed under
concatenation by any word in $A^*$ on the right).
For two words $u, v \in A^*$, we say that $u$ is a {\it prefix} of $v$ iff
$v = ux$ for some $x \in A^*$; we also write $u \geq_{\rm pref} v$ or
$v \leq_{\rm pref} u$; this is a partial order, related to set inclusion by
the fact that $v \leq_{\rm pref} u$ iff $vA^* \subseteq uA^*$.
We say that $u$ and $v$ are {\it prefix-comparable} iff 
$v \leq_{\rm pref} u$ or $u \leq_{\rm pref} v$;
we denote this by $u \lesseqgtr_{\rm pref} v$.
A {\it prefix code} over $A$ is defined to be a subset $C$ of $A^*$
such that no element of $C$ is a strict prefix of another element of $C$.
A {\it maximal prefix code} over an alphabet $A$ is a prefix
code over $A$ which is not a strict subset of any other prefix code over $A$.
For a right ideal $R$ of $A^*$, a set $\Gamma \subseteq R$ is called a set of
{\it right-ideal generators} of $R$ iff \ $R = \Gamma \cdot A^*$.
One can prove  that any right ideal $R$ of $A^*$ has a unique minimal (under 
inclusion) set of right-ideal generators, and this set of generators is a 
prefix code.
Right ideals of $A^*$ and prefix codes over $A$ are in one-to-one 
correspondence. A right ideal $R$ of $A^*$ is said to be {\it finitely 
generated} iff the prefix code corresponding to $R$ is finite.
A right ideal $R$ of $A^*$ is called {\it essential} iff $R$ has a non-empty
intersection with every right ideal of $A^*$. One can prove that a right 
ideal is essential iff its prefix code is a maximal prefix code.
 
A {\it right-ideal homomorphism} of $A^*$ is defined to be a function
$\varphi: R_1 \to R_2$ such that $R_1$ and $R_2$ are right ideals of $A^*$,
and such that for all $u \in R_1$ and all $x \in A^*$: \   
$\varphi(u) \cdot x = \varphi(ux)$.
A {\it right-ideal isomorphism} of $A^*$ is a  bijective right-ideal
homomorphism.
The set of all right-ideal homomorphisms (or isomorphisms) of $A^*$ is in 
one-to-one correspondence with the set of all functions (respectively  
bijections) between prefix codes of $A^*$.  For a right-ideal isomorphism 
$\varphi: P_1A^* \to P_2A^*$, where $P_1$ and $P_2$ are prefix codes, the 
restriction \ $\tau_{\varphi}: P_1 \to P_2$ \ is a bijection, and 
$\tau_{\varphi}$ determines $\varphi$ uniquely. Following Thompson, the 
restriction $\tau_{\varphi} : P_1 \to P_2$ of $\varphi$ will be called the
{\it table of} \ $\varphi$, and will be used to represent $\varphi$ by a 
traditional function table. (In \cite{Hig74} and \cite{ESc} this 
was called the ``symbol of $\varphi$''.)
The maximal prefix code $P_1$ is called the {\it domain code} of $\varphi$,
and $P_2$ is called the {\it image code} or {\it range code} of $\varphi$.
An {\it extension} of a right-ideal isomorphism
$\varphi: R_1 \to R_2$ is defined to be a right-ideal isomorphism
$\Phi: J_1 \to J_2$ where $J_1, J_2$ are right ideals such that
$R_1 \subseteq J_1, \ R_2 \subseteq J_2$, and $\Phi$ agrees with $\varphi$
on $R_1$ (i.e., $\Phi(x) = \varphi(x)$ for all $x \in R_1$). In that case we
also call $\varphi$ a {\it restriction} of $\Phi$.
A right-ideal isomorphism is said to be {\it maximal} iff it has no strict
extension in $A^*$; it is called {\it extendable} otherwise.
We denote the maximum extension of $\varphi$ by ${\sf max} \, \varphi$; \ 
one can prove (see \cite{ESc} or \cite{BiThomps}) that the maximum 
extension of an isomorphism between {\it essential} right ideals is unique.

The above concepts can be pictured using trees. The monoid $A^*$ can be 
described by the Cayley graph of the right regular representation of $A^*$
relative to the generating set $A$. We will simply call this {\it the tree 
of} $A^*$.  It is an infinite tree rooted at the empty word $\varepsilon$.
Every vertex has $|A|$ children. Every subset of $A^*$  is pictured as a 
set of vertices of this infinite tree.
A prefix code is pictured as a set of vertices, no two of which lie on a 
same directed path from the root.  For any prefix code $P \subset A^*$ 
($P \neq \emptyset$), the {\it prefix tree} of $P$ is defined to be the 
subtree of the tree of $A^*$, whose vertex subset consists of all the 
prefixes of words in $P$ (and whose root is still $\varepsilon$).
Hence, the set of leaves of this subtree is $P$.

One can prove (see \cite{ESc} or \cite{BiThomps}) that an isomorphism of 
finitely generated essential right ideals 
$\varphi: P_1 A^* \to P_2 A^*$, with $P_1$ and $P_2$ finite maximal prefix 
codes, is {\em extendable} iff there are $x_0, y_0 \in A^*$ such that
for every letter $\alpha \in A$: \   
$x_0 \alpha \in P_1$, $y_0 \alpha \in P_2$, and 
$\varphi(x_0 \alpha) = y_0 \alpha$. (If this condition holds, $\varphi$ can 
be extended by mapping $x_0$ to $y_0$.) 
More generally (see \cite{BiThomps}), an isomorphism of (not necessarily 
finitely generated) essential right ideals $\varphi: P_1 A^* \to P_2 A^*$, 
with $P_1$ and $P_2$ arbitrary maximal prefix codes, is {\em extendable} iff 
there are $x_0, y_0 \in A^*$ and there exists a maximal prefix code 
$Q \subseteq A^*$ with $|Q| > 1$ such that for all $q \in Q:$ \ \
$x_0q \in P_1$, \ $y_0q \in P_2$, and $\varphi(x_0q) = y_0q$.

We now define the Thompson groups, following the approach of Scott 
\cite{ESc} and Higman \cite{Hig74}. The tree representation of codes 
connects this definition and the definition by action on finite trees used 
in \cite{CFP}.  
The {\bf Thompson-Higman group} $G_{N,1}$ is the partial action group on
$A^*$ (for some fixed alphabet $A$ with $|A| = N$), consisting of all 
maximal isomorphisms between {\it finitely generated} essential right ideals 
of $A^*$.
The Thompson-Higman group ${\mathcal G}_{N,1}$ is the partial action group 
on $A^*$ consisting of {\it all} maximal isomorphisms between essential right 
ideals of $A^*$. 
Multiplication in ${\mathcal G}_{N,1}$, and hence in the subgroup $G_{N,1}$
and in any subgroup of ${\mathcal G}_{N,1}$, is defined as follows: For 
$\varphi, \psi \in {\mathcal G}_{N,1}$ the product
$\varphi \cdot \psi$ is \ {\sf max}$(\varphi \circ \psi)$ \ (i.e., the 
maximum extension of the composition of $\psi$ and $\varphi$, where $\psi$ 
is applied first). 
In general, in this paper, we apply (partial) functions 
on the left of the argument, and hence compose functions from right to left. 

In this paper we call any partial transformation subgroup of 
${\mathcal G}_{N,1}$ (for any integer $N \geq 2$) a {\it Thompson group}.
(This is a slight misnomer, since these groups are actually more than just 
groups; they are partial transformation groups.)
It is easy to see that every countable group is isomorphic to a Thompson 
group; in fact (see e.g. \cite{BiThomps}), every subgroup of 
${\mathfrak S}_{\mathbb Z}$ (the group of all permutations of the integers)
can be represented as a Thompson group.
It is remarkable that Thompson groups consist of partial transformations;
it is the uniqueness of the maximal extension that enables them, nevertheless,
to be groups.

\bigskip

\noindent {\large \bf Main results}

\medskip

In this paper we use {\it polynomial-time constant-arity conjunctive 
reduction} (instead of many-to-one reduction). This is defined in 
Definition \ref{conj_red_def}. 
The complexity classes P, NP, coNP, as well as most other common complexity 
classes containing P, are closed under this reduction.  

\begin{thm} \label{mainThm1} \ 
There exists a finitely presented group $G$ whose word problem is
coNP-complete (with respect to polynomial-time constant-arity conjunctive
reduction).

\smallskip

\noindent Moreover, we have: \\  
$\bullet$ \  The group $G$ is explicitly embedded into ${\mathcal G}_{3,1}$ 
as \ $G = $
$ \langle G_{3,1}^{\rm mod \, 3}(0,1;\#) \cup \{\kappa_{321}\} \rangle$
 \  (see Theorem \ref{finPrescoNPcompl}).  The subgroup 
$G_{3,1}^{\rm mod \, 3}(0,1;\#)$ of $G_{3,1}$ is defined in Definition
\ref{notation_G01} and at the end of Step 1 below, and is finitely presented. 
The element  
$\kappa_{321} = \kappa_3\kappa_2\kappa_1 \in {\mathcal G}_{3,1}$ is defined 
in Section 2. \\  
$\bullet$ \ $G$ is an HNN extension (by one stable letter) of 
$G_{3,1}^{\rm mod \, 3}(0,1;\#)$. Moreover, $G$ is isomorphic to a 
semidirect product \ $G_{3,1}^{\rm mod \, 3}(0,1;\#) \rtimes {\mathbb Z}$.
\end{thm}

\begin{thm} \label{mainThm2} \
There exists a finitely generated {\em simple} group $S$ whose word problem is
coNP-complete (with respect to polynomial-time constant-arity conjunctive
reduction).
\end{thm}
The group $S$ is explicitly embedded into ${\mathcal G}_{3,1}$ as \  
$S = \langle G_{3,1} \cup \{\kappa_0, \kappa_1, \kappa_2 \}\rangle'$, i.e., 
the commutator subgroup of
$\langle G_{3,1} \cup \{\kappa_0, \kappa_1, \kappa_2 \}\rangle$ \    
(Theorem \ref{finGenSimple}), where $\kappa_0$, $\kappa_1$, and $\kappa_2$ 
are elements of ${\mathcal G}_{3,1}$ defined in Section 2. Moreover, $S$ 
has finite index in 
$\langle G_{3,1} \cup \{\kappa_0, \kappa_1, \kappa_2 \}\rangle$.

We conjecture that 
$\langle G_{3,1} \cup \{\kappa_0, \kappa_1, \kappa_2 \}\rangle$, and
hence $S$, is not only finitely generated but also finitely
presented.  This would give us a finitely presented simple group with 
coNP-complete word problem.   

\bigskip

\noindent {\large \bf Overview of the paper}

\medskip

\noindent $\bullet$ {\bf Step 1} (Sections 2 and 3):  

Recall that $G_{3,1}$ is the Thompson-Higman group of right-ideal 
isomorphisms between finitely generated essential right ideals of the free 
monoid $\{0,1,\#\}^*$. It is well known that $G_{3,1}$ is finitely 
presented \cite{Hig74}; let $\Delta_{3,1}$ be a finite generating set
of $G_{3,1}$.
We give a polynomial-time many-to-one reduction of the circuit equivalence 
problem to the following ``{\it word problem with restriction}'' in the 
Thompson-Higman group $G_{3,1}$: 

\smallskip

\noindent {\sc Input}: Two words $u,v$ over 
$\Delta_{3,1}^{\pm 1} \cup \{ \tau_{i,i+1} : i \geq 0 \}$, where 
$\tau_{i,i+1}$ is the element of $G_{3,1}$ that transposes the bits in 
positions $i$ and $i+1$ in any string 
$x_0 x_1 \ldots x_i x_{i+1} \ldots \# \ \in \ \{0,1\}^*\#$. \\ 
{\sc Question}: Are the two elements of $G_{3,1}$, represented by $u,v$, 
equal when restricted to the subset 
$0 \, \{0,1\}^* \, \#$ \ of \ $\{0,1,\#\}^*$ \ ?

\smallskip

\noindent
In order to find the above reduction, we first represent the circuit 
components by elements of $G_{3,1}$: {\sc and}, {\sc or}, {\sc not}, as well 
as wire forking (i.e., duplication or copying of variables), and wire crossing 
(i.e., permutations of variables); wire crossings are described by the 
transpositions $\tau_{i,i+1}$.

Now let $C$ be any acyclic boolean circuit, with input-output function
$f_C: \{0,1\}^m \to \{0,1\}^n$. We simulate $C$ by a Thompson group element 
$\Phi_C \in G_{3,1}$ such that: \\  
- the action of $\Phi_C$ on the subset $0 \{0,1\}^* \#$ represents the 
function $f_C$ in the sense that for all $x_0, x_1, \ldots, x_m \in \{0,1\}$
and all $w \in \{0,1\}^*$:  

 \ \ \ \ \   
$\Phi_C(0 x_1 \ldots x_m \, w \, \#) = 0 x_1 \ldots x_m \, f_C(x_1, \ldots,
x_m) \, w \, \#$; \\ 
- the word-length of $\Phi_C$ over 
$\Delta_{3,1}^{\pm 1} \cup \{ \tau_{i,i+1} : i \geq 0 \}$, as well as the 
largest subscript of the $\tau_{i,i+1}$ used to represent $\Phi_C$, 
have a polynomial upper bound in terms of the circuit size $|C|$; \\  
- a word $w_C$ over $\Delta_{3,1}^{\pm 1} \cup \{ \tau_{i,i+1} : i \geq 0 \}$, 
representing $\Phi_C$, can be computed deterministically in polynomial time 
(in terms of $|C|$.

\smallskip

Note that although $G_{3,1}$ is finitely generated, we are using an infinite
generating set here in order to obtain the word $w_C$ with polynomial length;
in fact, $\tau_{i,i+1}$ has exponential word-length over $\Delta_{3,1}$.
Eventually we will want a finitely generated (and finitely presented) group
for representing $C$. For this we introduce elements 
$\kappa_i \in {\mathcal G}_{3,1}$ $i =0,1,2,3$  such that each $\tau_{i,i+1}$
has polynomial word length over $\Delta_{3,1} \cup \{\kappa_0, \kappa_1,
\kappa_2, \kappa_3 \}$. However, $\kappa_i$ does not belong to $G_{3,1}$.

\smallskip

We observe that $\Phi_C$ and the representatives of the circuit elements
belong to the following subgroup of $G_{3,1}$: 

\smallskip

\noindent $G_{3,1}^{\rm mod \, 3}(0,1;\#) \ = \ $

$ \ \ \  \{ \phi \in G_{3,1} : \ $ 
$\phi$ and $\phi^{-1}$ map $\{0,1\}^*$ \ to \ $\{0,1\}^*$ and 
map $\{0,1\}^*\#$ \ to \ $\{0,1\}^*\#$; 

\hspace{1in} $\phi$ and $\phi^{-1}$ are defined everywhere on $\{0,1\}^*\#$;

\hspace{1in} moreover, for all $x \in \{0,1\}^*$, \ 
    $|\phi(x)| \equiv |x|$ mod 3 when $\phi(x)$ is defined$\}$.

\smallskip

\noindent From now on we will usually use $G_{3,1}^{\rm mod \, 3}(0,1;\#)$, 
rather than $G_{3,1}$. Later we will prove that 
$G_{3,1}^{\rm mod \, 3}(0,1;\#)$ is finitely presented.
 
\medskip
 
\noindent $\bullet$ {\bf Step 2} (Section 4):  

It follows from step 1 that two circuits $C_1, C_2$ are equivalent iff 
$\Phi_{C_2}^{-1} \Phi_{C_1}$ fixes every point in $0 \{0,1\}^* \#$ \, on which
$\Phi_{C_2}^{-1} \Phi_{C_1}$ is defined. Thus, we have reduced the circuit 
equivalence problem to the generalized word problem of the subgroup 
pFix$(0 \{0,1\}^* \#)$ of $G_{3,1}^{\rm mod \, 3}(0,1;\#)$. Here, for any 
$S \subseteq \{0,1,\#\}^*$, pFix$(S)$ denotes the ``partial fixator''

\smallskip

pFix$(S) \ = \ \{ \phi \in G_{3,1}^{\rm mod \, 3}(0,1;\#) : \, \phi$ fixes 
all points of $S$ on which $\phi$ is defined$\}$.  

\smallskip

\noindent In these problems we still represent words over the infinite 
generating set $\Delta_{3,1}^{\pm 1} \cup \{ \tau_{i,i+1} : i \geq 0 \}$ for 
$G_{3,1}$.

\medskip

\noindent $\bullet$ {\bf Step 3} (Section 5): \\  
We show that for any $g \in G_{3,1}^{\rm mod \, 3}(0,1;\#)$:

\smallskip

$g \in$ pFix$(0 \{0,1\}^* \#)$ \ \ iff \ \ $gh = hg$ \ for all \ $h \in$
pFix$(\{1,\#\} \{0,1\}^* \#)$.

\medskip

\noindent $\bullet$ {\bf Step 4} (Section 6): \\
Moreover, pFix$(\{1,\#\} \{0,1\}^* \#)$ is finitely generated (and 
in fact finitely presented). The above commutation relation only needs to be 
checked between $g$ and the finitely many generators of 
pFix$(\{1,\#\} \{0,1\}^* \#)$. 
The group $G_{3,1}^{\rm mod \, 3}(0,1;\#)$ is also finitely presented.

As a consequence of steps 3 and 4, we have reduced the circuit equivalence 
problem to the word problem of $G_{3,1}^{\rm mod \, 3}(0,1;\#)$ (and hence 
also of $G_{3,1}$), via a polynomial-time constant-arity conjunctive reduction
(the arity being the number of generators of pFix$(\{1,\#\} \{0,1\}^* \#)$). 
The generating set of $G_{3,1}^{\rm mod \, 3}(0,1;\#)$ used for these problems 
is still the infinite set $\Delta \cup \{ \tau_{i,i+1} : i \geq 0 \}$, where
$\Delta$ is any finite generating set of $G_{3,1}^{\rm mod \, 3}(0,1;\#)$. 

\medskip

\noindent $\bullet$ {\bf Step 5} (Section 7): \\
We show that conjugation by $\kappa_i$ ($i = 0,1,2,3$) is an automorphism of 
$G_{3,1}^{\rm mod \, 3}(0,1;\#)$. (It is for this property that we needed the 
length-preservation mod 3 in the elements of $G_{3,1}^{\rm mod \, 3}(0,1;\#)$.)
Hence, the following HNN extension yields a group $H(0,1;\#)$ which contains 
$G_{3,1}^{\rm mod \, 3}(0,1;\#)$ and $\kappa_3\kappa_2\kappa_1$. 

\smallskip

$H(0,1;\#) \ = \ $
$\langle G_{3,1}^{\rm mod \, 3}(0,1;\#) \cup \{t\} \ : \ $
  $\{t \, g \, t^{-1} = g^{\kappa_3\kappa_2\kappa_1} : g \in $
  $G_{3,1}^{\rm mod \, 3}(0,1;\#) \} \rangle$.

\smallskip

\noindent Since in step 4 we saw that $G_{3,1}^{\rm mod \, 3}(0,1;\#)$ is 
finitely presented, $H(0,1;\#)$ is finitely presented.
The transpositions $\tau_{i,i+1}$ have linear word-length over the 
finite generating set of $H(0,1;\#)$. Hence, the circuit equivalence problem
reduces (via polynomial-time constant-arity conjunctive reduction) to the
word problem of the finitely presented group $H(0,1;\#)$ (over its finite
generating set). 

The group $H(0,1;\#)$ is isomorphic to the subgroup \ 
$\langle G_{3,1}^{\rm mod \, 3}(0,1;\#) \ \cup $
$\{\kappa_3\kappa_2\kappa_1\} \rangle$
 \ of ${\mathcal G}_{3,1}$, and also to the semidirect product 
$G_{3,1}^{\rm mod \, 3}(0,1;\#) \rtimes {\mathbb Z}$.

\medskip

\noindent $\bullet$ {\bf Step 6} (Section 8): \\  
We prove that the word problems of $H(0,1;\#)$ and, more generally, of 
$\langle G_{3,1} \cup \{ \kappa_0, \kappa_1, \kappa_2 \} \rangle$ are in coNP. 
As a consequence, the finitely presented group $H(0,1;\#)$ has a 
coNP-complete word problem (relative to polynomial-time constant-arity
conjunctive reduction); this is the group $G$ of Theorem \ref{mainThm1}. 

By results of Thompson and Scott, the commutator subgroup \ 
$\langle G_{3,1} \cup \{ \kappa_0, \kappa_1, \kappa_2 \} \rangle'$ \ is a 
simple group.  We prove that \   
$\langle G_{3,1} \cup \{ \kappa_0, \kappa_1, \kappa_2 \} \rangle'$ \ has 
finite index in 
$\langle G_{3,1} \cup \{ \kappa_0, \kappa_1, \kappa_2 \} \rangle$. 
Hence, $\langle G_{3,1} \cup \{ \kappa_0, \kappa_1, \kappa_2 \} \rangle'$ 
is a finitely generated simple group with coNP-complete word problem; this
is the group $S$ of Theorem \ref{mainThm2}.
Moreover, if $\langle G_{3,1} \cup \{\kappa_0, \kappa_1, \kappa_2\} \rangle$
is finitely presented (as we conjecture), 
$\langle G_{3,1} \cup \{ \kappa_0, \kappa_1, \kappa_2 \} \rangle'$ will also 
be finitely presented.  

\medskip

\noindent $\bullet$ {\bf Appendix} (Section 9): \\  
The first subsection of the Appendix contains 
properties of prefix codes, used in the paper.

Another  subsection of the Appendix shows that 
Theorems \ref{mainThm1} and \ref{mainThm2} and the Overview above also hold 
with $G_{3,1}^{\rm mod \, 3}(0,1;\#)$ replaced by another subgroup of 
$G_{3,1}$, namely by

\smallskip 

$G_{3,1}^{\rm mod \, 3}(0,1) \ = \  $
$ \{ \phi \in G_{3,1} : \ $
$\phi$ and $\phi^{-1}$ map $\{0,1\}^*$ \ to \ $\{0,1\}^*$ and

\hspace{1.2in} for all $x \in \{0,1\}^*$, \
    $|\phi(x)| \equiv |x|$ mod 3 when $\phi(x)$ is defined$\}$.

\smallskip 

\noindent The proofs for $G_{3,1}^{\rm mod \, 3}(0,1)$ are similar to (but 
somewhat more complicated than) the proofs for 
$G_{3,1}^{\rm mod \, 3}(0,1;\#)$, and appear in the Appendix.

A special property is shown: $G_{3,1}^{\rm mod \, 3}(0,1)$ is the largest
subgroup of $G_{3,1}$ closed under conjugation by 
$\kappa_3 \kappa_2 \kappa_1$.

%%%%%%%%%%%%%%%%%%%%%%%%%%%%%%%%%%%%%%%%%%%%%%%%%%%%%%%%
% Section 2 
%%%%%%%%%%%%%%%%%%%%%%%%%%%%%%%%%%%%%%%%%%%%%%%%%%%%%%%%

\section{Circuits and permutations of boolean variables}

Acyclic boolean circuits are a fundamental model of computation
\cite{Wegener}, \cite{Handb}, \cite{JESavage}. The {\it equivalence problem}
for acyclic boolean circuits, mentioned above, is a well-known example of 
a coNP-complete problem.

Circuits are traditionally built from the boolean functions {\sc and},
{\sc or}, and {\sc not}, with domains $\{0,1\}^2$ or $\{0,1\}$, and image
$\{0,1\}$. Moreover, circuits use the {\it fork} (or ``fan-out'', or
``duplication'') function \
{\sc fork}$: x \in \{0,1\} \mapsto (x,x) \in \{0,1\}^2$. The use of
{\sc fork} is usually tacit; in a circuit diagram, {\sc fork} appears whenever
a wire fans out (or forks, or splits) to become two wires that carry the same
boolean value. One can view an acyclic boolean circuit as a composition of
several copies of the functions {\sc and}, {\sc or}, {\sc not}, and
{\sc fork}. Since {\sc and}, {\sc or}, and {\sc fork} are multi-variable
functions, composition is complicated and requires a circuit diagram
(which is essentially an acyclic graph) to describe how the operations
are connected.  We will use {\sc and} and {\sc or} gates with
fan-in 2 only.

We will see that a circuit can be represented by {\it ordinary composition of
functions}, thanks to Thompson groups. We have seen that these groups can be 
described as partial action groups, acting on strings.
We will use this partial action to simulate circuits.

\medskip

The functions {\sc and}, {\sc or}, {\sc not}, and {\sc fork}
that make up acyclic circuits use one or two boolean variables (that range
over the set of boolean values $\{0,1\}$). An acyclic circuit has boolean
variables $(x_0, x_1, \ldots, x_{m-1})$ as input (ranging over all of
$\{0,1\}^m$), and boolean variables $(y_0, y_1, \ldots, y_{n-1})$ as output
(ranging over a subset of $\{0,1\}^n$); the circuit computes a function
$f: \{0,1\}^m \to \{0,1\}^n$. We will extend the function $f$ to the partial
function

\smallskip

 \ \ \ \ \ $x_0 x_1 \ldots x_{m-1} w \in \{0,1\}^* \ \longmapsto \ $
 $f(x_0, x_1, \ldots, x_{m-1}) \ w \in \{0,1\}^*$

\smallskip

\noindent for any $w \in \{0,1\}^*$. We let boolean functions operate on an
arbitrary (large enough) number of variables, rather than a fixed number.

In this paper, we always write functions on the left of their argument.
Also, we make the following convention: Let $\phi: A^* \to A^*$ be a partial 
map and $x \in A^*$; when we write $\phi(x)$ it is to be understood that
$\phi(x)$ is defined (i.e., $x \in$ Dom$(\phi)$). 

\medskip

Since we write the variables \ $x_0, x_1, \ldots, x_{m-1}, \ldots$ \  
in a fixed order we need to introduce maps that permute these variables. In 
circuit drawings this corresponds to {\it crossing of wires}. 
In particular, we use the {\it transposition of variables} $x_i$, $x_j$ 
(with $0 \leq i < j$), defined by \ $ux_ivx_jw \in \{0,1\}^* \longmapsto $
$ux_jvx_iw \in \{0,1\}^*$, where $|u| = i, |v| = j - i -1, w \in \{0,1\}^*$. 

The finite symmetric groups are generated by two elements, a 
transposition and a cyclic permutation. Here we also want to obtain a 
finite number of generators, but since we deal now with unbounded finite 
bit-strings, we need to consider new versions of the cyclic permutation.
This, in turn, requires the introduction of a new letter into the alphabet; 
the new letter, denoted \#, will act as a ``boundary marker'' for the cyclic 
permutations. A first idea of an unbounded cyclic permutation would be to take 
 \  $x_0 x_1 \ldots x_{m-1} \# w$ $\longmapsto$ $x_1 \ldots x_{m-1} x_0\# w$, 
 \ for all $x_0, x_1, \ldots, x_{m-1} \in \{0,1\}$, and 
$w \in \{0,1,\#\}^*$; but it turns out that this definition does not lead to 
good properties (some Thompson groups that we will work with are not closed 
under conjugation by this permutation). So we will use the following 
permutations of ${\mathbb N}$, written as infinite products of disjoint cyclic 
permutations. Recall that a cycle $(i|j|k)$ (for three distinct elements 
$i, j, k \in \mathbb N$), denotes the 
permutation $i \mapsto j \mapsto k \mapsto i$, and $x \mapsto x$ for 
$x \not\in \{i,j,k\}$. We denote the group of all permutations of 
${\mathbb N}$  by ${\mathfrak S}_{\mathbb N}$. 
Again, recall that we write maps to the left of the argument.

\bigskip

 \ \ \ \ \
$\gamma_0 \ = \ \ldots \ \ldots \ (3n \ | \ 3n+1 \ | \ 3n+2) \ $
$\ldots \ (3 \ | \ 4 \ | \ 5) \ (0 \ | \ 1 \ | \ 2)$,

\medskip

 \ \ \ \ \
$\gamma_1 \ = \ \ldots \ \ldots \ (3n+1 \ | \ 3n+2 \ | \ 3(n+1) ) \ $
$\ldots \ (4 \ | \ 5 \ | \ 6) \ (1 \ | \ 2 \ | \ 3) \ (0)$,

\medskip

 \ \ \ \ \
$\gamma_2 \ = \ \ldots \ \ldots \ (3n+2 \ | \ 3(n+1) \ | \ 3(n+1)+1 ) \ $
$\ldots \ (5 \ | \ 6 \ | \ 7) \ (2 \ | \ 3 \ | \ 4) \ (1) \ (0)$,

\medskip

 \ \ \ \ \
$\gamma_3 \ = \ \ldots \ \ldots \ (3n \ | \ 3n+1 \ | \ 3n+2) \ $
$\ldots \ (3 \ | \ 4 \ | \ 5) \ (2) \ (1) \ (0)$.

\bigskip

\noindent Based on these permutations of ${\mathbb N}$ we define the 
following elements $\kappa_0$, $\kappa_1$, $\kappa_2$, $\kappa_3$ 
$\in {\mathcal G}_{3,1}$. 

The effect of $\kappa_i$ ($i = 0,1,2,3$) on a string 
$x_0 x_1 \ldots x_m \, \# \, w$ (with $x_0, x_1, \ldots, x_m$
$\in \{0,1\}$, $w \in \{0,1,\#\}^*$) is to permute the bits 
$x_0 x_1 \ldots x_m$ according to $\gamma_i$; the bit $x_k$ at position $k$ 
($0 \leq k \leq m$) is moved to position $\gamma_i(k)$. Thus, \   
$\kappa_i(x_0 x_1 \ldots x_m\, \#\, w) \ = \ y_0 y_1\ldots y_m\, \# \, w$, 
where $y_{\gamma_i(k)} = x_k$. Equivalently, \  
$y_j = x_{\gamma_i^{-1}(j)}$ \ for $0 \leq j \leq m$.    
According to this definition, $\kappa_i$ is well defined on a string
$x_0 x_1 \ldots x_m \, \# \, w$ when $m \equiv i$ mod 3. To make $\kappa_i$
well defined on all strings in $\{0,1,\#\}^*$ we let $\kappa_i$ act as the
identity on the one or two right-most ``extra bits'', when $m$ is not 
$\equiv i$ mod 3. The detailed definition of $\kappa_i$ is as follows:

\medskip

\noindent $\bullet$ 
For \ $x \ = \ x_0 \ldots x_i \ldots x_{3n+2} \, r \, \#$, \  
where $n \in {\mathbb N}$, $x_i \in \{0,1\}$ ($0 \leq i \leq 3n+2$), and
$r \in \{0,1\}^{\leq 2}$, we define

\medskip

 \ \ \ \ \
$\kappa_0(x) \ = \ $
   $x_{\gamma_0^{-1}(0)} \ldots x_{\gamma_0^{-1}(i)} \ldots $
   $x_{\gamma_0^{-1}(3n+2)}  \, r \, \#$, \ \ \ and

\medskip

 \ \ \ \ \
$\kappa_3(x) \ = \ $
   $x_{\gamma_3^{-1}(0)} \ldots x_{\gamma_3^{-1}(i)} \ldots$
   $ x_{\gamma_3^{-1}(3n+2)} \, r \, \#$.

\medskip

\noindent $\bullet$ Similarly, for \
$x \ = \ x_0 \ldots x_i \ldots x_{3(n+1)} \, r \, \#$ \ we define

\medskip

 \ \ \ \ \
$\kappa_1(x) \ = \ $
 $x_{\gamma_1^{-1}(0)} \ldots x_{\gamma_1^{-1}(i)} \ldots $
  $x_{\gamma_1^{-1}(3(n+1))}  \, r \, \#$.

\medskip

\noindent $\bullet$ For \
$x \ = \ x_0 \ldots x_i \ldots x_{3(n+1)+1} \, r \, \#$ \ we define

\medskip

 \ \ \ \ \
$\kappa_2(x) \ = \ $
 $x_{\gamma_2^{-1}(0)} \ldots x_{\gamma_2^{-1}(i)} \ldots $
  $x_{\gamma_2^{-1}(3(n+1)+1)}  \, r \, \#$.

\bigskip

\noindent We will abbreviate \ $\kappa_3 \kappa_2 \kappa_1(\cdot)$ \ to 
$\kappa_{321}(\cdot)$. The element $\kappa_{321} \in {\mathcal G}_{3,1}$ 
will play an important role in this paper.

\bigskip

The introduction of the new letter \# in the boolean alphabet $\{0,1\}$ 
forces us to rethink the correspondence between the Thompson groups. 
We will now use the Thompson-Higman group $G_{3,1}$ of \cite{Hig74}, acting
on $\{0,1,\#\}^*$. The Thompson-Higman group  $G_{3,1}$ is isomorphic 
to a subgroup of the Thompson group $V$. 

\medskip

As a Thompson group element, the {\bf transposition} 
$\tau_{i,j} \in G_{3,1}$ of \ $x_i, x_j$ \ ($0 \leq i < j$) is defined as 
follows.
The domain and image prefix code of $\tau_{i,j}$ is the finite maximal 
prefix code

\medskip

${\rm domC}(\tau_{i,j}) = {\rm imC}(\tau_{i,j}) = $
$\ \{0,1\}^{j+1} \ \cup \ \{0,1\}^{\leq j} \, \#$ .

\medskip

\noindent  On an argument in $\{0,1\}^{j+1}$ (i.e., the number of ``boolean
variables'' in the argument is at least $j+1$) we define
   
\medskip 

 \ \ \ $\tau_{i,j}: \ u x_i v x_j \ \longmapsto \ u x_j v x_i$ 

\medskip 

\noindent where \ $x_i, x_j \in \{0,1\}$, \ $u \in \{0,1\}^i$, and
$v \in \{0,1\}^{j-i-1}$. 

\medskip

We also need to consider the case of an argument of the form $z \#$  where
$z = x_0x_1 \ldots x_{\ell-1} \in \{0,1\}^{\ell}$ with $\ell \leq j$. 
Here, the number of boolean variables in the argument is strictly less than 
$j+1$; in other words, the argument is ``too short'' for the transposition 
$\tau_{i,j}$. 
For those arguments we define $\tau_{i,j}$ in such a way that  

\smallskip

\noindent  $\bullet$ \ $\tau_{i,j}$ is be a permutation of the boolean 
variables $x_0$, $x_1$, $\ldots$, $x_{\ell-1}$; 

\smallskip

\noindent  $\bullet$ \ when $\ell = 0$, \ $\tau_{i,j}(\#) = \#$.

\smallskip

\noindent $\bullet$ \ when $0 < i < j$, $\tau_{i,j}$ fixes $x_0$, i.e.,
$\tau_{i,j}$ maps the set \ $0 \{0,1\}^* \ \cup \ 0 \{0,1\}^* \#$  
 \ into itself, and it maps \ $1\{0,1\}^* \ \cup \ 1\{0,1\}^* \#$  
 \ into itself.

\smallskip

\noindent The actual details of the definition when the argument is too short 
are a matter of convenience, and will be given later. 
However, we will completely define $\tau_{0,1}$ here, by letting it act as
the identity map on \ $\{0,1\}^{\leq 1} \#$; and of course,
$x_0x_1 \mapsto x_1x_0$ for all $x_0, x_1 \in \{0,1\}$. 
Similarly, we completely define $\tau_{1,2}$ by letting it act as
the identity map on \ $\{0,1\}^{\leq 2} \#$; and 
$x_0x_1x_2 \mapsto x_0x_2x_1$ for all $x_0, x_1, x_2 \in \{0,1\}$.
For all $i, j$ we define $\tau_{i,j}$ to mean the same thing as $\tau_{j,i}$.

\medskip

The classical formulas about transpositions are still true for this 
definition of transpositions. 
For all $i, j, k \geq 0$, and for all $x \in \{0,1\}^*$: 

\smallskip

 \ \ \ \ \ $\tau_{i,j}(x) \ = \ \tau_{i,k} \ \tau_{k,j}(x)$, \ \ \  
 if $|x| > {\rm max}\{i,j,k\}$

\smallskip

 \ \ \ \ \ $\tau_{i,j}(x) \ = \ $
 $ \tau_{i,i+1} \ \tau_{i+1,i+2} \ \ldots \ \tau_{j-2,j-1} \ \tau_{j-1,j}$
 $ \ \tau_{j-2,j-1} \ \ldots \ \tau_{i+1,i+2} \ \tau_{i,i+1}(x)$ \ \ \ 
 when $0 \leq i < j$, and $|x| > j$. 

\smallskip

For an argument $x \in \{0,1\}^* \#$  that is ``too short'',
we will simply {\em define} $\tau_{i,j}$ by the second of the above formulas. 
Recall that initially we picked $\tau_{i,j}(x\#)$ to be arbitrary (subject to 
the requirement that $\tau_{i,j}$ should be a permutation of its domain code,
and that $\tau_{i,j}$ should fix the left-most boolean variable when $0<i$). 
Now $\tau_{i,i+1}(x\#)$ is still arbitrary (for all $0 \leq i$, 
and $x \in \{0,1\}^{\leq i+1}$), but all other 
$\tau_{i,j}(x\#)$  (when $|j-i| > 1$)  are now defined in terms of the 
$\tau_{i,i+1}(x\#)$. 

The classical formulas about transpositions, are now true on a maximal
prefix code (see the Lemma below). For the first formula, the maximal prefix 
code is \ $\{0,1\}^{m+1} \ \cup \ \{0,1\}^{\leq m} \#$, where
$m = {\rm max}\{i,j,k\}$, and for the second formula the maximal prefix code
is \ $\{0,1\}^{j+1} \ \cup \ \{0,1\}^{\leq j} \#$.

\medskip

\noindent {\bf Definition and notation.} \ 
For a group $G$, a subset $\Delta \subseteq G$, and an element 
$g \in \langle \Delta \rangle_G$, we define the
{\em word length of $g$ over} $\Delta$ to be the length of the shortest word
over $\Delta^{\pm 1}$ that is equivalent to $g$ in $G$.
We denote the word length by $|g|_{\Delta}$.

\smallskip

In summary, we proved:

\begin{lem} \  
\label{classicalFormulas} \
As elements of the Thompson-Higman group $G_{3,1}$ the transpositions 
satisfy the following equalities for all $i, j, k \geq 0$:

\smallskip

 \ \ \ \ \ $\tau_{i,j} \ = \ \tau_{i,k} \ \tau_{k,j}$ \

\smallskip

 \ \ \ \ \ $\tau_{i,j} \ = \ $
 $ \tau_{i,i+1} \ \tau_{i+1,i+2} \ \ldots \ \tau_{j-2,j-1} \ \tau_{j-1,j}$
 $ \ \tau_{j-2,j-1} \ \ldots \ \tau_{i+1,i+2} \ \tau_{i,i+1}$ \ \ \
 (when $0 \leq i < j$).

\smallskip

\noindent So the word length of $\tau_{i,j}$ $(0 \leq i < j)$ over the 
alphabet $\{\tau_{k,k+1} : 0 \leq k\}$ is \, $\leq 2(j-i)-1$.
\end{lem}
We also have:

\begin{lem} \  
\label{tau_generated} \
Let $n \geq 0$ and \ $x \in \{0,1\}^*$.

\smallskip

$\tau_{3n+1,3n+2}(x\#) \ = \ \kappa_{321}^{-n} $
 $ \tau_{1,2} \ \kappa_{321}^n(x\#)$, \ \ \ if $|x| \geq 3(n+1)$;

\smallskip

$\tau_{3n+2,3(n+1)}(x\#) \ = \ \kappa_{321}^{-n} \ \kappa_1^{-1} $
            $\tau_{1,2} \ \kappa_1 \ \kappa_{321}^n(x\#)$

\smallskip

   \hspace{1.1in}
  $ \ = \ \kappa_{321}^{-n} \ \tau_{1,3} \ \kappa_{321}^n(x\#)$, \ \ \
if $|x| \geq 3(n+1)+1$;

\smallskip

$\tau_{3(n+1),3(n+1)+1}(x\#) \ = \ $
$\kappa_{321}^{-n} \ \kappa_1^{-1} \ \kappa_2^{-1} \ \tau_{1,2} $
$\kappa_2 \ \kappa_1 \ \kappa_{321}^n(x\#)$ 

\smallskip

  \hspace{1.1in} 
 $ \ = \ \kappa_{321}^{-n} \ \tau_{3,6} \ \kappa_{321}^n(x\#)$, \ \ \
if $|x| \geq 3(n+1)+2$.

\smallskip

\noindent 
Every transposition $\tau_{i-1,i}$  ($i > 0$) has word length $< 2i$ over
$\{ \tau_{0,1}, \tau_{1,2}, \kappa_1, \kappa_2, \kappa_3 \}$, 
and has word length $\leq \lceil \frac{2i}{3} \rceil$ over \
 $\{\tau_{0,1}, \tau_{1,2}, \tau_{3,6}, \kappa_{321} \}$.
\end{lem}
{\bf Proof.} \ On an input $x\#$ as above, we can verify that

\smallskip

$x\# \ = \ x_0 \ x_1 x_2 x_3 \ x_4 x_5 x_6 \ x_7 x_8 x_9 \ \ldots \ $
  $x_{3k+1} x_{3k+2} x_{3(k+1)} \ \ldots \ \#$ \ \ \
$ \stackrel{\kappa_{321}}{\longmapsto} \ $

\smallskip

 \ \ \ \ \  \ \ \ \ \  $x_0 \ x_4 x_5 x_1 \ x_7 x_8 x_2 \ \ldots \ $
$x_{3(k+1)+1} x_{3(k+1)+2} x_{3(k-1)} \ \ldots \ \#$

\medskip

\noindent For any $n \geq 0$ we can then verify the first formula:

\smallskip

$x\# \ = \ x_0 \ x_1 x_2 x_3 x_4 \ \ldots \ x_{3n+1} x_{3n+2} \ \ldots \#$
 \ \ \ $\stackrel{\kappa_{321}^n}{\longmapsto} \ $

\smallskip
$x_0 \ x_{3n+1} x_{3n+2} x_? \ \ldots \ x_? x_? \ \ldots \#$
 \ \ \ $\stackrel{\tau_{1,2}}{\longmapsto} \ $

\smallskip

$x_0 \ x_{3n+2} x_{3n+1}  x_? \ \ldots \ x_? x_? \ \ldots \#$
 \ \ \ $\stackrel{\kappa_{321}^{-n}}{\longmapsto} \ $

\smallskip

$x_0 \ x_1 x_2 x_3 \ \ldots \ x_{3n+2} x_{3n+1} \ \ldots \#$.

\medskip

\noindent Note that $\kappa_3$, $\kappa_2$, $\kappa_1$, and $\tau_{1,2}$
do not change $x_0$.
For the other two formulas the proof is very similar.

\smallskip

For arguments $x\#$ that are ``too short'' we will {\em define}
$\tau_{i,i+1}(x\#)$ by the above formulas (when $1 < i$).
 \ \ \ $\Box$

\bigskip

\noindent {\bf Remark on the definition of the transpositions:}
We defined $\tau_{1,2}$ and $\tau_{0,1}$ earlier, and we gave
formulas that define any $\tau_{i,j}$ in terms of transpositions of the
form $\tau_{n,n+1}$ ($n \geq 0$). So, since the above Lemma defines 
$\tau_{i,i+1}(x\#)$ when $x\#$ is ``too short'', all transpositions are now 
completely defined as elements of $G_{3,1}$.

\medskip

\noindent {\bf Remark on the role of the transpositions:} The transpositions 
are elements of $G_{3,1}$, and $G_{3,1}$ is finitely generated; let 
$\Delta_{3,1}$ be a finite generating set for $G_{3,1}$. So we can write 
each $\tau_{i,i+1}$ as a finite word over $\Delta_{3,1}^{\pm 1}$. 
Why do want to use a generator like $\kappa_{321}$ which doesn't belong to 
$G_{3,1}$? The reason is complexity: 
Over $\Delta_{3,1} \cup \{\kappa_{321}\}$,
the word length of $\tau_{i,i+1}$ has a linear upper bound, but over 
$\Delta_{3,1}$ alone, the word length of $\tau_{i,i+1}$ has a lower bound 
which is exponential in $i$ (as we will prove in Lemma 
\ref{distortionVinTboolLoweBound} and Theorem \ref{distortionVinTbool}).

%%%%%%%%%%%%%%%%%%%%%%%%%%%%%%%%%%%%%%%%%%%%%%%%%%%%%%%%
% Section 3
%%%%%%%%%%%%%%%%%%%%%%%%%%%%%%%%%%%%%%%%%%%%%%%%%%%%%%%%

\section{Simulation of a boolean function by a group element}

One problem in trying to simulate circuits by group elements is that 
the input-output function of a circuit is not necessarily a permutation.
Obtaining permutations is a slightly stronger requirement than the classical
problem of constructing injective (a.k.a.~``reversible'') circuits. See
e.g.\ \cite{Lec}, \cite{Ben73}, \cite{Ben89} for the construction of 
injective Turing machines, and \cite{FredToff} for injective circuits;
the latter reference contains insightful comments on the physical 
significance of injective computing.

To do injective computing with non-injective functions, we apply the  
following transformation from functions to permutations. 
For a function $A \stackrel{f}{\longmapsto} B$,
let $\Gamma_f = \{(x,f(x)) : x \in A\}$ be the graph of the function.
Consider the transformation $\pi$ defined by

\smallskip

 \ \ \ \ \   \ \ \   $\pi: \ (A \stackrel{f}{\longrightarrow} B)$
 \ \ \ \  $\longmapsto$ \ \ \ \  
$(A \cup \Gamma_f \ \stackrel{\pi(f)}{\longrightarrow} \ A \cup \Gamma_f)$,

\smallskip

\noindent where $\pi(f)$ is defined by \
$x \in A \longmapsto (x, f(x)) \in \Gamma_f$, \ and \
$(x, f(x)) \in \Gamma_f \longmapsto x \in A$. Note that $\pi(f)$ is a 
permutation of the set $A \cup \Gamma_f$, for any function $f$.

In programming, functions $f$ are often tacitly replaced by $\pi(f)$
because when an output is computed, people also want to remember the input.
Note also that for two functions  $f_1$ and $f_2$ with same domain set
$A$ and same image set $B$, we have \ $f_1 = f_2$ \ iff \    
$\pi(f_1) = \pi(f_2)$.

\smallskip

In this section we first associate elements of the Thompson-Higman 
group $G_{3,1}$ with the elementary circuit components {\sc not}, 
{\sc or}, {\sc and}, and {\sc fork}. We base this on the above 
transformation $\pi$. Then we define ``simulation'' of an acyclic 
circuit by an element of $G_{3,1}$; an element of $G_{3,1}$ is described
by a sequence of generators. Finally we prove that every acyclic circuit
can be simulated by an element of $G_{3,1}$; moreover, this simulation 
provides a polynomial-time reduction of the equivalence problem of circuits
to the equality problem of elements of $G_{3,1}$, restricted to the subset 
$0 \{0,1,\#\}^*$ of $\{0,1,\#\}^*$ (the {\it word problem with restriction}). 
In the next section we will go further and we reduce the word problem with
restriction to the actual word problem. 

\medskip

With the boolean functions {\sc not}, {\sc or}, and {\sc and},
we associate the following elements of $G_{3,1}$ (described by tables).

\[ \varphi_{\neg} \ = \ \left[ \begin{array}{ccc}
0 \ &\ 1 \ &\ \#  \\
1 \ &\ 0 \ &\ \# 
\end{array}        \right]
\]

\[ \varphi_{\vee} \ = \ \left[ \begin{array}{ccc}
0x_1x_2            \ &\ 1x_1x_2  \ &\ {\rm identity} \\
               \ &\         \ &\   {\rm on} \\
(x_1\vee x_2)x_1x_2\ &\ (\, {\overline{x_1\vee x_2}} \, )x_1x_2 \ &\
\{0,1\}^{\leq 2}\#
\end{array}        \right]
\]

\[ \varphi_{\wedge} \ = \ \left[ \begin{array}{ccc}
0x_1x_2              \ &\ 1x_1x_2 \ &\ {\rm identity} \\
               \ &\         \ &\   {\rm on} \\
(x_1\wedge x_2)x_1x_2\ &\ (\, {\overline{x_1\wedge x_2}} \, )x_1x_2 \ &\
\{0,1\}^{\leq 2}\#
\end{array}        \right]
\]
where $x_1, x_2$ range over $\{0,1\}$. Hence the domain and image codes of
$\varphi_{\vee}$ and $\varphi_{\wedge}$ are all equal to \ 
$\{0,1\}^3 \ \cup \ \{0,1\}^{\leq 2} \#$. 

The three functions above are length-preserving: \ 
 $|\varphi_{\wedge}(x)| = |x|$ for all 
$x \in {\rm Dom}(\varphi_{\wedge})$, and similarly for $\varphi_{\vee}$ 
and $\varphi_{\neg}$.

\bigskip

\noindent In order to represent the {\sc fork} function in circuits by an
element of $G_{3,1}$ a first idea would be to define a ``0-fork'' element 
of $G_{3,1}$ (which duplicates a leading 0), as follows: \  

\smallskip
 
\( \varphi_{\rm 0f} \ = \ \left[ \begin{array}{ccccc}
0  \ & \ \#  \ &\ 10 \ & \ 1\# \ & \ 11   \\ 
00 \ & \ 0\# \ &\ 01 \ & \ \#  \ & \ 1
\end{array}        \right] \). 

\smallskip

\noindent Then, \ 
 $\tau_{0,1} \ \varphi_{\vee} \ \varphi_{\rm 0f}(0x) = 0xx$ \
(for all $x \in \{0,1\}$), so we could use this as a way to represent the
{\sc fork} operation in a circuit. 

\smallskip

However, it will turn out later that what we need is a forking operation 
that preserves the string length modulo 3. Thus, we define a 
``four-fold 0-fork'' element of $G_{3,1}$ (which turns a leading 0 into four
leading 0s).
\[ \varphi_{\rm 0f,4} \ = \ \left[ \begin{array}{ccc ccc ccc}
0 \ & \ \# \ &\ 10 \ & \ 1\#\ &\ 1^20 \ &\ 1^2\#\ &\ 1^30\ &\ 1^3\#\ &\ 1^4 \\  
0^4\ &\ 0^3\#\ &\ 01\ &\ 0\#\ &\ 0^21\ &\ 0^2\#\ &\ 0^31\ & \   \# \ &\ 1 
\end{array}        \right] \]
We have \ 
${\rm domC}(\varphi_{\rm 0f,4}) \ = \ 1^{\leq 3}\{0,\#\} \cup \{1^4\}$, 
and \  
${\rm imC}(\varphi_{\rm 0f,4}) \ = \ 0^{\leq 3}\{1,\#\} \cup \{0^4\}$. From 
the definitions one immediately verifies the following.
\begin{lem}  \  
The maps $\tau_{i,j}$ (where $0 \leq i < j$), $\varphi_{\rm 0f,4}$,
$\varphi_{\neg}$, $\varphi_{\vee}$, $\varphi_{\wedge}$ belong to the 
Thompson-Higman group $G_{3,1}$, they stabilize the sets $\{0,1\}^*$ and  
$\{0,1\}^*\#$, they preserve lengths modulo $3$, they map $0 \, \{0,1\}^*$
into itself, and they map $0 \, \{0,1\}^* \#$ into itself.
\end{lem}

\noindent {\bf Notation:} \ Let  $G \subseteq {\mathcal G}_{3,1}$;  note that 
``$\subseteq$'' means that $G$ isn't just a subgroup, but a particular 
embedding into ${\mathcal G}_{3,1}$ is considered.
By $G^{\rm mod \, 3}$ we denote the subgroup

\smallskip

 \ $\{ \varphi \in G \ : \ \forall x \in \{0,1\}^*, \ $
    $|\varphi(x)| \equiv |x| \ {\rm mod} \ 3\}$, \ 

\smallskip

\noindent i.e., the elements of $G$ that, when restricted to $\{0,1\}^*$,
preserve the length of strings modulo 3. 
In particular, we will use the notation $G_{3,1}^{\rm mod \, 3}$ for the 
corresponding subgroup of the Thompson-Higman group $G_{3,1}$.
 
\bigskip

We point out that $\varphi_{\neg}$, $\varphi_{\vee}$, $\varphi_{\wedge}$, 
and all $\tau_{i,j}$ $(0 \leq i \leq j)$ are length-preserving, and that 
$\varphi_{\rm 0f,4}$ preserves length modulo 3. We will not use any other 
elements of $G_{3,1}$ in the constructions and proofs in this Section.

In order to obtain computational results we describe boolean functions by
acyclic circuits, and we describe elements of $G_{3,1}$ by words.
Let us choose a finite set of generators $\Delta_{3,1}$ of the group
$G_{3,1}$.  For  $G_{3,1}$ we also use the infinite generating set \ 
$\Delta_{3,1} \cup \{\tau_{i,i+1} : 0 \leq i \}$.

Let $C$ be an acyclic boolean circuit with $m$ input variables
$x_1, \ldots, x_m$ and $n$ output variables $y_1, \ldots, y_n$. Let \
$f_C: \{0,1\}^m \to \{0,1\}^n$ \ be the input-output function of $C$.
Hence, two circuits $C_1$ and $C_2$ are equivalent iff $f_{C_1} = f_{C_2}$.

Our definition of ``simulation'' is a variation of the above transformation
$\pi$.
\begin{defn}
\label{DEFsimulat}  \          An element 
$\Phi_f \in G_{3,1}^{\rm mod \, 3}$ {\bf simulates} a boolean function \ 
$f: \{0,1\}^m \to \{0,1\}^n$ \  iff

\smallskip

\noindent $\bullet$ \ the domain code and the image code of $\Phi_f$ are
subsets of \ $\{0,1\} \, \{0,1\}^* \, \cup \, \{0,1\}^* \, \#$  

\smallskip

\noindent $\bullet$ \ $\Phi_f$ maps $0\{0,1\}^m$ into 
$0^{1 + i(n)}\{0,1\}^{n+m}$ in such a way that 

\smallskip

 \ \ \ \ \  \ \ \ \ \  $\Phi_f(0 \, x_1 \, \ldots \, x_m) \ = \ $
  $0^{1 + i(n)} \ f(x_1, \ldots, x_m) \ \ x_1 \ \ldots \ x_m $  

\smallskip

\noindent where $i(n) \in \{0,1,2\}$ is such that $1+ n +i(n)$ is a multiple 
of $3$ (i.e., \ $i(n) \equiv -(1+n)$ {\rm mod 3}); so the role of $i(n)$ is 
to make $\Phi_f$ preserve lengths modulo $3$;

\medskip

\noindent
$\bullet$ \ $\Phi_f$ and $\Phi_f^{-1}$ map the set $\{0,1\}^*$ into itself,
and map $\{0,1\}^*\#$ into itself; moreover,  
 \ $\Phi_f$ maps the set $0\{0,1\}^*$ into itself, and $\Phi_f^{-1}$ maps
the set $1\{0,1\}^*$ into itself.

\medskip

\noindent
When $\Phi_f$ is represented by a {\em word} $w_f$ over \  
$\Delta_{3,1}^{\pm 1} \cup \{\tau_{i,i+1} : 0 \leq i\}$ \ we say that $w_f$
simulates $f$.
\end{defn}
A boolean function $f$ can be simulated by many elements of $G_{3,1}$. 

\medskip

By the above definition, if \ $w \in$ 
$ \{0,1\}^{\geq m} \cup \{0,1\}^{\geq m}\#$ \   
then $\Phi_f(0 w)$ tells us the value of $f$ on input $x_1 \ldots x_m$ (where 
$x_1 \ldots x_m$ is the prefix of length $m$ of $w$). The definition does 
not give any connection between $\Phi_f(0 \, x_1 \ldots x_k \, \#)$ and $f$ 
when $k < m$ (where $x_1, \ldots , x_k \in \{0,1\}$);
we call this the ``case when the input is too short''. In some applications 
we want such a connection, hence we will need the definition of ``strong
simulation'' below. 
(We cannot do much about the fact that $\Phi_f(1 w)$ has no connection with 
$f$; since $\Phi_f$ is an element of $G_{3,1}$, it is a bijection between 
maximal prefix codes, whereas $f$ need not be injective nor surjective.
So there has to be a big difference between $\Phi_f$ and $f$ somewhere.)

\begin{defn}
\label{DEFstrong_simulat}  \
We say that $\Phi_f$ {\bf strongly simulates} $f$ iff in addition to the 
conditions of simulation (Definition \ref{DEFsimulat}), we have for all
$0 \leq k < m$: \  $\Phi_f(0 \, x_1 \, \ldots \, x_k \, \#)$ is defined
for all $x_1 \ldots x_k \in \{0,1\}^k$.

So for strong simulation, $\Phi_f(0 \, x_1 \, \ldots \, x_k \, \#)$ depends
only on the function $f$ and on $k$ and on $x_1 \ldots x_k$; it does not 
depend on any particular circuit used to compute $f$.
\end{defn}

The next Lemma follows immediately from the definition of simulation. It
gives a connection between the equivalence problem of circuits and the 
{\it word problem with restriction} of $G_{3,1}$.
For a Thompson group $G$ $(\subset {\mathcal G}_{3,1})$ with generating set 
$A$, and a subset $S \subseteq \{0,1,\#\}^*$, the word problem with 
restriction is defined as follows: \\  
{\sc Input:} Two words $u,v$ over $A^{\pm 1}$. \\
{\sc Question:} Are the partial functions described by $u$ and $v$ the same
when restricted to $S$?  

\medskip

\noindent  We denote the restriction of a partial function $F$ to a set $S$
by $F|_S$.  The next Lemma follows immediately from Definitions 
\ref{DEFsimulat} and \ref{DEFstrong_simulat}.

\begin{lem} \label{simulatVSfix} \   
Let $f$ and $g$ be any boolean functions with the same number of input
variables and the same number of output variables.
If $f$ and $g$ are simulated by $\Phi_f$, respectively $\Phi_g$, then
we have

\medskip

 \ \ \ \ \    $f = g$ \ \ \  iff \ \ \  
      $(\Phi_f)|_{0 \{0,1,\# \}^*} \ = \ (\Phi_g)|_{0 \{0,1,\# \}^*}$

\medskip

\noindent In the case of strong simulation we have, in addition,  

\medskip

 \ \ \ \ \  $f = g$  \ \ \  iff \ \ \  
 $(\Phi_f)|_{\{0,\#\}\{0,1,\# \}^*} \ = \ $
     $(\Phi_g)|_{\{0,\#\}\{0,1,\# \}^*}$ \ \   
\end{lem}

\medskip

Let $\Delta_{3,1}$ be a finite set of generators of the group $G_{3,1}$.  
For  $G_{3,1}$ we also use the infinite generating set 
 \ $\Delta_{3,1} \cup \{\tau_{i,i+1} : 0 \leq i \}$.
With every acyclic boolean circuit $C$ we want to associate a word $w_C$   
over the alphabet $\Delta_{3,1}^{\pm 1} \cup \{\tau_{i,i+1} : 0 \leq i \}$, 
and we want the correspondence $C \mapsto w_C$ to be polynomial-time 
computable. For every word $w$ over \  
$\Delta_{3,1}^{\pm 1} \cup \{\tau_{i,j} : 0 \leq i < j\}$ \ we denote the 
length of $w$ by $|w|$, and we denote the largest subscript in any 
$\tau_{i,j}$ occurring in $w$ by $J_w$.
 
\medskip

The {\it size} of an acyclic boolean circuit $C$ is denoted by $|C|$; if $C$
has $k_1$ gates of type {\sc not} or {\sc fork},
$k_2$ gates of type {\sc and} or {\sc or}, and $n$ output variables,
the size of $C$  is defined to be \ $|C| = k_1 + 2 \cdot k_2 + n$.
Equivalently, $|C|$ is the number of connections (wires, or edges in the
circuit graph) between gates or from an input/output port to a gate (for 
that reason, gates with two input variables are counted twice). 
(Our definition of the size $|C|$ is slightly different from the traditional
definition, which just counts {\sc not}, {\sc and}, {\sc or} gates and I/O
ports, but it is linearly related to the traditional definition.)

In an acyclic circuit every gate, and also every input or output variable,
can be assigned a {\it level} (or ``layer'', or ``depth'').
The input variables of the circuit have level 0.
A gate or an output variable has level 1 iff only input variables of the
circuit feed into it.
A gate or an output variable has level $\ell$ iff it receives input from 
levels $< \ell$ only, and at least one of its inputs comes from level 
$\ell - 1$. The maximum level of any output variable is called the 
{\it depth} of the circuit.

\begin{thm} \label{reduction} \ 
There is an injective function  \ $C  \mapsto w_C$ \ 
from the set of acyclic boolean circuits to the set of words over the 
alphabet $\Delta_{3,1}^{\pm 1} \cup \{\tau_{i,i+1} : 0 \leq i \}$ (where 
$\Delta_{3,1}$ is a finite generating set of $G_{3,1}$), with the 
following properties:

\smallskip

\noindent 
(1) \ \ $w_C$ strongly simulates $f_C$.

\smallskip

\noindent 
(2) \ \  The length of $w_C$ satisfies \ \ $|w_C|< c \, |C|^4 + c$ \ \ 
  (for some positive constant $c$), 

 \ \ \ and the largest subscript $J_{w_C}$ of any $\tau_{i,i+1}$ in $w_C$ 
   satisfies \ \   $J_{w_C} \leq c \, |C|^2 + c$.

\smallskip

\noindent 
(3) \ \ $w_C$ is computable from $C$ in polynomial time, as a
polynomial in $|C|$. 
\end{thm}
To make sense of the phrase ``$w_C$ is computable $\ldots$'', we need to
represent any transposition $\tau_{i,i+1}$ (with $i \in {\mathbb N}$)
by a string over a finite alphabet; we simply write the integer $i$ in unary 
notation (i.e., $i$ is represented by the string $0^i$).

\smallskip

\noindent {\bf Proof.} \ 
We assume that the elements $\varphi_{\neg}$, $\varphi_{\vee}$, 
$\varphi_{\wedge}$, $\varphi_{\rm 0f,4}$ and $\tau_{0,1}$ belong to 
$\Delta_{3,1}$. 
If this were not the case, we could express these by fixed words over 
another finite generating set of $G_{3,1}$.

We can assume that our acyclic circuits are {\it strictly layered}, i.e.,
a gate or an output variable at level $\ell$ only receives inputs from level 
$\ell-1$. Hence, all the output variables of the circuit are at the same 
level $L$ ($L$ is the depth of the circuit). 
If the layering of a circuit $C$ is not strict, we can insert 
{\it identity gates} to enforce strictness. An identity gate has one input
variable and one output variable, connected by a wire; the two variables carry
the same boolean value. In the present proof we will
count these identity gates as gates in the definition of circuit size.
In order to make a circuit $C$ strictly layered, fewer than $|C|^2$ identity
gates need to be introduced. (Indeed, for each gate $g$ we add at most as 
many identity gates as the depth of this gate $g$; so, in total we add at 
most \ $|C| \cdot {\rm depth}(C) \ (\leq \ |C|^2)$ \ identity gates).
So the size increase is polynomially bounded. Moreover, identity gates will
not affect $w_C$, as we will see in the construction of $w_C$. 

A circuit $C$ has input variables $x_1, \ldots, x_m$, output variables 
$y_1, \ldots, y_n$, and {\it internal variables} which correspond to the 
boolean values carried by internal wires (between gates or between a gate 
and an input or an output port). The internal variables at level $\ell$
are denoted $y_1^{\ell}$, $y_2^{\ell}$, $\ldots$, $y_{n_{\ell}}^{\ell}$.
When $\ell = L$ (output level) we have $n_L = n$ and $y_i^L = y_i$; and when 
$\ell = 0$ (input level) we have $n_0 = m$ and $y_i^0 = x_i$. 
For every level $\ell$ $(1 \leq \ell \leq L)$, we consider a circuit 
$C_{\ell}$ (called the {\it slice of} $C$ at level $\ell$). The input
variables of $C_{\ell}$ are $y_1^{\ell-1}$, $\ldots$,
$y_{n_{\ell-1}}^{\ell-1}$, and the output variables are $y_1^{\ell}$, 
$\ldots$, $y_{n_{\ell}}^{\ell}$; the gates of $C_{\ell}$ are all the gates of
$C$ at level $\ell$.  

It will be convenient to use the notation \ 
$Y^{\ell} = y_1^{\ell} y_2^{\ell} \ \ldots y_{n_{\ell}}^{\ell}$ \ 
(concatenation of all the variables $y_i^{\ell}$), for $0 \leq \ell \leq L$. 

\medskip

In order to define $w_C$ let us first consider the case when 
$L = 1$, i.e, the circuit consists of just one slice. 

Let $k \geq 0$ and assume that for every circuit $C$ of depth 1 and of size 
$|C| \leq k$ (where identity gates are counted as well), we can 
compute a word $w_C$ (over the alphabet 
$\Delta_{3,1}^{\pm 1} \cup \{\tau_{i,i+1} : 0 \leq i \}$).

Any circuit $C$ of depth 1 and of size $k+1$ can be viewed as a circuit $K$ 
of depth 1 and of size $\leq k$, with an additional gate ({\sc and}, 
{\sc or}, {\sc not}, identity, or {\sc fork}). 
Let $x_1, \ldots, x_m$ be the input 
variables and let $y_1, \ldots, y_n$ be the output variables of $K$.
 
\medskip

{\sc Case 1:} \  Suppose our circuit $C$ is obtained from $K$ by adding an 
identity gate or a {\sc not} gate, with new input variable $x_{m+1}$ and 
new output variable $y_{n+1}$. Note that only one wire can be connected to
an input variable $x_i$; we use explicit {\sc fork} operations when we 
want to duplicate a variable.   
In case a {\sc not} gate is added, the input-output function of the new 
circuit is \ 
$f_C(x_1, \ldots, x_m, x_{m+1}) = $
$(y_1, \ldots, y_n, \overline{x_{m+1}})$,
where $f_K(x_1, \ldots, x_m) = (y_1, \ldots, y_n)$. 
The boolean function $f_C$ is to be simulated by a Thompson group element \
$\Phi_f: \{0,1\}^* \to \{0,1\}^*$ \ such that \

\smallskip

$\Phi_f(0 \, x_1 \ldots x_m, x_{m+1}) \ = \ $
$0^{1+i(n+1)} \, y_1 \ldots y_n \ \overline{x_{m+1}} \ $
$x_1 \ldots x_m x_{m+1}$ \

\smallskip

\noindent for all $x_1, \ldots, x_m, x_{m+1} \in \{0,1\}$, and such that
$\Phi_f$ has the stability properties of Definition \ref{DEFsimulat}; recall
(as we saw in the Definition of ``simulation'') that $i(n) \equiv -(n+1)$  
mod 3, $i(n) \in \{0,1,2\}$.

Let $w_K$ and $\Phi_{f_K} \in G_{3,1}$ be the simulation of $f_K$, 
which exists by induction. We proceed as follows:

\smallskip

$0\, x_1 \ldots x_m \ x_{m+1} \ $
$\stackrel{\Phi_{f_K} }{\longmapsto} \ $
$0^{1+i(n)} \ y_1 \ldots y_n \ x_1 \ldots x_m \ x_{m+1} \ $

\medskip

\noindent Case $i(n) = 1$: \    
In this case we continue the simulation of $f_C$ as follows. 

\smallskip

$0 \ 0 \ y_1 \ y_2 \ldots y_n \ x_1 \ldots x_m \ x_{m+1} \ $
$\stackrel{\tau_{2,n+m+2}}{\longmapsto} \ $
$0 \ 0 \ x_{m+1} \ y_2 \ldots y_n \ x_1 \ldots x_m \ y_1$
$\stackrel{\varphi_{\vee}}{\longmapsto} \ $

\smallskip

$x_{m+1} \ 0 \ x_{m+1} \ y_2 \ldots y_n \ x_1 \ldots x_m \ y_1$
$\stackrel{\varphi_{\neg}}{\longmapsto} \ $
$\overline{x_{m+1}} \ 0 \ x_{m+1} \ y_2 \ldots y_n \ x_1 \ldots x_m \ y_1$
$\stackrel{\tau_{2,n+m+2}}{\longmapsto} \ $

\smallskip

$\overline{x_{m+1}} \ 0 \ y_1y_2 \ldots y_n \ x_1 \ldots x_m \ x_{m+1}$

\smallskip

\noindent Applying \ 
$\tau_{n+1,n+2} \ \tau_{n,n+1} \ \dots \ \tau_{1,2} \ \tau_{0,1} (\cdot)$ 
 \ then yields 

\smallskip

$0 \ y_1 \ldots y_n \ \overline{x_{m+1}} \ x_1 \ldots x_m x_{m+1}$.

\smallskip

\noindent Thus our circuit $C$ is simulated by the following word

\smallskip

$w_C \ = \ \tau_{n+1,n+2} \ \dots \ \tau_{1,2} \ \tau_{0,1}$
 $ \ \varphi_{\neg} \ \varphi_{\vee} \ $ 
 $ \tau_{2,n+m+2} \ w_K$.

\medskip

\noindent Case $i(n) = 2$: \    
In this case we continue the simulation of $f_C$ as follows.

\smallskip

$0 \ 0 \ 0 \ y_1 \ldots y_n \ x_1 \ldots x_m \ x_{m+1} \ $
$\stackrel{\tau_{2,n+m+3}}{\longmapsto} \ $
$0 \ 0 \ x_{m+1} \ y_1 \ldots y_n \ x_1 \ldots x_m \ 0$
$\stackrel{\varphi_{\vee}}{\longmapsto} \ $

\smallskip

$x_{m+1} \ 0 \ x_{m+1} \ y_1 \ldots y_n \ x_1 \ldots x_m \ 0$
$\stackrel{\varphi_{\neg}}{\longmapsto} \ $
$\overline{x_{m+1}} \ 0 \ x_{m+1} \ y_1 \ldots y_n \ x_1 \ldots x_m \ 0$
$\stackrel{\tau_{2,n+m+3}}{\longmapsto} \ $

\smallskip

$\overline{x_{m+1}} \ 0 \ 0y_1 \ldots y_n \ x_1 \ldots x_m \ x_{m+1}$

\noindent Applying \
$\tau_{n+1,n+2} \ \tau_{n,n+1} \ \dots \ \tau_{1,2} \ \tau_{0,1} (\cdot)$
 \ then yields

\smallskip

$0 \ 0 \ y_1 \ldots y_n \ \overline{x_{m+1}} \ x_1 \ldots x_m x_{m+1}$.

\medskip

\noindent Case $i(n) = 0$: \    
In this case we continue the simulation of $f_C$ as follows.

\smallskip

$0 \ y_1 \ldots y_n \ x_1 \ldots x_m \ x_{m+1} \ $
$\stackrel{\varphi_{\rm 0f, 4}}{\longmapsto} \ $
$0 \ 0 \ 0 \ 0 \ y_1 \ldots y_n \ x_1 \ldots x_m \ x_{m+1} \ $
$\stackrel{\tau_{2,n+m+4}}{\longmapsto} \ $

\smallskip

$0 \ 0 \ x_{m+1} \ 0 \ y_1 \ldots y_n \ x_1 \ldots x_m \ 0$
$\stackrel{\varphi_{\vee}}{\longmapsto} \ $
$x_{m+1} \ 0 \ x_{m+1} \ 0 \ y_1 \ldots y_n \ x_1 \ldots x_m \ 0$
$\stackrel{\varphi_{\neg}}{\longmapsto} \ $

\smallskip

$\overline{x_{m+1}} \ 0 \ x_{m+1} \ 0 \ y_1 \ldots y_n \ x_1 \ldots x_m \ 0$
$\stackrel{\tau_{2,n+m+4}}{\longmapsto} \ $

\smallskip

$\overline{x_{m+1}} \ 0 \ 0 \ 0 \ y_1 \ldots y_n \ x_1 \ldots x_m \ x_{m+1}$

\smallskip

\noindent Applying \
$\tau_{n+2,n+3} \ \tau_{n+1,n+2} \ \dots \ \tau_{1,2} \ \tau_{0,1} (\cdot)$
 \ then yields

\smallskip

$0 \ 00 \ y_1 \ldots y_n \ \overline{x_{m+1}} \ x_1 \ldots x_m x_{m+1}$.

\bigskip

\noindent The case where, instead of a {\sc not} gate, an identity gate 
is added is similar (except that we simply omit $\varphi_{\neg}$). 

In any case the length of $w_C$ over the alphabet 
$\Delta_{3,1}^{\pm 1} \cup \{\tau_{i,i+1} : 0 \leq i \}$ is at most \  
 $|w_K|+ 2 \, |\tau_{2,n+m+4}|+ 4 + n +2$. By Lemma \ref{classicalFormulas}, 
 \ $|\tau_{2,n+m+4}| \leq 2 \, (n+m+2) -1$. Hence, \ 
$|w_C| \leq |w_K|+ 4m + 5n +12$. 
Moreover, the subscripts of the transpositions appearing in $w_C$ are \ 
$\leq {\rm max}\{n+m+4, J_K\}$, where $J_K$ is the largest subscript in any
transposition appearing in $w_K$.
 
In case we want to change the positions of the added variables $x_{m+1}$
and $y_{n+1}$ (so that $x_{m+1}$ is the $i$th input variable and $y_{n+1}$
is the $j$th output variable), we apply other appropriate permutations 
(instead of  \ 
$\tau_{n+2,n+3} \ \tau_{n+1,n+2} \ \dots \ \tau_{1,2} \ \tau_{0,1}$ \ 
and $\tau_{2,n+m+4}$ above). This does not change our upper bound on $|w_C|$.

\bigskip

{\sc Case 2:} \  Suppose our circuit $C$ (still of depth 1) is obtained by 
adding an {\sc and} gate or an {\sc or} gate to $K$, with new output variable 
$y_{n+1}$ and new input variables $x_{m+1}, x_{m+2}$. Recall that only 
one wire can be connected to an input variable $x_i$.  We only deal with the 
{\sc or} case (the {\sc and} case being practically the same).
The input-output function of the new circuit is

\smallskip

$f_C(x_1, \ldots, x_m, x_{m+1}, x_{m+2}) \ = \ $
$(y_1, \ldots, y_n, \ x_{m+1} \vee x_{m+2})$,

\smallskip

\noindent where $f_K(x_1, \ldots, x_m) = (y_1, \ldots, y_n)$.
The boolean function $f_C$ is to be simulated by a Thompson group element \
$\Phi_f: \{0,1\}^* \to \{0,1\}^*$ \ such that \

\smallskip

$\Phi_f(0 \, x_1 \ldots x_m \ x_{m+1} x_{m+2}) \ = \ $
$0^{1+i(n+1)} \ y_1 \ldots y_n \ (x_{m+1} \vee x_{m+2}) \ $
 $ x_1 \ldots x_m x_{m+1} x_{m+2}$ \

\smallskip

\noindent for all $x_1, \ldots, x_m, x_{m+1}, x_{m+2} \in \{0,1\}$, and such 
that $\Phi_f$ has the stability properties of Definition \ref{DEFsimulat}.
Let $w_K$ and $\Phi_{f_K} \in G_{3,1}$ be the simulation of $f_K$, which
exists by induction. Then

\smallskip

$0\ x_1 \ldots x_m \ x_{m+1} x_{m+2} \ $
$\stackrel{\Phi_{f_K} }{\longmapsto} \ $
$0^{1+i(n)} \ y_1 \ldots y_n \ x_1 \ldots x_m \ x_{m+1} x_{m+2} \ $

\medskip

\noindent Case $i(n) = 1$: \ The simulation continues as follows.

\smallskip

$00 \ y_1 \, y_2 \ldots y_n \ x_1 \ldots x_m \ x_{m+1} x_{m+2} \ $
$\stackrel{\tau_{1,n+m+2}}{\longmapsto} \ $
$ \stackrel{\tau_{2,n+m+3}}{\longmapsto} \ $
$0 \ x_{m+1} x_{m+2} \ y_2 \ldots y_n \ x_1 \ldots x_m \ 0 \ y_1 \ $
$\stackrel{\varphi_{\vee}}{\longmapsto} $

\smallskip

$(x_{m+1} \vee x_{m+2}) \ x_{m+1} x_{m+2} \ y_2 \ldots y_n \ $
$x_1 \ldots x_m \ 0 \ y_1 \ $
$\stackrel{\tau_{1,n+m+2}}{\longmapsto} \ $
$\stackrel{\tau_{2,n+m+3}}{\longmapsto} \ $

\smallskip

$(x_{m+1} \vee x_{m+2}) \ 0 \ y_1 y_2 \ldots y_n \ x_1 \ldots x_m $
   $x_{m+1} x_{m+2}$ .

\smallskip

\noindent By applying \ 
$ \tau_{n+1,n+2} \ \ldots \ \tau_{1,2} \ \tau_{0,1}$ \  
we obtain

\smallskip

$0 \ y_1 y_2 \ldots y_n \ (x_{m+1} \vee x_{m+2}) \ x_1 \ldots x_m $
$x_{m+1} x_{m+2}$. 

\medskip

\noindent Case $i(n) = 2$: \ The simulation continues as follows.

\smallskip

$000 \ y_1 \, y_2 \ldots y_n \ x_1 \ldots x_m \ x_{m+1} x_{m+2} \ $
$\stackrel{\tau_{1,n+m+3}}{\longmapsto} \ $
$ \stackrel{\tau_{2,n+m+4}}{\longmapsto} \ $
$0 \ x_{m+1} x_{m+2} \ y_1 \ldots y_n \ x_1 \ldots x_m \ 0 \ 0 $
$\stackrel{\varphi_{\vee}}{\longmapsto} $

\smallskip

$(x_{m+1} \vee x_{m+2})\ x_{m+1}x_{m+2}\ y_1 \ldots y_n$
                               $ \ x_1 \ldots x_m\ 0 \ 0$
$\stackrel{\tau_{1,n+m+3}}{\longmapsto} \ $
$\stackrel{\tau_{2,n+m+4}}{\longmapsto} \ $

\smallskip

$(x_{m+1} \vee x_{m+2}) \ 00 \ y_1 y_2 \ldots y_n \ x_1 \ldots x_m $
   $x_{m+1} x_{m+2}$ .

\smallskip

\noindent By applying \
$\tau_{n+2,n+3} \ \ldots \ \tau_{1,2} \ \tau_{0,1}$ \ 
we obtain

\smallskip

$0 \ 0 \ y_1 y_2 \ldots y_n \ (x_{m+1} \vee x_{m+2}) \ x_1 \ldots x_m $
$x_{m+1} x_{m+2}$.

\medskip

\noindent Case $i(n) = 0$: \ The simulation continues as follows.

\smallskip

$0 \ y_1 \, y_2 \ldots y_n \ x_1 \ldots x_m \ x_{m+1} x_{m+2} \ $
$\stackrel{\varphi_{\rm 0f,4}}{\longmapsto} \ $

\smallskip

$0 000 \ y_1 \, y_2 \ldots y_n \ x_1 \ldots x_m \ x_{m+1} x_{m+2} \ $
$\stackrel{\tau_{1,n+m+4}}{\longmapsto} \ $
$\stackrel{\tau_{2,n+m+5}}{\longmapsto} \ $

\smallskip

$0 \ x_{m+1} x_{m+2} \ 0 \ y_1 \ldots y_n \ x_1 \ldots x_m \ 0 \ 0 $
$\stackrel{\varphi_{\vee}}{\longmapsto} $

\smallskip

$(x_{m+1} \vee x_{m+2})\ x_{m+1}x_{m+2}\ 0 y_1 \ldots y_n \ $
              $x_1 \ldots x_m \ 0 \ 0 $
$\stackrel{\tau_{1,n+m+4}}{\longmapsto} \ $
$\stackrel{\tau_{2,n+m+5}}{\longmapsto} \ $

\smallskip

$(x_{m+1} \vee x_{m+2})\ 000 \ y_1 \ldots y_n \ $
              $x_1 \ldots x_m \ x_{m+1}x_{m+2}$

\smallskip

\noindent By applying \
$\tau_{n+3,n+4} \ \ldots \ \tau_{1,2} \ \tau_{0,1}$ \
we obtain

\smallskip

$0 \ 00 \ y_1 y_2 \ldots y_n \ (x_{m+1} \vee x_{m+2}) \ x_1 \ldots x_m $
$x_{m+1} x_{m+2}$.

\medskip

\noindent Thus our circuit $C$ is simulated by the word $w_C$ of length \ 
$\leq |w_K| + 8m + 9n + 15$ \  over the alphabet
$\Delta_{3,1}^{\pm 1} \cup \{\tau_{i,i+1} : 0 \leq i \}$.
Moreover, the subscripts of the transpositions appearing in $w_C$ are \ 
$\leq {\rm max}\{n+m+5, J_K\}$, where $J_K$ is the largest subscript in any
transposition appearing in $w_K$.

\smallskip

In case we want to change the positions of the added variables $x_{m+1}$,
$x_{m+2}$ and $y_{n+1}$ (so that $x_{m+1}$ is the $i_1$th input variable,
$x_{m+2}$ is the $i_2$th input variable, and $y_{n+1}$ is the $j$th output 
variable), we apply other appropriate permutations (instead of \ 
$\tau_{n+3,n+4} \ \ldots \ \tau_{1,2} \ \tau_{0,1}$, 
$\tau_{2,n+m+5}$, and $\tau_{1,n+m+4}$). This will not change our upper
bounds on $|w_C|$ and $J_C$. 
 
\bigskip 

{\sc Case 3:} \ Suppose our circuit $C$ (still of depth 1) is obtained by 
adding a {\sc fork} gate with a new input variable $x_{m+1}$ and two new 
output variables $y_{n+1}$ and $y_{n+2}$. The input-output function of the 
new circuit is

\smallskip

$f_C(x_1, \ldots, x_m, x_{m+1}) \ = \ (y_1, \ldots, y_n, x_{m+1}, x_{m+1})$,

\smallskip

\noindent where $f_K(x_1, \ldots, x_m) = (y_1, \ldots, y_n)$.
The boolean function $f_C$ is to be simulated by a Thompson group element \
$\Phi_f$ \ such that \

\smallskip

$\Phi_f(0 \, x_1 \ldots x_m x_{m+1}) \ = \ $
$ 0^{1+i(n+2)} \ y_1 \ldots y_n\ x_{m+1}x_{m+1} \ x_1 \ldots x_m x_{m+1}$ \

\smallskip

\noindent for all $x_1, \ldots, x_m, x_{m+1} \in \{0,1\}$, 
$i(n+2) = -n$ mod 3,  and such that
$\Phi_f$ has the stability properties of Definition \ref{DEFsimulat}.
Let $w_K$ and $\Phi_{f_K} \in G_{3,1}$ be the simulation of $f_K$, which
exists by induction. Then

\smallskip

$0\, x_1 \ldots x_m \ x_{m+1} \ $
$\stackrel{\Phi_{f_K} }{\longmapsto} \ $
$0^{1+i(n)} \ y_1 \ldots y_n \ x_1 \ldots x_m \ x_{m+1} \ $

\medskip

\noindent Case $i(n) = 2$: We continue the simulation with 

\smallskip

$000 \ y_1 \ldots y_n \ x_1 \ldots x_m \ x_{m+1} \ $
$\stackrel{\tau_{2,n+m+3}}{\longmapsto} \ $
$00 \ x_{m+1} \ y_1 \ldots y_n \ x_1 \ldots x_m \ 0 \ $
$\stackrel{\varphi_{\vee}}{\longmapsto} \ $

\smallskip

$x_{m+1} \ 0 \ x_{m+1} \ y_1 \ldots y_n \ x_1 \ldots x_m \ 0 \ $
$\stackrel{\tau_{0,1}}{\longmapsto} \ $
$0 \ x_{m+1}x_{m+1} \ y_1 \ldots y_n \ x_1 \ldots x_m \ 0 \ $
$\stackrel{\varphi_{\vee}}{\longmapsto} \ $

\smallskip

$x_{m+1}x_{m+1}x_{m+1} \ y_1 \ldots y_n \ x_1 \ldots x_m \ 0 \ $
$\stackrel{\tau_{2,n+m+3}}{\longmapsto} \ $
$x_{m+1}x_{m+1} \ 0 \ y_1 \ldots y_n \ x_1 \ldots x_m \ x_{m+1} \ $

\smallskip

\noindent Finally we apply \   
$\tau_{n+1,n+2} \ \ldots \ \tau_{3,4} \ \tau_{1,2}$ \ and \  
$\tau_{n,n+1} \ \ldots \ \tau_{2,3} \ \tau_{0,1}$ \   
to obtain 

\smallskip

$0 \ y_1 y_2 \ldots y_n \ x_{m+1} x_{m+1} \ x_1 \ldots x_m x_{m+1}$. 

\medskip

\noindent Case $i(n) = 0$: We continue the simulation with 

\smallskip

$0\ y_1 \ldots y_n \ x_1 \ldots x_m \ x_{m+1} \ $
$\stackrel{\varphi_{\rm 0f,4}}{\longmapsto} \ $
$0000 \ y_1 \ldots y_n \ x_1 \ldots x_m \ x_{m+1} \ $
$\stackrel{\tau_{2,n+m+4}}{\longmapsto} \ $

\smallskip

$00 \ x_{m+1} \ 0 \ y_1 \ldots y_n \ x_1 \ldots x_m \ 0 \ $
$\stackrel{\varphi_{\vee}}{\longmapsto} \ $
$x_{m+1} \ 0 \ x_{m+1} \ 0\ y_1 \ldots y_n \ x_1 \ldots x_m \ 0 \ $
$\stackrel{\tau_{0,1}}{\longmapsto} \ $

\smallskip

$0 \ x_{m+1}x_{m+1} \ 0\ y_1 \ldots y_n \ x_1 \ldots x_m \ 0 \ $
$\stackrel{\varphi_{\vee}}{\longmapsto} \ $
$x_{m+1}x_{m+1}x_{m+1} \ 0 \ y_1 \ldots y_n \ x_1 \ldots x_m \ 0 \ $
$\stackrel{\tau_{2,n+m+4}}{\longmapsto} \ $

\smallskip

$x_{m+1}x_{m+1} \ 00 \ y_1 \ldots y_n \ x_1 \ldots x_m \ x_{m+1} \ $

\smallskip

\noindent Applying \ 
$\tau_{n+2,n+3} \ \ldots \ \tau_{3,4} \ \tau_{1,2}$ \ and \
$\tau_{n+1,n+2} \ \ldots \ \tau_{2,3} \ \tau_{0,1}$ \
we obtain

\smallskip

$0 \ 0 \ y_1 y_2 \ldots y_n \ x_{m+1} x_{m+1} \ x_1 \ldots x_m x_{m+1}$.

\medskip

\noindent Case $i(n) = 1$: We continue the simulation with 

\smallskip

$00 \ y_1 \ldots y_n \ x_1 \ldots x_m \ x_{m+1} \ $
$\stackrel{\varphi_{\rm 0f,4}}{\longmapsto} \ $
$00000 \ y_1 \ldots y_n \ x_1 \ldots x_m \ x_{m+1} \ $
$\stackrel{\tau_{2,n+m+5}}{\longmapsto} \ $

\smallskip

$00 \ x_{m+1} \ 00 \ y_1 \ldots y_n \ x_1 \ldots x_m \ 0 \ $
$\stackrel{\varphi_{\vee}}{\longmapsto} \ $
$x_{m+1} \ 0 \ x_{m+1} \ 00 \ y_1 \ldots y_n \ x_1 \ldots x_m \ 0 \ $
$\stackrel{\tau_{0,1}}{\longmapsto} \ $

\smallskip

$0 \ x_{m+1}x_{m+1} \ 00 \ y_1 \ldots y_n \ x_1 \ldots x_m \ 0 \ $
$\stackrel{\varphi_{\vee}}{\longmapsto} \ $
$x_{m+1}x_{m+1}x_{m+1} \ 00 \ y_1 \ldots y_n \ x_1 \ldots x_m \ 0 \ $
$\stackrel{\tau_{2,n+m+5}}{\longmapsto} \ $

\smallskip

$x_{m+1}x_{m+1} \ 000 \ y_1 \ldots y_n \ x_1 \ldots x_m \ x_{m+1} \ $

\smallskip

\noindent Applying \ 
$\tau_{n+3,n+4} \ \ldots \ \tau_{3,4} \ \tau_{1,2}$ \ and \
$\tau_{n+2,n+3} \ \ldots \ \tau_{2,3} \ \tau_{0,1}$ \
we obtain

\smallskip

$0 \ 00 \ y_1 y_2 \ldots y_n \ x_{m+1} x_{m+1} \ x_1 \ldots x_m x_{m+1}$.

\medskip

\noindent The above gives us a word $w_C$ of length \ 
$|w_C| \leq |w_K| + 4m + 6n +20$ \ over the alphabet 
$\Delta_{3,1}^{\pm 1} \cup \{\tau_{i,i+1} : 0 \leq i \}$, simulating $f_C$. 
Moreover, the subscripts of the transpositions appearing in $w_C$ are \ 
$\leq {\rm max}\{n+m+5, J_K\}$, where $J_K$ is the largest subscript in any
transposition appearing in $w_K$.

In case we want to change the positions of the added variables $x_{m+1}$,
$y_{n+1}$, and $y_{n+2}$ (so that $x_{m+1}$ is the $i$th input variable,
$y_{n+1}$ is the $j_1$th output variable, and $y_{n+2}$ is the $j_2$th output
variable), we apply appropriate other permutations (instead of \ 
$\tau_{n+3,n+4} \ \ldots \ \tau_{3,4} \ \tau_{1,2}$, 
$\tau_{n+2,n+3} \ \ldots \ \tau_{2,3} \ \tau_{0,1}$, and $\tau_{2,n+m+5}$).
This does not change our upper bounds on $|w_C|$ and $J_C$.  

\medskip

In each of the three cases, the circuit $C$ of depth 1 is simulated by a 
word \ $w_C$ over the alphabet 
$\Delta_{3,1}^{\pm 1} \cup \{\tau_{i,i+1} : 0 \leq i \}$,
of length $|w_C| \leq 9 \, |K| + 20 + |w_K|$. 
After $\leq |C|$ construction steps (starting with $K$ being the empty 
circuit, and ending with $K$ being $C$), the length of $w_C$ will be 
 \  $|w_C| \ \leq \ \frac{9}{2} \, |C|^2 + 25 |C|$.
The transpositions occurring in $w_C$ have maximum subscript $\leq |C| + 5$.  
The above construction of each word $w_C$ from $C$ is a
polynomial-time algorithm.

\bigskip

{\sc Inductive step:} Assume that $C$ has depth $L > 1$.  In order to define 
$w_C$ we can use the fact that we have already defined the words 
$w_{C_{\ell}}$ $(1 \leq \ell \leq L)$ for the slices $C_{\ell}$ of $C$.
Indeed, each slice has depth 1, so the base of the induction applies.
Each word $w_{C_{\ell}}$ has all the properties claimed in the Theorem for 
circuit $C_{\ell}$. In particular, $w_{C_{\ell}}$
defines the map 

\medskip

$\Phi_{C_{\ell}} : 0 \ Y^{\ell-1} \ \longmapsto \ $
$0 \ 0^{i(n_{\ell})} \ Y^{\ell} \ Y^{\ell-1}$.

\medskip

\noindent
Hence, since $\Phi_{C_{\ell}}$ is a right ideal isomorphism, we also have  

\medskip

$0 \ Y^{\ell-1} \ 0^{i(n_{\ell-1})} \ Y^{\ell-2} \ 0^{i(n_{\ell-2})} \ $
 $\ldots \ Y^1 \ 0^{i(n_1)} \ x_1 \ldots x_m$ 
 $  \ \stackrel{\Phi_{C_{\ell}}}{\longmapsto} \ $

\medskip

$0 \ 0^{i(n_{\ell})} \ Y^{\ell} \ $
  $Y^{\ell-1} \ 0^{i(n_{\ell-1})} \ Y^{\ell-2} \ 0^{i(n_{\ell-2})} $
 $ \ \ldots \ Y^1 \ 0^{i(n_1)} \ x_1 \ldots x_m \ $

\medskip

\noindent Applying \ $(\sigma_{1,n_{\ell}})^{i(n_{\ell})}$ \ to this word
yields

\medskip

$0 \ Y^{\ell} \ 0^{i(n_{\ell})} \ Y^{\ell-1} \ 0^{i(n_{\ell-1})} \ $
  $Y^{\ell-2} \ 0^{i(n_{\ell-2})} \ \ldots \ Y^1 \ 0^{i(n_1)}$
     $ \ x_1 \ldots x_m$.

\medskip

\noindent where, in general, $\sigma_{i,j}$ denotes the permutation \ 
$\tau_{j-1,j} \, \tau_{j-2,j-1} \ \ldots \ $
$\tau_{i+1,i+2} \, \tau_{i,i+1} (\cdot)$ \ (for all $0 \leq i < j$).  
Therefore, 

\medskip

$(\sigma_{1,n_L})^{i(n_L)} \ w_{C_L} \ $
$(\sigma_{1,n_{L-1}})^{i(n_{L-1})} \ w_{C_{L-1}} \ $
$ \ldots \ $
$(\sigma_{1,n_{\ell}})^{i(n_{\ell})} \ w_{C_{\ell}} \ $ 
$\ldots \ $
$(\sigma_{1,n_1})^{i(n_1)} \ w_{C_1}$ \

\medskip

\noindent defines the map 

\smallskip

 $0 x_1 \ldots x_m \ \ \longmapsto \ \ $

\smallskip

$0 \ y_1 \ldots y_n \ 0^{i(n)}$
$Y^{L-1} \ 0^{i(n_{L-1})} \ \ldots \ Y^{\ell} \ 0^{i(n_{\ell})} \ \ldots \ $
$Y^2 \ 0^{i(n_2)} \ Y^1 \ 0^{i(n_1)} \ x_1 \ldots x_m \ \ (=_{\rm def} \ Z).$ 

\smallskip

\noindent Note that the length of the word $Z$ is \  
$|Z| \leq 1 + |C| + 2L \leq 3 \cdot |C|$. Indeed, the total number of 
variables in the circuit (i.e., $n_L+ \ldots +n_1 + m$) is equal to the total 
number of wires (i.e., $|C|$); the ``$+1$'' comes from the leading letter $0$; 
the ``$2L$'' comes from $i(n)$, $i(n_{L-1})$, $\ldots$, $i(n_1)$. 
Recall that $y_1 \ldots y_n = Y^L$, and $n_L = n$. 

\smallskip

\noindent Now the permutation \  $\pi_1 \ = \ (\sigma_{1,|Z|})^{n + i(n)}$ \ 
transforms the word $Z$ into

\smallskip

$0 \ Y^{L-1} \ 0^{i(n_{L-1})} \ \ldots \ Y^{\ell} \ 0^{i(n_{\ell})}$
$ \ \ldots \ Y^2 \ 0^{i(n_2)} \ Y^1 \ 0^{i(n_1)} \ x_1 \ldots x_m$ \
$y_1 \ldots y_n \  0^{i(n)}$.

\smallskip

\noindent Note that the word length of $\pi_1$ is less than \ 
$(n+2) \cdot |Z| \leq 3(n+2) \cdot |C| \leq 3 \, |C|^2$  \ over the alphabet 
$\{\tau_{i,i+1} : 0 \leq i\}$.

\smallskip

\noindent Next (and this is a crucial idea in reversible computing), 
applying 

\medskip

$[(\sigma_{1,n_{L-1}})^{i(n_{L-1})} \ w_{C_{L-1}} \ $
$\ldots \ (\sigma_{1,n_{\ell}})^{i(n_{\ell})} \ w_{C_{\ell}} \ \ldots \ $
$ (\sigma_{1,n_2})^{i(n_2)} \ w_{C_2}\ $
$ (\sigma_{1,n_1})^{i(n_1)} \ w_{C_1}]^{-1}$ \

\medskip

\noindent yields  \ \ $0 \ x_1 \ldots x_m \ y_1 \ldots y_n \ 0^{i(n)}$.

\smallskip

\noindent Finally, applying the permutation \ 
$\pi_2 \ = \ (\sigma_{1,n+m+i(n)})^m$ \ 
produces the desired final output
 
\smallskip

 \ \ \ \ \ $0 \ 0^{i(n)} \ y_1 \ldots y_n \ x_1 \ldots x_m$.

\smallskip

\noindent Therefore we can define $w_C$ (over the alphabet 
$\Delta_{3,1}^{\pm 1} \cup \{\tau_{i,i+1} : 0 \leq i \}$) by

\bigskip

$w_C \ = \  $
  $ \pi_2 \  $
  $[(\sigma_{1,n_{L-1}})^{i(n_{L-1})} \ w_{C_{L-1}} \ \ldots \ $ 
        $ (\sigma_{1,n_1})^{i(n_1)} \ w_{C_1}]^{-1} \ $
 $ \pi_1 \ \cdot  $
 
\smallskip

\hspace{.55in}  
 $ (\sigma_{1,n_L})^{i(n_L)} \ w_{C_L} \ $ 
 $ (\sigma_{1,n_{L-1}})^{i(n_{L-1})} \ w_{C_{L-1}} \ \ldots \ $ 
 $ (\sigma_{1,n_1})^{i(n_1)} \ w_{C_1} $     

\bigskip

\noindent For the length we have therefore 

\medskip

$|w_C| \leq |\pi_2| + \sum_{\ell=1}^{L-1} |w_{C_{\ell}}| + $
       $\sum_{\ell = 1}^{L-1} i(n_{\ell}) \, |\sigma_{1,n_{\ell}}| + $
     $|\pi_1| +  \sum_{\ell=1}^L |w_{C_{\ell}}| + $
       $\sum_{\ell=1}^L i(n_{\ell}) \, |\sigma_{1,n_{\ell}}|$.

\medskip

\noindent Since \   
$|w_{C_{\ell}}| \leq \frac{9}{2} \, |C_{\ell}|^2 + 25 \, |C_{\ell}|$ 
  \  (for $1 \leq \ell \leq L$), 
and \ $\sum_{\ell = 1}^L |C_{\ell}| = |C|$, we have \ 
$\sum_{\ell=1}^L |C_{\ell}|^2 \leq |C|^2$. Also, $i(n_{\ell}) \leq 2$, 
and $|\sigma_{1,n_{\ell}}| \leq n_{\ell}$, so \   
$\sum_{\ell=1}^L i(n_{\ell}) \, |\sigma_{1,n_{\ell}}| \leq 2 \, |C|$. 
Thus \ $|w_C| \leq c \cdot |C|^2$, for some positive constant $c$.
Also, the largest subscript in any permutation is $\leq |Z| \leq  3 \,|C|$.
Since $|C|$ was squared in order to obtain strict layering, the above 
bounds become

\medskip

$|w_C| \leq c \, |C|^4$,  

\medskip

$J_C \leq 3 \, |C|^2$.

\medskip

\noindent
The word $w_C$ can be written down in linear time, based on the words 
$w_{C_{\ell}}$ ($1\leq \ell \leq L$), and we saw that each $w_{C_{\ell}}$
can be computed in polynomial time from $C_{\ell}$.  

\bigskip

In order to obtain a word that {\it strongly} simulates $f_C$ we need to 
make two additions to $w_C$: A pre-processing step $w_0$ is attached at the 
beginning (the right side) of $w_C$, to make sure inputs that are 
``too short'' are handled correctly. 
A post-processing step $w_{L+1}$ is attached at the end (the left side) of 
$w_C$, in order to remove excess letters introduced during pre-processing. 
Recall that we write functions to the left of the argument.
The word that strongly simulates $f_C$ is denoted by $W_C$ and defined by 
         \[ W_C = w_{L+1} \, w_C \, w_0  \] 
We define $w_0$ by 

\medskip 

$w_0 \ = \ $
$\tau_{1,3(n+m)+1} \ \ldots \ $
$\tau_{j,3(n+m)+j} \ \ldots \ \tau_{m,3(n+m)+m} \ $
$ (\varphi_{\rm 0f,4})^{n+m}(\cdot)$ 

\medskip 
  
\noindent So, $|w_0|$ is bounded from above by a quadratic function in 
$n+m$, and $J_{w_0}$ is linearly bounded in $n+m$.  We have 

\smallskip

$0 \ x_1 \ldots x_m  \ \stackrel{w_0}{\longmapsto} \ $
$0 \ x_1 \ldots x_m \ 0^{3(n+m)}$
$\stackrel{w_C}{\longmapsto} \ $
$0^{1+i(n)} \ f_C(x_1, \ldots, x_m) \ x_1 \ldots x_m \ 0^{3(n+m)}$ . 

\medskip

\noindent For $0 \leq k < m$ we have on an input that is ``too short'': 

\smallskip

$0 \ x_1 \ldots x_k \ \# \ \stackrel{w_0}{\longmapsto} \ $
$0 \ z_1 \ldots z_{k+3(n+m)} \ \#$ 
$\stackrel{w_C}{\longmapsto} \ $
$0^{1+i(n)} \ f_C(z_1, \ldots, z_m) \ z_1 \ldots z_{k+3(n+m)} \ \#$,

\smallskip

\noindent where $z_1 \ldots z_{k+3(n+m)}$ is a permuted version of
$x_1 \ldots x_k \ 0^{3(n+m)}$; this permutation depends only on
the number $k+3(n+m)$. So the outcome \ 
$0^{1+i(n)} \ f_C(z_1, \ldots, z_m) \ z_1 \ldots z_{k+3(n+m)} \#$ \ does 
not depend on the  circuit $C$ that was used to implement the function $f_C$. 

\medskip

\noindent Finally, it is also easy to verify that 

\smallskip

$ \# \ \stackrel{w_0}{\longmapsto} \ 0^{3(n+m)} \ \#$
$\stackrel{w_C}{\longmapsto} \ $
$0 \ f_C(0, \ldots, 0) \ 0^{3(n+m)-1} \ \#$ .

\bigskip
 
\noindent We define $w_{L+1}$ by 

\medskip

$w_{L+1} \ = \ $
$(\varphi_{\rm 0f,4})^{-n-m} \ $
$\tau_{n+m,3(n+m)+n+m} \ \ldots \ \tau_{j,3(n+m)+j} \ $
$\ldots \ \tau_{1,3(n+m)+1} (\cdot)$. 

\medskip

\noindent So, $|w_{L+1}|$ is bounded from above by a quadratic function in
$n+m$, and $J_{w_{L+1}}$ is linearly bounded in $n+m$. 
One can verify easily that 

\smallskip

$0^{1+i(n)} y_1 \ldots y_n x_1 \ldots x_m 0^{3(n+m)} \ $
$\stackrel{w_{L+1}}{\longmapsto} \ $
$0^{1+i(n)} y_1 \ldots y_n x_1 \ldots x_m$. 

\medskip

\noindent For $0 \leq k < m$, and $x_1, \ldots, x_k \in \{0,1\}$, 
let $z_1 \ldots z_{k+3(n+m)}$ be the permuted version of
$x_1 \ldots x_k \ 0^{3(n+m)}$ considered above. Let \  
$f_C(z_1 \ldots z_m) = y_1 \ldots y_n$; note that this string does not depend
on the circuit $C$ that was used to implement the function $f_C$.

Then the sequence of transformations $w_{L+1}$ will be applied to \   
$0 \ y_1 \ldots y_n \ z_1 \ldots z_{k+3(n+m)} \ \#$. 
This will produce a new string ($\in \{0,1\}^* \#$) which does not
depend on the circuit $C$ that was used to implement the function $f_C$.

\medskip 

\noindent Also, recall that on argument $\#$, the outcome of the sequence 
of transformations $w_0 w_C$ is  \ $0 \ y_1 \ldots y_n \ 0^{3(n+m)-1} \#$, 
where  \ $f_C(0, \ldots, 0) = y_1 \ldots y_n$.  Then, applying $w_{L+1}$
yields a string ($\in \{0,1\}^* \#$) which does not
depend on the circuit $C$ that was used to implement the function $f_C$.
 \ \ \ $\Box$ 

\bigskip

\noindent {\bf Remarks}: \ The length of $w_C$ (over the infinite alphabet
$\Delta_{3,1}^{\pm 1} \cup \{\tau_{i,i+1} : 0 \leq i \}$), and the largest 
subscript $J_{w_C}$ (in any transposition occurring in $w_C$) are bounded 
from above by polynomials in $|C|$. 
Hence, if we write subscripts of transpositions in unary notation, the length 
of $w_C$ remains bounded from above by a polynomial in $|C|$.  

The group $G_{3,1}$ is finitely generated, so one may wonder what 
the word length of $w_C$ would be if $w_C$ were expressed over such a finite 
generating set; we will see that it is exponential (Lemma 
\ref{distortionVinTboolLoweBound}, Theorem \ref{distortionVinTbool}).

\medskip

In the next section we reduce the above problem to a certain generalized word 
problem of $G_{3,1}$, still over the infinite generating set 
$\Delta_{3,1} \cup \{\tau_{i,i+1} : 0 \leq i \}$.

%%%%%%%%%%%%%%%%%%%%%%%%%%%%%%%%%%%%%%%%%%%%%%%%%%%%%%%%%%%%%%%%%%%
% Section 4
%%%%%%%%%%%%%%%%%%%%%%%%%%%%%%%%%%%%%%%%%%%%%%%%%%%%%%%%

\section{Reduction to a generalized word problem \\   
  (over an infinite generating set)}

We will now restate the above reduction as a reduction to a generalized word
problem of a Thompson group, over an infinite generating set. 
In the following definitions we represent elements of ${\mathcal G}_{3,1}$
by right ideal isomorphisms between essential right ideals of $\{0,1,\#\}^*$.
We will extend the classical concepts of stabilizers and fixators to the case
of partial permutations. 

\medskip

\noindent {\bf Definitions.} \  
We say that $g$ {\em partially stabilizes} a set of words
$S \subseteq \{0,1,\#\}^*$ \ iff \ $g(S) \cup g^{-1}(S) \subseteq S$.
So $g$ maps $S$ into itself wherever $g$ is defined, and
similarly for $g^{-1}$.
For a subgroup $G \subseteq {\mathcal G}_{3,1}$, the {\bf partial stabilizer} 
(in $G$) of $S$ is

\medskip

 \ \ \ \ \  ${\rm pStab}_G(S) \ = \ $
  $\{ g \in G \ : \ g(S) \ \cup \ g^{-1}(S) \ \subseteq S \}.$

\medskip

\noindent We say that $g$ {\em totally stabilizes} a set of words $S$ \ iff \
$g(S) \cup g^{-1}(S) \subseteq S$, and  in addition,
$S \subseteq {\rm Dom}(g) \cap {\rm Im}(g)$. So $g$ totally stabilizes $S$ 
iff $g$ partially stabilizes $S$ and moreover, $g$ and $g^{-1}$ are defined 
everywhere on $S$.
For a subgroup $G \subseteq {\mathcal G}_{3,1}$, the {\bf total stabilizer} 
(in $G$) of $S$ is

\medskip

 \ \ \ \ \  ${\rm tStab}_G(S) \ = \ $
  $\{ g \in G \ : \ g(S) \cup g^{-1}(S) \subseteq S$
    $\subseteq {\rm Dom}(g) \cap {\rm Im}(g) \}.$

\medskip

\noindent We say that $g$ {\em partially fixes} a set $S$ \ iff \
$g(x) = x$ \ for every $x \in S \cap {\rm Dom}(g) \cap {\rm Im}(g)$;
this is also called partial ``pointwise stabilization''.
For $G \subseteq {\mathcal G}_{3,1}$, the {\bf partial fixator} (in $G$) 
of $S$ is

\medskip

 \ \ \ \ \  ${\rm pFix}_G(S) \ = \ \{ g \in G : \ $
$(\forall x \in S \cap {\rm Dom}(g) \cap {\rm Im}(g)) \ \ g(x) = x \}$

\medskip

\noindent i.e., the elements $g$ of $G$ that fix every point in $S$ on which
$g$ and $g^{-1}$ are defined.
We can also define the {\bf total fixator} by

\medskip

 \ \ \ \ \  ${\rm tFix}_G(S) \ = \ \{ g \in G : \ $
$ S \subseteq {\rm Dom}(g) \cap {\rm Im}(g) $ \ and \
$(\forall x \in S)\ g(x) = x \}$,

\medskip

\noindent i.e., the elements $g$ of $G$ that fix every point in $S$ and 
such that $g$ and $g^{-1}$ are defined on every point of $S$. 
This completes the definitions of stabilizers and fixators.

\bigskip

Observe that when $R \subseteq \{0,1,\#\}^*$ is a right ideal generated by 
a maximal prefix code $P$ (over the alphabet $\{0,1,\#\}$), then

\medskip

 \ \ \ \ \  ${\rm pFix}_G(R) \ = \ {\rm tFix}_G(R) \ = \ {\rm tFix}_G(P)$.

\medskip

\noindent So, for right ideals, the notions of partial fixator and total
fixator coincide.
Moreover, for every right ideal $S \subseteq R$ such that $S$ is
{\it essential in $R$} (i.e., $S$ has a non-empty intersection with every
right ideal contained in $R$), we have:

\medskip

 \ \ \ \ \  ${\rm pFix}_G(S) \ = \ {\rm pFix}_G(R)$.

\medskip

\noindent It is easy to see that ${\rm tStab}_G(X)$ and ${\rm tFix}_G(X)$ 
are always groups (for any group $G \subseteq {\mathcal G}_{3,1}$ and any set 
$X \subseteq \{0,1,\#\}^*$). However, ${\rm pStab}_G(X)$ and 
${\rm pFix}_G(X)$ are not always groups.
For this paper, all we need is the next Lemma.
\begin{lem} \label{pStabpFixGroups} \  
Let $G \subseteq {\mathcal G}_{3,1}$. For any set $X$ of words over 
$\{0,1,\#\}$,
 \ ${\rm tStab}_G(X)$ and ${\rm tFix}_G(X)$ \ are subgroups of $G$.   
For any right ideal $R$ of $\{0,1,\#\}^*$, \ ${\rm pFix}_G(R)$ is a subgroup 
of $G$.    

If $G \subseteq G_{3,1}$ and if $B^*$ is any free submonoid of 
$\{0,1,\#\}^*$ (generated as a submonoid by a set of words 
$B \subseteq \{0,1,\#\}^*$), then ${\rm pStab}_G(B^*)$ and 
${\rm pFix}_G(B^*)$ are subgroups of $G$.
\end{lem}
{\bf Proof.} \ The sets ${\rm tStab}_G(X)$, ${\rm tFix}_G(X)$,
${\rm pFix}_G(R)$, and ${\rm pStab}_G(B^*)$ are closed under inverse, by
definition. The closure under multiplication is obvious for
${\rm tStab}_G(X)$ and ${\rm tFix}_G(X)$. And when $R$ is a right
ideal we saw that ${\rm pFix}_G(R) = {\rm tFix}_G(R)$

If $x \in B^*$ and $\varphi_2, \, \varphi_1 \in {\rm pStab}_G(B^*)$, and 
if $({\rm max} \varphi_2 \varphi_1)(x)$ is defined,  we need to show that 
$({\rm max} \varphi_2 \varphi_1)(x) \in B^*$.
Note that $\varphi_1(x)$ and $\varphi_2 \varphi_1(x)$ might be undefined; 
but in any case, there exists $w \in B^*$ such that
$\varphi_2 \varphi_1(xw)$ is defined; we just need to take $w$ long
enough. Then we also have $\varphi_2 \varphi_1(xw) \in B^*$ and
$\varphi_2 \varphi_1(xw) = ({\rm max} \varphi_2 \varphi_1)(xw)$
$= ({\rm max} \varphi_2 \varphi_1)(x) \cdot w$. Therefore, since
$w$ and $({\rm max} \varphi_2 \varphi_1)(x) \cdot w$ belong to $B^*$, and
since $B^*$ is free, we conclude that 
$({\rm max} \varphi_2 \varphi_1)(x) \in B^*$. The proof for 
${\rm pFix}_G(B^*)$ is very similar.
 \ \ \  $\Box$

\medskip 

\noindent With this terminology we can restate Lemma \ref{simulatVSfix}:
\begin{lem} \label{simulatVSfixators} \   
Let $f$ and $g$ be any boolean functions such that $f$ and $g$ have the same 
number of input variables, and $f$ and $g$ have the same number of output 
variables.
If $f$ and $g$ are simulated by $\Phi_f$, respectively $\Phi_g$, then the
following are equivalent:

\smallskip

\noindent $\bullet$ \ \ \  $f = g$

\smallskip

\noindent $\bullet$ \ \ \
          $\Phi_f^{-1} \, \Phi_g \in {\rm pFix}_{G_{3,1}}(0\{0,1,\# \}^*)$

\smallskip

\noindent $\bullet$ \ \ \  $\Phi_f^{-1} \, \Phi_g \in $
${\rm pFix}_{G_{3,1}^{\rm mod \, 3}}(0\{0,1,\#\}^*)$

\medskip

\noindent In the case of strong simulation the following are equivalent:

\smallskip

\noindent $\bullet$ \ \ \  $f = g$

\smallskip

\noindent $\bullet$ \ \ \
 $\Phi_f^{-1} \, \Phi_g \in {\rm pFix}_{G_{3,1}}(\{0,\#\} \{0,1,\# \}^*)$

\smallskip

\noindent $\bullet$ \ \ \ $\Phi_f^{-1} \, \Phi_g \in $
    ${\rm pFix}_{G_{3,1}^{\rm mod \, 3}}(\{0,\#\} \{0,1,\# \}^*)$
\end{lem}

Theorem \ref{reduction} and Lemma \ref{simulatVSfixators} give a 
polynomial-time one-to-one reduction from the equivalence problem for acyclic 
circuits to the generalized word problem of
${\rm pFix}_{G_{3,1}^{\rm mod \, 3}}(0\{0,1,\#\}^*)$ in $G_{3,1}$, with
elements of $G_{3,1}$ written over the set of generators
$\Delta_{3,1} \cup \{\tau_{i,i+1} : 0 \leq i \}$ \ (where $\Delta_{3,1}$ is 
a finite generating set of $G_{3,1}$).
It follows that this generalized word problem is coNP-hard.
Because of the existence of a strong simulation, we also have a
polynomial-time one-to-one reduction from the equivalence problem for acyclic
circuits (with last output variable 0 when the inputs are all 0) to the
generalized word problem of
${\rm pFix}_{G_{3,1}^{\rm mod \, 3}}(\{0,\#\}\{0,1,\#\}^*)$ in $G_{3,1}$
over the set of generators $\Delta_{3,1} \cup \{\tau_{i,i+1} : 0 \leq i \}$.
Hence we have:

\begin{cor}  {\bf (co-NP hard generalized word problem).} \
The generalized word problems of \ 
 ${\rm pFix}_{G_{3,1}^{\rm mod \, 3}}(0\{0,1,\#\}^*)$ \ and of \ 
 ${\rm pFix}_{G_{3,1}^{\rm mod \, 3}}(\{0,\#\} \{0,1,\#\}^*)$, as subgroups 
of $G_{3,1}$ (with generating set \  
$\Delta_{3,1} \cup \{\tau_{i,i+1} : 0 \leq i \}$) \ are coNP-hard.
\end{cor}

\noindent The following subgroups of $G_{3,1}$ will play a major role.

\begin{defn} \label{notation_G01} \   
The groups of {\em bit-preserving} (or $\{0,1\}$-preserving) transformations, 
$G_{3,1}(0,1)$ and $G_{3,1}^{\rm mod \, 3}(0,1)$, are defined by

\bigskip

\noindent 
$G_{3,1}(0,1) \ = \ {\rm pStab}_{G_{3,1}}(\{0,1\}^*)$

\smallskip

\hspace{.4in} $= \ \{ \phi \in G_{3,1} \ : \ $
   $ \phi(\{0,1\}^*) \subseteq \{0,1\}^*$ \ {\rm and} \    
   $\phi^{-1}(\{0,1\}^*) \subseteq \{0,1\}^* \, \}$, 

\medskip

\noindent $G_{3,1}^{\rm mod \, 3}(0,1) \ = \ $
${\rm pStab}_{G_{3,1}^{\rm mod \, 3}}(\{0,1\}^*)$

\smallskip

\hspace{.55in} 
 $= \ \{ \phi \in G_{3,1}(0,1) \ : \ $

\hspace{1in}
$|\phi(x)| \equiv |x| \ {\rm mod} \ 3,$
 {\rm for all} $x \in \{0,1\}^*$ {\rm for which} $\phi(x)$  
 {\rm is defined}$\}$. 

\bigskip

\noindent The groups of \, {\em \#-preserving} transformations, 
$G_{3,1}(0,1;\#)$ and $G_{3,1}^{\rm mod \, 3}(0,1;\#)$, are defined by

\bigskip

\noindent
$G_{3,1}(0,1;\#) \ = \ {\rm pStab}_{G_{3,1}}(\{0,1\}^*) \ $
   $\cap \ \, {\rm tStab}_{G_{3,1}}(\{0,1\}^* \, \#)$.

\smallskip

\hspace{.6in}
 $= \ \{ \phi \in G_{3,1}(0,1) \ : \ $ 
 {\rm for all} $x \in \{0,1\}^*$,

\hspace{1in}
 $\phi(x\#)$ {\rm and} $\phi^{-1}(x\#)$ {\rm are defined and} 
 $\phi(x\#), \, \phi^{-1}(x\#) \, \in \,  \{0,1\}^* \#  \, \}$, 

\medskip

\noindent $G_{3,1}^{\rm mod \, 3}(0,1;\#) \ = \ $
${\rm pStab}_{G_{3,1}^{\rm mod \, 3}}(\{0,1\}^*) \ $
   $ \cap \ \, {\rm tStab}_{G_{3,1}^{\rm mod \, 3}}(\{0,1\}^* \, \#)$

\smallskip

\hspace{.7in}
$= \ \{ \phi \in G_{3,1}(0,1;\#) \ : \ $

\hspace{1in}
$ |\phi(x)| \equiv |x| \ {\rm mod} \ 3$
 {\rm for all} $x \in \{0,1\}^*$ {\rm for which} $\phi(x)$ 
{\rm is defined}$\}$.
\end{defn}
It follows from Lemma \ref{pStabpFixGroups} that $G_{3,1}(0,1)$, 
$G_{3,1}(0,1;\#)$, $G_{3,1}^{\rm mod \, 3}(0,1)$, and 
$G_{3,1}^{\rm mod \, 3}(0,1;\#)$ 
are indeed groups.

All the elements of $G_{3,1}$ that we have used in the proof of Theorem
\ref{reduction} are generated by $\varphi_{\neg}$, $\varphi_{\vee}$, 
$\varphi_{\wedge}$, $\tau_{i,j}$ $(0 \leq i \leq j)$, and 
$\varphi_{\rm 0f,4}$. These elements also belong to 
$G_{3,1}^{\rm mod \, 3}(0,1;\#)$ 
($\subset G_{3,1}^{\rm mod \, 3}(0,1)$). 
Hence, the above Corollary implies the following, where $\Delta_{\#}$ is a 
finite generating set of $G_{3,1}^{\rm mod \, 3}(0,1;\#)$, and 
$\Delta_{(0,1)}$ is a finite generating set of $G_{3,1}^{\rm mod \, 3}(0,1)$:

\begin{cor}  {\bf (co-NP hard generalized word problem).} \
The generalized word problems of \, ${\rm pFix}_G(0 \, \{0,1,\#\}^*)$
and of \, ${\rm pFix}_G(\{0,\#\} \{0,1,\#\}^*)$ as subgroups of
$G = G_{3,1}^{\rm mod \, 3}(0,1;\#)$ are coNP-hard.
Here the generating set used for $G_{3,1}^{\rm mod \, 3}(0,1;\#)$ is 
$\Delta_{\#} \cup \{\tau_{i,i+1} : 0 \leq i \}$.

The generalized word problems of \, ${\rm pFix}_G(0 \, \{0,1,\#\}^*)$
and of \, ${\rm pFix}_G(\{0,\#\} \{0,1,\#\}^*)$, as subgroups of
$G = G_{3,1}^{\rm mod \, 3}(0,1)$ are coNP-hard.
Here the generating set used for $G_{3,1}^{\rm mod \, 3}(0,1)$ is 
$\Delta_{(0,1)} \cup \{\tau_{i,i+1} : 0 \leq i \}$.
\end{cor}
We will see later that $G_{3,1}^{\rm mod \, 3}(0,1;\#)$ and 
$G_{3,1}^{\rm mod \, 3}(0,1)$ are finitely presented, so the finite 
generating sets $\Delta_{\#}$ and $\Delta_{(0,1)}$ exist.

\bigskip

\noindent Here is a more concrete view of the subgroup 
$G_{3,1}^{\rm mod \, 3}(0,1;\#)$:
\begin{lem} \label{tableForm} \  
The group $G_{3,1}^{\rm mod \, 3}(0,1;\#)$ consists of the elements of 
$G_{3,1}$ that have tables of the form 
\[
\left[ \begin{array}{ccc ccc}
x_1 & \ldots & x_n & x_1'\# & \ldots & x_m'\#  \\
y_1 & \ldots & y_n & y_1'\# & \ldots & y_m'\#
\end{array}  \right],  
\]
for some positive integers $n, m$, with $x_1$, $\ldots$, $x_n$,
$x_1'$, $\ldots$, $x_m'$, $y_1$, $\ldots$, $y_n$, $y_1'$, $\ldots$, $y_m'$
$\in \{0,1\}^*$, and $|x_i| \equiv |y_i|$ mod 3 (for all $i = 1, \ldots, n$).
Moreover,  
$\{ x_1, \ldots, x_n \}$ $\cup$ $\{ x_1', \ldots, x_n'\} \#$ \ and \   
$\{ y_1, \ldots, y_n \}$ $\cup$ $\{ y_1', \ldots, y_n'\} \#$ \ are
maximal prefix codes over $\{0,1,\#\}$. 
\end{lem}
{\bf Proof.} \ From the shape of the above table we see immediately that 
the corresponding element $\phi$ of $G_{3,1}$, as well as $\phi^{-1}$,
map $\{0,1\}^*$ into $\{0,1\}^*$, and $\{0,1\}^*\#$ into $\{0,1\}^*\#$.
On $\{0,1\}^*$, $\phi$ preserves length modulo 3. Thus, $\phi \in$
pStab$_{G_{3,1}^{\rm mod \, 3}}(\{0,1\}^*)$.
Moreover, since \ $\{x_1, \ldots, x_n \}$ $\cup$ $\{ x_1', \ldots, x_n'\} \#$
 \ is a maximal prefix code, $\phi(w\#)$ is defined for all $w \in \{0,1\}^*$.
Similarly, $\phi^{-1}(w\#)$ is always defined. Thus, $\phi \in$
tStab$_{G_{3,1}^{\rm mod \, 3}}(\{0,1\}^*\#)$. 

Conversely, if $\phi \in G_{3,1}^{\rm mod \, 3}(0,1;\#)$ then the domain code 
of $\phi$ is a subset of \ $\{0,1\}^* \cup \{0,1\}^* \#$ (since $\phi$ 
partially stabilizes $\{0,1\}^*$ and totally stabilizes $\{0,1\}^* \#$). For 
the same reason, the image code of $\phi$ is a subset of \ 
 $\{0,1\}^* \cup \{0,1\}^* \#$.  By Lemma \ref{endmarker_code}, and the 
definition of $G_{3,1}^{\rm mod \, 3}(0,1;\#)$
it now follows immediately that $\phi$ has a table of the above form.  
 \ \ \ $\Box$

\begin{lem}
\label{endmarker_code} \ {\bf (1)} \ 
If \ $P \subset \{0,1\}^* \cup \{0,1\}^* \#$ \ is a maximal prefix code 
over $\{0,1,\# \}$ then \ $P = P_1 \, \cup \, P_2 \#$ \ for some
$P_1, P_2 \subset \{0,1\}^*$, with the following properties:

\smallskip

$\bullet$  \ $P_1$ is a maximal prefix code over $\{0,1\}$;

\smallskip

$\bullet$ \ $P_2 \ = \ \{ p \in \{0,1\}^* \ : \ $
       $p$ is a strict prefix of some element of $P_1  \}$.

\smallskip

\noindent When $P_1$ is finite, this last property implies: \
 $|P_2| = |P_1| - 1$.

\smallskip

\noindent {\bf (2)} \ 
 Conversely, if \ $P \ = \ P_1 \ \cup \ P_2 \#$ \ for some
$P_1, P_2 \subset \{0,1\}^*$ with the above two properties, then $P$ is a
maximal prefix code over $\{0,1,\# \}$.
\end{lem}
{\bf Proof.}  The proof is not difficult and appears in the Appendix.
 \ \ \ $\Box$

\bigskip

\noindent Similarly, $G_{3,1}^{\rm mod \, 3}(0,1)$ has a concrete description.

\begin{lem} \label{tableForm01} \
The group $G_{3,1}^{\rm mod \, 3}(0,1)$ consists of the elements of $G_{3,1}$ 
that have a table of the form
\[  \left[ \begin{array}{ccc ccc}
x_1 & \ldots & x_n & x_1'\#s_1 & \ldots & x_m'\#s_m  \\
y_1 & \ldots & y_n & y_1'\#t_m & \ldots & y_m'\#t_m
\end{array}        \right],
\]
\noindent for some positive integers $n, m$, with $x_1$, $\ldots$, $x_n$,
$x_1'$, $\ldots$, $x_m'$, $y_1$, $\ldots$, $y_n$, $y_1'$, $\ldots$, $y_m'$
$\in \{0,1\}^*$, and $s_1$, $\ldots$, $s_m$, $t_1$,  $\ldots$, $t_m$ 
$\in \{0,1,\#\}^*$, and $|x_i| \equiv |y_i|$ mod 3 \
for all $i = 1, \ldots, n$.
Moreover,
$\{ x_1, \ldots, x_n \}$ $\cup$ $\{ x_1'\#s_1, \ldots, x_n'\#s_m \}$ \ and 
 \ $\{ y_1, \ldots, y_n \}$ $\cup$ $\{ y_1'\#t_1, \ldots, y_n'\#t_m \}$
 \ are maximal prefix codes over $\{0,1,\#\}$.
\end{lem}
{\bf Proof.} \ The proof is similar to the proof of the corresponding Lemma
for $G_{3,1}^{\rm mod \, 3}(0,1;\#)$.  \ \ \  $\Box$

\bigskip

In the next section we will reduce the above generalized word problems to the 
word problem of $G_{3,1}$ (still over the infinite generating set
$\Delta_{3,1} \cup \{\tau_{i,i+1} : 0 \leq i \}$).

%%%%%%%%%%%%%%%%%%%%%%%%%%%%%%%%%%%%%%%%%%%%%%%%%%%%%%%%%%%%%%%%%%%
% Section 5
%%%%%%%%%%%%%%%%%%%%%%%%%%%%%%%%%%%%%%%%%%%%%%%%%%%%%%%%

\section{Reduction to the word problem of a Thompson group \\  
  (over an infinite generating set)}

We will give a linear-time $k$-ary conjunctive reduction (for a constant $k$)
from the generalized word problem of 
${\rm pFix}_{G_{3,1}^{\rm mod \, 3}}(\{0, \#\}\{0,1,\#\}^*)$ to 
the word problem of $G_{3,1}$, over the infinite generating set \ 
$\Delta_{3,1} \cup \{\tau_{i,i+1} : 0 \leq i \}$.

\begin{defn} \label{conj_red_def} \  
A polynomial-time $k$-ary {\bf conjunctive reduction} from a language 
$L \subseteq A^*$ to a language $W \subseteq B^*$ is a function 
$f: A^* \ \to \ (B^*)^k$ \ such that $f(x)$ is computable in time bounded 
by a polynomial in $|x|$, and such that we have:

\smallskip

$x \in L$ \ \ iff \ \ $f(x) = (y_1, \ldots, y_k)$ \ with $y_i \in W$ 
for all $i$ $(1 \leq i \leq k)$.

\smallskip

\noindent
Any polynomial-time $k$-ary conjunctive reduction, for some constant $k$,
is called polynomial-time constant-arity conjunctive reduction.
\end{defn}

\noindent The conjunctive reductions used in this paper will have a 
constant arity. More general definitions of 
polynomial-time conjunctive reductions are possible (where the arity $k$ 
is a polynomial function of $|x|$), but we will not need this here. 
Conjunctive reductions are a special case of truth-table reductions. 
Note that the classes P, NP, and coNP are closed under polynomial-time 
constant-arity conjunctive reduction.

\medskip

In the classical theory of permutation groups there are many results of the
following form: Let $G$ be a permutation group acting on a set $X$ (i.e., 
$G \subseteq {\mathfrak S}_X$), and let $Q_1$, $Q_2$ be two ``complementary''
subsets of $X$. Then for all $g \in G$ we have: 

\medskip

\noindent (C) \hspace{1.5in}     $g \in {\rm Fix}_G(Q_1)$
 \ iff \ $gh = hg$ for all $h \in {\rm Fix}_G(Q_2)$. 

\medskip

\noindent We call property (C) the {\it commutation test} for 
the generalized word problem of ${\rm Fix}_G(Q_1)$. 
The left-to-right implication is obvious. For the right-to-left
implication to be true, special assumptions have to be made on $G$, on 
its action (i.e., on the embedding $G \hookrightarrow {\mathfrak S}_X$), and 
on the meaning of ``complementary''.  

What is interesting about the commutation test (C) is that it reduces 
the generalized word problem of ${\rm Fix}_G(Q_1)$ (as a subgroup of $G$) 
to $N$ instances of the word problem of $G$, where $N$ is the minimum number 
of generators of ${\rm Fix}_G(Q_2)$; indeed, $g$ commutes with all elements 
$h$ in Fix$_G(Q_1)$ iff $g$ commutes with all the members of a generating 
set of Fix$_G(Q_1)$. So, if ${\rm Fix}_G(Q_2)$ is finitely generated then we 
obtain a constant-arity conjunctive reduction of the generalized word 
problem of ${\rm Fix}_G(Q_1)$ to the word problem of $G$. 

In this Section we prove our version of the commutation test, namely 
Theorem \ref{commutation_test} below. Since we deal with partial actions 
(Thompson groups), everything is somewhat different from the classical case. 
We first introduce some concepts about prefix codes and fixators. 

We make the following convention: Let $\phi: A^* \to A^*$ be a partial
map and $x \in A^*$; when we write $\phi(x)$ it is to be understood that
$\phi(x)$ is defined (i.e., $x \in$ Dom$(\phi)$).

\begin{defn} \   
Let $A$ be a finite alphabet. 
Two prefix codes $P, P' \subset A^*$ are {\em complementary prefix codes} 
iff $P \cup P'$ is a maximal prefix code over $A$, and
 \ $P A^* \ \cap \ P' A^* = \emptyset$.
\end{defn}

\begin{defn} \ 
Let $A$ be a finite alphabet with $|A| = n$, and let 
$G \subseteq {\mathcal G}_{n,1}$ (i.e., $G$ is a subgroup with a particular 
embedding). The fixator ${\rm pFix}_G(P'A^*)$ is called {\em maximal} iff 
there exists $P \subset A^*$ such that $P, P'$ are complementary prefix 
codes, and such that we have:

\smallskip

for all $x \in PA^*$ there is $h \in {\rm pFix}_G(P'A^*)$ such that
$h(x) \neq x$.  

\smallskip

\noindent Equivalently: The fixator of a right ideal $P'A^*$
is maximal iff it does not fix any larger right ideal than $P'A^*$.
\end{defn} 
Recall our convention that when we write $\phi(x)$ (for a partial map $\phi$)
it is to be understood that $\phi(x)$ is defined (i.e., $x \in$ Dom$(\phi)$).

\medskip

\noindent In analogy with $G_{3,1}^{\rm mod \, 3}(0,1)$ and
$G_{3,1}^{\rm mod \, 3}(0,1;\#)$ we use the notation

\smallskip

${\mathcal G}_{3,1}^{\rm mod \, 3}(0,1) \ = \ $
${\rm pStab}_{{\mathcal G}_{3,1}^{\rm mod \, 3}}(\{0,1\}^*)$,

\smallskip

${\mathcal G}_{3,1}^{\rm mod \, 3}(0,1;\#) \ = \ $
${\rm pStab}_{{\mathcal G}_{3,1}^{\rm mod \, 3}}(\{0,1\}^*) \ \cap \ $
${\rm tStab}_{{\mathcal G}_{3,1}^{\rm mod \, 3}}(\{0,1\}^*\#)$.

\begin{defn} \label{def_max_sep_endm} \
Let \ $G \subset {\mathcal G}_{3,1}^{\rm mod \, 3}(0,1;\#)$ be a group.
Let $P, P'$ be complementary prefix codes over $\{0,1,\#\}$,  with
$P \cap \{0,1\}^* \neq \emptyset$, $P' \cap \{0,1\}^* \neq \emptyset$,
and \  $P, P' \subset \{0,1\}^* \cup \{0,1\}^*\#$. So,
$P = P_1 \cup P_2\#$, and $P' = P'_1 \cup P'_2\#$, according to Lemma
\ref{endmarker_code}.

\smallskip

\noindent The fixator ${\rm pFix}_G(P'\{0,1,\#\}^*)$ is {\em separating} 
on $P \{0,1,\#\}^*$ iff the following hold:

\smallskip

\noindent $\bullet$ \ For any ordered pair of prefix-incomparable words  
$(x,y)$ with \ $x, y \in P_1 \{0,1\}^*$, there exists 
$h \in {\rm pFix}_G(P' \{0,1,\#\}^*)$ and there exists $u \in \{0,1\}^*$ 
such that

\smallskip

 \ \ \ \ \    $h(xu) = xu$ \ and \ $h(yu) \neq yu$.

\smallskip

\noindent $\bullet$ \ For any ordered pair of prefix-incomparable words  
$(x,y)$ with \ $x, y \in P_1 \{0,1\}^*\# \ \cup \ P_2 \#$ \ there exists 
$h \in {\rm pFix}_G(P' \{0,1,\#\}^*)$ such that

\smallskip

 \ \ \ \ \   $h(x) = x$ \ and \ $h(y) \neq y$.
\end{defn}
We will not need any explicit separation requirements in the case
where $x \in \{0,1\}^*$ and $y \notin \{0,1\}^*$, or the case where
$x \notin \{0,1\}^*$ and $y \in \{0,1\}^*$. Also, note that for words
$x, y \in \{0,1\}^* \#$, $x,y$ are prefix-incomparable iff $x \neq y$.

\begin{thm} \label{commutation_test} 
{\bf (Commutation test for} ${\rm pFix}_G(0 \, \{0,1,\#\}^*)$ {\bf).} \ 
Let \ $G \ = \ G_{3,1}^{\rm mod \, 3}(0,1;\#)$. Then for any $g \in G$ we 
have:

\smallskip

$g \in {\rm pFix}_G(0 \, \{0,1,\#\}^*)$ \ \ iff \ \ $gh = hg$
 \ for all $h \in {\rm pFix}_G( \{1,\# \}\{0,1,\#\}^*)$.
\end{thm}
This Theorem follows immediately from the following two Propositions, 
\ref{max_sep_end_comm} and \ref{G_31_endm_sep}.

\begin{pro} \label{max_sep_end_comm} \
Suppose ${\rm pFix}_G(P'\{0,1,\#\}^*)$ is {\em separating} on $P \{0,1,\#\}^*$, 
where $G$, $P$, and $P'$ are as in Definition \ref{def_max_sep_endm}.
Then for all $g \in G$ we have:

\smallskip

If $g$ commutes with all elements of ${\rm pFix}_G(P'\{0,1,\#\}^*)$ then
$g \in {\rm pFix}_G(P\{0,1,\#\}^*)$.
\end{pro}

\begin{pro} \label{G_31_endm_sep} \
Let $P$, and $P'$ be as in Definition \ref{def_max_sep_endm}, and let 
\ $G = G_{3,1}^{\rm mod \, 3}(0,1;\#)$.
Then the  fixator \ ${\rm pFix}_G(P'\{0,1,\#\}^*)$ \ is separating 
on $P \{0,1,\#\}^*$.
\end{pro}
Before proving Propositions \ref{max_sep_end_comm} and \ref{G_31_endm_sep}
we need some lemmas.

\begin{lem} \label{non_id} \ 
Let $P \subset A^*$ be any prefix code, where $|A| = n \geq 2$.  Assume 
$\varphi \in {\rm pStab}_{{\mathcal G}_{n,1}}(P A^*)$, but  
$\varphi \not\in {\rm pFix}_{{\mathcal G}_{n,1}}(P A^*)$.
Then there exists $x \in P A^*$ such that $x$ and 
$\varphi(x)$ are not prefix-comparable. 

In particular, if $\varphi \in {\mathcal G}_{n,1}$ is not the identity element
then there exists $x \in {\rm domC}(\varphi)$ such that $x$ and 
$\varphi(x)$ are not prefix-comparable.
\end{lem}
{\bf Proof.} \ The proof is in the Appendix dedicated to properties of
prefix codes.
 \ \ \ $\Box$  

\begin{lem} \label{extending_to_max_pref_code} \
Suppose $P, P' \subset A^*$ are complementary finite prefix codes.
Let $x_1, \ldots, x_k \in P A^*$ (for any positive integer $k$), and assume 
$x_1$, $\ldots$, $x_k$ are two-by-two prefix-incomparable. 
Then for all $n$ of the form \ $n = 1 + i \, (|A|-1)$, with \   
$n \geq \ |P| - k + (|A|-1) \, (|x_1| + \ldots + |x_k|)$, 
there exists a prefix code $Q$  such that

\smallskip

$\bullet$ \ \ $Q \cup \{x_1, \ldots, x_k\}$ and $P'$ are complementary 
prefix codes, with \ $Q \cup \{x_1, \ldots, x_k \} \subset P A^*$;

\smallskip

$\bullet$ \ \ $|Q| = n$.

\smallskip

$\bullet$ \ \ The set of prefixes of $P$ is a subset of the set of prefixes of
$Q \cup \{x_1, \ldots, x_k\}$.  
\end{lem}
{\bf Proof.} \ The proof is in the Appendix dedicated to properties of
prefix codes. 
 \ \ \ $\Box$
  
\begin{lem} \label{sep_implies_max_endm} \
Let $G$, $P$, and $P'$ be as in Definition \ref{def_max_sep_endm}.
If ${\rm pFix}_G(P'\{0,1,\#\}^*)$ is separating on $P \{0,1,\#\}^*$ then it 
is a maximal fixator.
\end{lem}
{\bf Proof.} \ Suppose by contradiction that there exists 
$x_0 \in P\{0,1,\#\}^*$ such that $h(x_0) = x_0$ for all 
$h \in {\rm pFix}_G(P'\{0,1,\#\}^*)$. The prefix code $P$ is of the form  
$P_1 \cup P_2\#$, with $P_1, P_2 \subset \{0,1\}^*$, by Lemma 
\ref{endmarker_code}. 

\smallskip

\noindent Case 1: \ $x_0 \in P_1 \{0,1\}^*$. 

Choose $x = x_00$ and $y = x_01$. Then $x$ and $y$ are prefix incomparable,
hence by the separation property of the fixator, there exists 
$h_0 \in {\rm pFix}_G(P'\{0,1,\#\}^*)$ and $u_0 \in \{0,1\}^*$ with

\smallskip
 
 \ \ \ \ \   $h_0(xu_0) = xu_0$, \ $h_0(yu_0) \neq yu_0$.

\smallskip

\noindent However, $h_0(yu_0) \neq yu_0$ contradicts the fact that  
$h_0(x_0) = x_0$.

\smallskip

\noindent Case 2: \ $x_0 \in P_1 \{0,1\}^* \# $, \ or \  
$x_0 \in P_2 \# $ \ with $|P_2| \geq 2$. 

Let $x_0 = v_0 \#$.
Let $w_0 \in P_2$ with $w_0 \neq v_0$, and choose $x = w_0\# $ 
and $y = v_0\# $. 
Then $x$ and $y$ are prefix incomparable, and both are in
$\{0,1\}^* \# $; so there exists 
$h_0 \in {\rm pFix}_G(P'\{0,1,\#\}^*)$ with

\smallskip
 
 \ \ \ \ \   $h_0(x) = x$, \ $h_0(y) \neq y$.

\smallskip

\noindent However, $h_0(y) \neq y$ contradicts the fact that  
$h_0(x_0) = x_0$.

\smallskip

\noindent Case 3: \ $x_0 \in P_2 \#$ \ and $|P_2| = 1$. (Obviously the case 
$|P_2| = 0$ cannot occur when $x_0 \in P_2 \#$.)

Then $P_2 = \{v_0 \}$, so we have \ $x_0 = v_0\#$. Let $z_0 \in P_1$ 
(recall that in the Definition \ref{def_max_sep_endm} we assume that 
$P_1 \neq \emptyset$). 
Let $x = z_0\#$ and $y = x_0 = v_0\#$. Since $z_0 \neq v_0$, \ 
$x$ and $y$ are prefix incomparable, and both are in $\{0,1\}^* \#$; 
so there exists $h_0 \in {\rm pFix}_G(P'\{0,1,\#\}^*)$ with

\smallskip

 \ \ \ \ \   $h_0(x) = x$, \ $h_0(y) \neq y$.

\smallskip

\noindent Again, $h_0(y) \neq y$ contradicts the fact that
$h_0(x_0) = x_0$.
  \ \ \ $\Box$

\bigskip

\noindent
{\bf Proof of Proposition \ref{max_sep_end_comm}.} \ Let $g \in G$ and 
assume $g$ commutes with all elements of ${\rm pFix}_G(P'\{0,1,\#\}^*)$. 
We want to show that $g \in {\rm pFix}_G(P\{0,1,\#\}^*)$. We first prove:

\smallskip

\noindent {\sc Claim:} \ $g$ stabilizes $P' \{0,1,\#\}^*$ and 
$P \{0,1,\#\}^*$.

\smallskip

\noindent Proof of the Claim: \ Assume by contradiction that $g(x') = y$ for 
some $x' \in P' \{0,1,\#\}^*$ and $y \in P \{0,1,\#\}^*$. Since $g$ commutes 
with all elements of the fixator we have for all 
$h \in {\rm pFix}_G(P'\{0,1,\#\}^*)$: 
 \ $gh(x') = hg(x') = g(x') = y$, i.e., $h(y) = y$. This contradicts the 
maximality of the fixator ${\rm pFix}_G(P'\{0,1,\#\}^*)$, proved in Lemma
\ref{sep_implies_max_endm}. So $g$ maps $P'\{0,1,\#\}^*$ into itself. 

In a similar way one proves that $g^{-1}$ maps $P'\{0,1,\#\}^*$ into itself. 
It follows from this that $g$ also maps $P \{0,1,\#\}^*$ into itself.
Indeed, if we had $g(x) = y'$ for some $x \in P \{0,1,\#\}^*$ and
$y' \in P' \{0,1,\#\}^*$ then $g^{-1}(y') = x$, contradicting the fact that
$g^{-1}$ maps $P'\{0,1,\#\}^*$ into itself.

Similarly, $g^{-1}$ maps $P \{0,1,\#\}^*$ into itself. This proves the Claim.

\medskip

\noindent Assume now by contradiction that $g$ does not fix some element 
$x_1 \in P \{0,1,\#\}^*$: \ $g(x_1) = y_1 \neq x_1$.
By the Claim, $y_1 \in P \{0,1,\#\}^*$. 

By Lemma \ref{non_id} there exist $x, y \in P \{0,1,\#\}^*$ such that $x$ 
and $y$ are prefix incomparable and $g(x) = y$.  And since $g$ commutes with 
the fixator, we have for all $h \in {\rm pFix}_G(P'\{0,1,\#\}^*)$: \   
$gh(x) =$ $hg(x) = h(y)$.

On the other hand, the separation property of the fixator implies that there 
exists $h_0 \in {\rm pFix}_G(P'\{0,1,\#\}^*)$ and $u_0 \in \{0,1,\#\}^*$
(with $u_0$ empty if $x, y \in \{0,1\}^*\#$), such that \    
 $h_0(yu_0) \neq yu_0$ and $h_0(xu_0) = xu_0$.

The equality $gh(x) = h(y)$ implies $gh_0(xu_0) = h(yu_0)$; this, together 
with $h_0(xu_0) = xu_0$, implies $yu_0 = gh_0(xu_0) = h(yu_0)$. But this
contradicts $h_0(yu_0) \neq yu_0$. 
       \ \ \ $\Box$

\begin{lem} \label{making_pref_codes} \  {\bf (1)} \  
For all $x, y \in \{0,1\}^*$ there exist letters $\ell_1, \ell_2 \in \{0,1\}$
such that $x \ell_1$, and $y \ell_2$ are prefix incomparable.

\noindent {\bf (2)} \  
For all $x, y, z \in \{0,1\}^*$ there exist letters $\ell_1, \ldots, \ell_6$ 
$\in \{0,1\}$ such that $x \ell_1 \ell_3$, $y \ell_2 \ell_4$, and 
$z \ell_5 \ell_6$, are prefix incomparable. 
\end{lem}
{\bf Proof.} \ The proof is in the Appendix dedicated to properties of 
prefix codes. 
       \ \ \ $\Box$

\bigskip

\noindent 
{\bf Notation:} \ When $S \subseteq A^*$, \   

\smallskip

$\geq_{{\rm pref}} \!\! (S) = $
$\{p \in A^* : p \geq_{\rm pref} s$, for some $s \in S \}$, 

\smallskip

\noindent i.e., $\geq_{{\rm pref}} \!\! (S)$ is the set of all prefixes of 
words of $S$.

\smallskip

$>_{{\rm pref}} \!\! (S) = $
$ \{p \in A^* : p >_{\rm pref} s$, for some $s \in S \}$, 

\smallskip

\noindent i.e., $>_{{\rm pref}} \!\! (S)$ is the set of all strict prefixes 
of words of $S$.

\bigskip

\noindent{\bf Proof of Proposition \ref{G_31_endm_sep}.} \ 
Let $x, y \in P_1 \{0,1\}^*$ and assume $x$ and $y$ are prefix incomparable.
We want to find $h_0 \in {\rm pFix}_G(P'\{0,1,\#\}^*)$ and 
$u_0 \in \{0,1\}^*$ such that $h_0(xu_0) = xu_0$ and 
$h_0(yu_0) \neq yu_0$. If $x, y \in \{0,1\}^*\#$ then $u_0$ is empty.

\medskip

\noindent {\sc Case 1:} \ $x, y \in P_1 \{0,1\}^*$.

\smallskip

The words $x, y0, y1$ are prefix-incomparable two-by-two (for $x$ and $y0$,
use Lemma \ref{pref_comp}, and similarly for $x$ and $y1$). Now use Lemma
\ref{extending_to_max_pref_code} to construct a maximal prefix code
$Q \cup \{x, y0, y1\} \cup P'$, with $Q \subset P \{0,1,\#\}^*$.

Define $h_0 \in G = G_{3,1}^{\rm mod \, 3}(0,1;\#)$ by 

\smallskip

$h_0(y0) = y1$, \  $h_0(y1) = y0$, \ $h_0(x) = x$, and $h$ is the identity 
on $Q \cup P'$.

\smallskip

\noindent So, $Q \cup \{x, y0, y1\} \cup P'$ is the domain code and image 
code of $h_0$. Note that $h_0$ preserves lengths.
Then $h_0 \in {\rm pFix}_G (P' \{0,1,\#\}^*)$, $h_0(y0) \neq y0$, and
$h_0(x 0) = x0$ (since $h_0(x) = x$). So here, $0$ plays the role of $u_0$
in the separation property.

\medskip

\noindent {\sc Case 2:} \   $x, y \in P_1 \{0,1\}^*\# \ \cup \ P_2 \#$.

\smallskip

Let $x = x_0\#$ and $y = y_0\#$

\smallskip

\noindent {\sc Case 2.1:} \ $y_0 \in P_1 \{0,1\}^*$.

\smallskip

Either $x_0$ is different from both $y_0$ and $y_0 0$, or $x_0$ is different
from both $y_0$ and $y_0 1$.
We only consider the case where $x_0$ is different from both $y_0$ and 
$y_0 0$; the other case is similar. 

\smallskip

\noindent $\bullet$ \ Assume $x_0 \in P_2$. 

By Lemma \ref{extending_to_max_pref_code} over the 
alphabet $A = \{0,1\}$, there is a finite prefix code 
$Q_1 \subset P_1 \{0,1\}^*$ such that
$Q_1 \cup \{y_0 00 \}$ and $P_1'$ and complementary prefix codes 
(over $\{0,1\}$). Therefore the following set \  
$C \subset \{0,1\}^* \cup \{0,1\}^*\#$ \ will be a finite 
maximal prefix code over $\{0,1,\#\}$: 

\smallskip

 \ \ \ \ \  $C \ = \ Q_1 \cup \{y_0 00\} \cup P'_1$
        $\cup \ >_{{\rm pref}}\!\!(Q_1 \cup \{y_0 00\} \cup P'_1) \ \#$,

\smallskip

\noindent Now we define $h_0$, with domain code and image code $C$, by

\smallskip

 \ \ \ \ \  $h_0(y_0\#) = y_0 0 \#$,  \ $h_0(y_0 0 \#) = y_0 \#$, \ 
            and $h_0$ is the identity everywhere else on $C$. 

\smallskip

\noindent Thus, $h_0(y) \neq y$. 
Moreover, $h_0 \in {\rm pFix}_G(P'\{0,1,\#\}^*)$  because $y_0, y_0 0$
$\notin P_2'$; indeed, $y_0, y_0 0$ $\in P_1 \{0,1\}^*$ 
$\subset P \{0,1,\# \}^*$. 

And $h_0$ preserves the length of strings in $\{0,1\}^*$ (since $h_0$ is 
the identity on $\{0,1\}^*$ wherever $h_0$ is defined).

We also claim that $h_0(x) = x$. 
Indeed, $x_0$ belongs to $P_2$, which is contained in 
$>_{{\rm pref}}\!\!(P_1 \cup P'_1)$; moreover, 
$>_{{\rm pref}}\!\!(P_1) \ \subset \ >_{{\rm pref}}\!\!(Q_1)$, by the 3rd 
point of Lemma \ref{extending_to_max_pref_code}. Therefore, $x_0 \#$ belongs 
to $C$. On the other hand, $x_0$ is different from $y_0$ and $y_0 0$.

\smallskip

\noindent $\bullet$ \ Assume $x_0 \in P_1 \{0,1\}^*$. 

Then, by Lemma \ref{making_pref_codes}, there are $\ell_1, \ell_2$ 
$\in \{0,1\}$ such that $x_0 \ell_1$ and $y_0 \ell_2$ are prefix
incomparable; also, $x_0 \ell_1$, $y_0 \ell_2$ $\in P_1 \{0,1\}^*$. 
By applying Lemma \ref{extending_to_max_pref_code} over 
the alphabet $A = \{0,1\}$ we obtain a finite prefix code
$Q_1 \subset P_1 \{0,1\}^*$ such that
$Q_1 \cup \{x_0 \ell_1, y_0 \ell_2 0\} $ and $P_1'$ and complementary 
prefix codes (over $\{0,1\}$). Therefore the following set \
$C \subset \{0,1\}^* \cup \{0,1\}^*\#$ \ will be a finite
maximal prefix code over $\{0,1,\#\}$:

\smallskip

 \ \ \ \ \  $C \ = \ Q_1 \cup \{x_0 \ell_1, y_0 \ell_2 0\} \cup P'_1$
  $\cup \ $
  $>_{{\rm pref}}\!\!(Q_1 \cup \{x_0 \ell_1, y_0 \ell_2 0\} \cup P'_1) \ \#$.

\smallskip

\noindent Now we define $h_0$, with domain code and image code $C$, by

\smallskip

 \ \ \ \ \ $h_0(y_0\#) = y_0 \ell_2 \#$, \ $h_0(y_0 \ell_2 \#) = y_0 \#$,
          \  and $h_0$ is the identity everywhere else on $C$.

\smallskip

\noindent Thus, $h_0(y) \neq y$.  Moreover, 
$h_0 \in {\rm pFix}_G(P'\{0,1,\#\}^*)$, because $y_0, y_0 \ell_2$
$\notin P_2'$; indeed, $y_0, y_0 \ell_2$ $\in P_1 \{0,1\}^*$
$\subset P \{0,1,\#\}^*$. 

And $h_0$ preserves the 
length of strings in $\{0,1\}^*$ (since $h_0$ is the identity on 
$\{0,1\}^*$ wherever it is defined).  Also, $h_0(x) = x$, since 
$x_0 \#$  belongs to $C$ (since $x_0$ is a strict prefix of 
$x_0 \ell_1$), and since $x_0$ is different from $y_0$ and 
$y_0\ell_2$. 

\medskip

\noindent {\sc Case 2.2:} \ $y_0 \in P_2$.

\smallskip

Since $P_1 \neq \emptyset$, there exists $w_0 \in P_1$; hence
$y_0$ is different from $w_0$, $w_0 0$, and $w_0 00$.    Also,
$x_0$ is different from $w_0 0$ or from $w_0 00$ (or from both).
Let $z_0 0$ be one of $w_0 0$ or $w_0 00$, so that $z_0 0\neq x_0$.
We still have $z_0 0 \neq y_0$ and $z_0 0 \in P_1 \{0,1\}^*$.

\smallskip

\noindent $\bullet$ \ Assume $x_0 \in P_2$.

By Lemma \ref{extending_to_max_pref_code} over the
alphabet $A = \{0,1\}$, there is a finite prefix code
$Q_1 \subset P_1 \{0,1\}^*$ such that
$Q_1 \cup \{z_00\}$ and $P_1'$ and complementary prefix codes
(over $\{0,1\}$). Therefore the following set \  
$C \subset \{0,1\}^* \cup \{0,1\}^*\#$  \ will be a finite maximal 
prefix code over $\{0,1,\#\}$:

\smallskip

 \ \ \ \ \  $C \ = \ Q_1 \cup \{z_0 0\} \cup P'_1 \ \cup \ $
 $>_{{\rm pref}}\!\!(Q_1 \cup \{z_0 0\} \cup P'_1) \ \#$.

\smallskip

\noindent Now we define $h_0$, with domain code and image code $C$, by

\smallskip

 \ \ \ \ \  $h_0(y_0\#) = z_0 \#$,  \ $h_0(z_0\#) = y_0\#$, \ 
            and $h_0$ is the identity everywhere else on $C$.

\smallskip

\noindent Thus, $h_0(y) \neq y$ and $h_0(x) = x$. Note that $h_0(x_0\#)$
and $h_0(y_0\#)$ are defined since $x_0, y_0 \in P_2 \ \subset \ $ 
$>_{{\rm pref}}\!\!(P_1 \cup P'_1)$; moreover, 
$>_{{\rm pref}}\!\!(P_1) \ \subset \  >_{{\rm pref}}\!\!(Q_1)$, by the 3rd 
point of Lemma \ref{extending_to_max_pref_code}. Therefore, $x_0\#$ and 
$y_0\#$ belong to $C$.

Also, $h_0 \in {\rm pFix}_G(P'\{0,1,\#\}^*)$, because $y_0, z_0$
$\notin P_2'$; indeed, $y_0, z_0$ $\in P_2 \cup P_1\{0,1\}^*$ 
$\subset P \{0,1,\#\}^*$.

Also, $h_0$ preserves the length of strings in $\{0,1\}^*$ since $h_0$ 
is the identity on $\{0,1\}^*$ wherever it is defined.

\smallskip

\noindent $\bullet$ \ Assume $x_0 \in P_1 \{0,1\}^*$.

Then, by Lemma \ref{making_pref_codes}, there are $\ell_1, \ell_2$
$\in \{0,1\}$ such that $x_0 \ell_1$ and $z_0 \ell_2$ are prefix
incomparable. By applying Lemma \ref{extending_to_max_pref_code} over
the alphabet $A = \{0,1\}$ we obtain a finite prefix code
$Q_1 \subset P_1 \{0,1\}^*$ such that
$Q_1 \cup \{x_0 \ell_1, z_0 \ell_2 \} $ and $P_1'$ and complementary
prefix codes (over $\{0,1\}$). Therefore the following set \
$C \subset \{0,1\}^* \cup \{0,1\}^*\#$ \ will be a finite
maximal prefix code over $\{0,1,\#\}$:

\smallskip

 \ \ \ \ \  $C \ = \ Q_1 \cup \{x_0 \ell_1, z_0 \ell_2 \} \cup P'_1$
  $\cup \ $
  $>_{{\rm pref}}\!\!(Q_1 \cup \{x_0 \ell_1, z_0 \ell_2 \} \cup P'_1) \ \#$.

\smallskip

\noindent Now we define $h_0$, with domain code and image code $C$, by

\smallskip

 \ \ \ \ \ $h_0(y_0\#) = z_0 \#$,  \ $h_0(z_0\#) = y_0\#$, \
            and $h_0$ is the identity everywhere else on $C$.

\smallskip

\noindent Thus, $h_0(y) \neq y$, and $y \in C$ (since 
$y_0 \in P_2 \ \subset \ $ $>_{{\rm pref}}\!\!(P_1 \cup P'_1)$ $\ \subset$ 
$ \ >_{{\rm pref}}\!\!(Q_1 \cup P'_1)$ ).

Moreover, $h_0 \in {\rm pFix}_G(P'\{0,1,\#\}^*)$, because $y_0, z_0$
$\notin P_2'$; indeed, $y_0, z_0$ $\in P_2 \cup P_1\{0,1\}^*$         
$\subset P \{0,1,\#\}^*$.

And $h_0$ preserves the length of strings in $\{0,1\}^*$ (since $h_0$ is 
the identity on $\{0,1\}^*$ wherever it is defined).  Also, $h_0(x) = x$, 
since $x_0 \#$  belongs to $C$ (since $x_0$ is a strict prefix of
$x_0 \ell_1$), and since $x_0$ is different from $y_0$ and $z_0$.
 \ \ \ $\Box$

\bigskip

As we observed near the beginning of this Section, the circuit equivalence
problem reduces to the generalized word problem of 
${\rm Fix}_G(0 \, \{0,1,\#\}^*)$, as subgroup of 
$G = G_{3,1}^{\rm mod \, 3}(0,1;\#)$. 
The generating set used for $G$ is \ 
$\Delta_{0,1;\#} \cup \{\tau_{i,i+1} : 0 \leq i \}$, where $\Delta_{0,1;\#}$ 
is a fixed finite generating set of $G_{3,1}^{\rm mod \, 3}(0,1;\#)$.
We will prove in the next Section, and independently of this Section, that 
$G_{3,1}^{\rm mod \, 3}(0,1;\#)$ is finitely presented. 

Theorem \ref{commutation_test} reduces the circuit equivalence problem to 
the word problem of $G_{3,1}^{\rm mod \, 3}(0,1;\#)$. 
The reduction is an unbounded 
conjunctive reduction, namely, the conjunction of all word problems 
``$gh=hg$'', as $h$ ranges over ${\rm Fix}_G( \{1,\# \}\{0,1,\#\}^*)$,
where $G = G_{3,1}^{\rm mod \, 3}(0,1;\#)$.

However, Proposition \ref{Fix_iso_V_endm} below implies that 
${\rm pFix}_G( \{1,\# \}\{0,1,\#\}^*)$ is isomorphic to 
$G = G_{3,1}^{\rm mod \, 3}(0,1;\#)$.  This and the fact that $G$ is finitely
generated (proved in Proposition \ref{V_endm_fin_gen}) implies that only 
the finitely many generators of ${\rm pFix}_G( \{1,\# \}\{0,1,\#\}^*)$ need 
to be used in the role of ``$h$'' in the Commutation Test. This then yields:

\begin{cor} \ 
 The circuit equivalence problem reduces to the word problem of 
$G_{3,1}^{\rm mod \, 3}(0,1;\#)$, and hence to the word problem of $G_{3,1}$ 
(over an infinite generating set), 
by a polynomial-time $k$-bounded conjunctive reduction. Here, $k$ is the 
minimum number of generators of $G_{3,1}^{\rm mod \, 3}(0,1;\#)$.  
\end{cor}
  
\begin{pro} \label{Fix_iso_V_endm}  \   
For $G = G_{3,1}^{\rm mod \, 3}(0,1;\#)$, the subgroup 
${\rm pFix}_G(\{1,\# \}\{0,1,\#\}^*)$ is isomorphic to $G$.    
\end{pro} 
{\bf Proof.} \ An element 
$\varphi \in  G = G_{3,1}^{\rm mod \, 3}(0,1;\#)$ belongs 
to  ${\rm Fix}_G(\{1,\# \}\{0,1,\#\}^*)$ iff $\varphi$ has 
a table of the form 

\[ \varphi \ = \   \left[ \begin{array}{cc ccc ccc}
1 & \# & 0x_1 & \ldots & 0x_n & 0x_1'\# & \ldots & 0x_m'\#  \\
1 & \# & 0y_1 & \ldots & 0y_n & 0y_1'\# & \ldots & 0y_m'\# 
\end{array}        \right]
\]
where $x_i,y_i, x'_j, y'_j$ range over $\{0,1\}^*$, and $|x_i| \equiv |y_i|$
mod 3 (for $i = 1, \ldots, n$).
The isomorphism to $G_{3,1}^{\rm mod \, 3}(0,1;\#)$ simply maps this table to 
\[ \psi \ = \  \left[ \begin{array}{ccc ccc}
x_1 & \ldots & x_n & x_1'\# & \ldots & x_m'\#  \\
y_1 & \ldots & y_n & y_1'\# & \ldots & y_m'\# 
\end{array}        \right]
\]
It is straightforward to see that $\psi$ preserves lengths mod 3 on 
$\{0,1\}^*$ if $\varphi$ does, and that $\varphi \mapsto \psi$ is an 
isomorphism.
 \ \ \ $\Box$

\bigskip

The commutation test not only works for certain fixators in 
$G_{3,1}^{\rm mod \, 3}(0,1;\#)$, but also for the analogous fixators in 
$G_{3,1}$, $G_{3,1}^{\rm mod \, 3}$, and $G_{3,1}^{\rm mod \, 3}(0,1)$. 
This is proved in the Appendices A2 and A3. 

\bigskip

Our next task will be to reduce this non-standard word problem of $G_{3,1}$
(over an infinite generating set) to the word problem of a finitely
generated group; we will actually obtain a finitely presented group.

%%%%%%%%%%%%%%%%%%%%%%%%%%%%%%%%%%%%%%%%%%%%%%%%%%%%%%%%%%%%%%%%%%%
% Section 6
%%%%%%%%%%%%%%%%%%%%%%%%%%%%%%%%%%%%%%%%%%%%%%%%%%%%%%%%%%%%%%%%%%%

\section{Finite presentations}

We will now prove that the groups $G_{3,1}(0,1;\#)$ and
$G_{3,1}^{\rm mod \, 3}(0,1;\#)$
are finitely generated, and in fact finitely presented. Higman's technique 
(see pp.~24-33 of \cite{Hig74}) can be applied rather directly to these 
groups, once we have proved certain properties of prefix codes. 
We will use Higman's notation 
\[  \left[ \begin{array}{ccc}
x_1 & \ldots & x_n \\
y_1 & \ldots & y_n \\
z_1 & \ldots & z_n   \end{array}        \right]
\]
for a composite of the form  
\[  
\left[ \begin{array}{ccc}
y_1 & \ldots & y_n \\
z_1 & \ldots & z_n \end{array} \right]  \cdot 
\left[ \begin{array}{ccc}
x_1 & \ldots & x_n \\
y_1 & \ldots & y_n \end{array} \right] (\cdot) 
\]
where $\{x_1, \ldots, x_n\}$, \ $\{y_1, \ldots, y_n\}$, \  
 $\{z_1, \ldots, z_n\}$ are three maximal prefix codes of cardinality $n$. 

\smallskip

\noindent A remark on terminology: \ 
Higman uses the word ``depth'' of a prefix code to refer to the number of 
vertices of the inner tree (he has a different point of view, and does not 
talk about prefix codes or trees explicitly). We will not follow Higman's 
terminology and use the word {\em depth} for the actual depth of a tree, i.e., 
the number of edges in a longest path from the root to a leaf.  

\smallskip

We first give a lemma concerning the particular maximal prefix codes used 
in $G_{3,1}(0,1;\#)$. Recall that, for an alphabet $A$, the {\it tree of 
the free monoid  $A^*$} consists of the vertex set $A^*$ and the edge set 
 \  $\{ (w, wa) : w \in A^*, \ a \in A \}$; the tree is rooted, with the
empty word  $\varepsilon$ as the root. For a prefix code $P \subset A^*$, 
the {\it prefix tree} of $P$ consists of the vertex set \ 

\smallskip
 
 \ \ \ \ \ $\{ w \in A^* : w$ is a prefix of some element of $P\}$, 

\smallskip

\noindent with root $\varepsilon$. The edge set is  

\smallskip

 \ \ \ \ \ $\{ (w, wa) : a \in A, \ wa$ is a prefix of some element of $P\}$.  

\smallskip

\noindent So the elements of $P$ are the leaves of the prefix 
tree of $P$. The {\it inner (or internal) vertices} of a rooted tree are, 
by definition, the vertices that are not leaves (i.e., a vertex $v$ is 
internal iff there exists an edge $(v,w)$ in the tree, for some vertex $w$). 
The tree spanned by the inner vertices is called the {\em inner tree}. 
We will denote the inner tree of the prefix tree of a prefix code $P$
by $T_{\rm in}(P)$. 

\begin{lem} \label{V_endm_two_leaves} \ 
{\bf (0)} \ Every finite maximal prefix code $P$ over an alphabet $A$ 
(e.g., $A = \{0,1,\#\}$) has cardinality \ $|P| = 1 + (|A| -1) \, i_P$, 
where $i_P$ is the number of inner vertices of the prefix tree of $P$. 

If $|P| > 1$ then $P$ contains a subset 
of the form \ $uA$ (for some word $u \in A^*$). 

Also, for every integer $i \geq 0$, 
there exists a maximal prefix code $P$ over an alphabet $A$ of cardinality \  
$1 + (|A| -1) \, i$.

\medskip

\noindent {\bf (1)} \ If \ $P \subset \{0,1\}^* \, \{\varepsilon,\# \}$, 
and $|P| > 1$, then $P$ contains a subset of the form \ 

\smallskip

 \ \ \ \ \  $u \ \{0,1,\# \}$, \ \ \  for some $u \in \{0,1\}^*$

\medskip

\noindent {\bf (2)} \ For every integer $i \geq 3$ there is a maximal 
prefix code \ $P \subset \{0,1\}^* \, \{\varepsilon, \# \}$, with 
\ $|P| = 1 + 2i$, and with the following property: 

\smallskip

 $P$ contains a subset of the form \ $\{u, v\} \, \{0,1,\# \}$,
  \ \  for some $u, v \in \{0,1\}^*$, \ $u \neq v$.

\medskip

\noindent {\bf (3)} \ For every integer $i \geq 5$, there is a maximal
prefix code \ $P \subset \{0,1\}^* \, \{\varepsilon, \#\}$, with
 \ $|P| = 1 + 2i$, and with the following property:

\smallskip

 $P$ contains a subset of the form \ $\{u, v, w\} \, \{0,1,\# \}$,
  \ \  for some $u, v, w \in \{0,1\}^*$, 

  with $u, v, w$ distinct two-by-two.  
\end{lem}
{\bf Proof.} \ The proof is in the appendix dedicated to properties of
prefix codes. 
   \ \ \ $\Box$    

\medskip

For elements of $G_{3,1}$ that preserve length modulo $3$, the following
concept and lemma are important.

\begin{defn} \  
Let $P \subset A^*$ be finite set, where $A$ is any finite
alphabet. The {\em mod 3 cardinality} of $P$ is the triple 
$(n_0,n_1,n_2) \in {\mathbb N}^3$, such that (for $i = 0, 1, 2$): 

\smallskip

$n_i \ = \ |P \ \cap \ \{w \in A^* : |w| \equiv i \ {\rm mod} \ 3 \}|$. 
\end{defn}
Note that if $(n_0,n_1,n_2)$ is the mod 3 cardinality of $P$ then
$n_0 + n_1 + n_2 = |P|$.

\medskip

\noindent {\bf Observation:} \  
By Lemma \ref{V_endm_two_leaves} (1), if a prefix code $Q \subset A^*$ has 
cardinality 2 or more, its inner tree $T_{\rm in}(Q)$ has at least one leaf.  
Moreover, if $|Q|$ is large enough then either $T_{\rm in}(Q)$ has a second 
leaf, or it has two (or more) one-child vertices, both having equivalent
depths modulo 3; for this to hold, it suffices that $T_{\rm in}(Q)$ has 
depth $\geq 4$.
More generally, if $|Q|$ is large enough then $T_{\rm in}(Q)$ has one 
of the following:  \\    
(1) either $T_{\rm in}(Q)$ has three leaves (or more); \\   
(2) or it has two leaves and two (or more) one-child vertices, both having 
equivalent depths modulo 3;  \\    
(3) or it has one leaf, and two (or more) one-child vertices, both having 
equivalent depths modulo 3, and two additional one-child vertices (or more), 
both having equivalent depths modulo 3.  (For one of these three properties
to be true it suffices that $T_{\rm in}(Q)$ has depth $\geq  6$.)  

\begin{lem} \label{mod_3_code}.

\noindent $\bullet$ \ Suppose that there exists a maximal prefix
code $Q$ over the alphabet $\{0,1\}$, whose inner tree $T_{\rm in}(Q)$ has 
two one-child vertices at depths $\equiv i$ {\rm mod 3} (for some 
$i \in \{0,1,2\}$). 
Then there exists a maximal prefix code $P \subset \{0,1\}^*$ with the same
{\rm mod 3} cardinality as $Q$, and with the following property:

\smallskip

 there is a word $u \in \{0,1\}^*$ such that \
 $u \cdot \{0,1\} \subseteq P$ \ and \ $|u| \equiv i$ {\rm mod 3}.
 
\smallskip

\noindent 
Equivalently, the inner tree of the prefix code $P$ has a leaf at
depth $\equiv i$ {\rm mod 3}.

\medskip

\noindent  $\bullet$ \ More generally, let $k \geq 2$, let $i_1, \ldots, i_k$ 
$\in \{0,1,2\}$, and suppose that $T_{\rm in}(Q)$ has the following property:
For every $\lambda$ ($1 \leq \lambda \leq k$), \ $T_{\rm in}(Q)$ has a leaf of 
depth $\equiv i_{\lambda}$ or it has two one-child  
vertices at depths $\equiv i_{\lambda}$ {\rm mod 3}.

Then there exists a maximal prefix code $P \subset \{0,1\}^*$ with the same 
{\rm mod 3} cardinality as $Q$, and with the following property:

\smallskip

there are $k$ different words $u_1, \ldots, u_k$ $\in \{0,1\}^*$ such that \   
$\{u_1, \ldots, u_k \} \cdot \{0,1\} \subseteq P$ \ and \ 

\smallskip

$|u_1| \equiv i_1$, \ \ldots \ , \  $|u_k| \equiv i_k$,  \  {\rm mod 3}. 

\smallskip

\noindent Equivalently, the inner tree of the prefix code $P$ has at 
least $k$ leaves that have depths respectively  $\equiv i_1$, $\ldots$, 
$\equiv i_k$ {\rm mod 3}.
\end{lem}
{\bf Proof.} \ The proof is in the appendix dedicated to properties of
prefix codes.  
 \ \ \ $\Box$

\begin{lem} \label{V_endm_fin_gen} \   {\bf (1)} \    
The group $G_{3,1}(0,1;\#)$ is generated by elements of table-size 
$\leq 7$.  \\ 
{\bf (2)} \ The group $G_{3,1}^{\rm mod \, 3}(0,1;\#)$ is generated by 
elements of table-size $\leq 61$. 

Hence these groups are finitely generated
\end{lem}
{\bf Proof.} \ {\bf (1)} Higman's proof  that $G_{N,r}$ is finitely generated
can be applied directly (see \cite{Hig74}, Lemma 4.2, pp.~26-27).
By Lemma \ref{V_endm_two_leaves} (1), every element 
$\varphi \in G_{3,1}(0,1;\#)$ of table-size \ $\|\varphi\| = n > 1$ \ (in 
particular when $n > 7$) has a table of the form
\[  \left[ \begin{array}{ccc ccc}
x 0 & x 1 & x \# & x_4 & \ldots & x_n \\
y_1 & y_2 & y_3  & y_4 & \ldots & y_n \end{array} \right].
\]
where $x \in \{0,1\}^*$; $x$ is a leaf of the inner tree of the domain code
domC$(\varphi)$.
The image code $\{y_1, \ldots, y_n\}$ also contains 3 words of the form  \ 
$y_{i_1} = y 0$, \ $y_{i_2} = y 1$, \ $y_{i_3} = y \#$, where 
$y \in \{0,1\}^*$.
The three indices $i_1, i_2, i_3$ are in $\{1, \ldots, n\}$, but any order 
relation between $i_1, i_2, i_3$ is possible. For the relation between 
$\{1, 2, 3\}$ and $\{i_1, i_2, i_3\}$ we have two cases, just as in 
\cite{Hig74}.

\medskip

\noindent {\sc Case 1} --- The column index sets $\{1, 2, 3\}$ and 
$\{i_1, i_2, i_3\}$ are disjoint:

\smallskip

\noindent
By permuting columns (if necessary) we can make $(i_1,i_2,i_3) = (4,5,6)$;
then the table of $\varphi$ has the form
\[  \left[ \begin{array}{ccc ccc ccc}
x0 & x1 & x\# & x_4 & x_5 & x_6 & x_7 & \ldots & x_n \\
y_1 & y_2 & y_3 & y0 & y1 & y\# & y_7 & \ldots & y_n 
 \end{array}        \right]
\]
If $n \geq 7$, we can apply Lemma \ref{V_endm_two_leaves}(2) to obtain a 
maximal prefix code $P_1$ over $\{0,1\}$ with two leaves and with the same 
cardinality as domC$(\varphi) \cap \{0,1\}^*$. Then (by Lemma 
\ref{endmarker_code}), $P_1$ determines a maximal prefix code 
$P \subset \{0,1\}^* \cup \{0,1\}^*\#$ with the same cardinality as 
domC$(\varphi)$. If we appropriately insert the code $P$ as a row we get 
\[  \left[ \begin{array}{ccc ccc ccc}
x0 & x1 & x\# & x_4 & x_5 & x_6 & x_7 & \ldots & x_n \\  
u 0 & u1 & u\# & v 0 & v1 & v\# & z_7 & \ldots & z_n \\  
y_1 & y_2 & y_3 & y0 & y1 & y\# & y_7 & \ldots & y_n
 \end{array}        \right]
\]
Thus, we can write the original element $\varphi \in G_{3,1}(0,1;\#)$ as a 
composite of two elements of $G_{3,1}(0,1;\#)$. 
Each of these two factors of $\varphi$ contains 3 columns in ``reducible''
form: Each of these two factors can be extended to a table of size $\leq n-2$,
obtained by replacing the three columns \
\(\left[ \begin{array}{ccc}
x0 & x1 & x\# \\
u 0 & u1 & u\#  \end{array}        \right] \)
 \ by the column \
\(\left[ \begin{array}{c}
x \\
u  \end{array}        \right] \), and similarly for \
\(\left[ \begin{array}{c}
v \\
y  \end{array}        \right] \).

Let us check that these two factors of $\varphi$ belong to $G_{3,1}(0,1;\#)$ 
(and not just to $G_{3,1}$). First, the inserted row corresponds to a maximal 
prefix code in \  $\{0,1\}^* \, \{\varepsilon, \#\}$. Since 
$\varphi \in G_{3,1}(0,1;\#)$, the table of $\varphi$ has the following
property: Words in $\{0,1\}^*$ line up (column-wise) with words in 
$\{0,1\}^*$, and words in $\{0,1\}^*\#$ line up (column-wise) with words in 
in $\{0,1\}^*\#$. The inserted row has the same size as the table of 
$\varphi$, and for maximal prefix codes in \  
$\{0,1\}^* \, \{\varepsilon, \#\}$ \ we know that 
the cardinality of the code determines the number of words in $\{0,1\}^*$ 
(or in $\{0,1\}^*\#$); see Lemma \ref{endmarker_code}. 
Thus we can correctly line up the elements of the new row with the two rows 
of $\varphi$, in such a way that the two factors belong to $G_{3,1}(0,1;\#)$. 

\medskip

\noindent {\sc Case 2} --- The column index sets $\{1, 2, 3\}$ and
$\{i_1, i_2, i_3\}$ have a non-empty intersection: 

\smallskip

\noindent
Then, we use Lemma \ref{V_endm_two_leaves} (2) to create a code, and 
insert it into the table of $\varphi$ as two rows, exactly as on p.~27 of 
\cite{Hig74}:
\[  \left[ \begin{array}{ccc ccc cc cccc}
x0 & x1 & x\# & x_4    & \ldots  & \ldots & \ldots & \ldots & 
                               x_{n-3} & x_{n-2} & x_{n-1} & x_n \\
u0 & u1 & u\# & \ldots & \ldots  & \ldots & \ldots & \ldots &
 \ldots  & v0 & v1 & v\#  \\
\ldots & ua_1 & \ldots & ua_2 & \ldots & ua_3 & \ldots & \ldots & 
 \ldots & v0 & v1 & v\#  \\
\ldots & ya_1 & \ldots & ya_2 & \ldots & ya_3 & \ldots & \ldots &
 \ldots & y_{n-2} & y_{n-1} & y_n 
\end{array}        \right]
\]
where $(a_1, a_2, a_3)$ is a permutation of $(0,1,\#)$; $ua_1$ is in 
any column from 1 through 3 (not necessarily in column 2, as drawn on the 
picture); $ua_2$ is in any column from the one just right of the column of
$ua_1$ through $n-4$, and $ua_3$ is in any column to the right of column 
$ua_2$ through column $n-3$.
This case is possible whenever $n$ is large enough so that there are 3 copies 
of the triple $(0,1,\#)$ with one overlap: \ $n = 1+2i \geq 3 \cdot |A| -1$
$= 3 \cdot 3 -1 = 8$, i.e., $n \geq 9$ \ (where $A = \{0,1,\#\}$). 
This will lead to a factorization of $\varphi \in G_{3,1}(0,1;\#)$ as a 
composition of three elements, each of which can be extended to a table 
of size $\leq n-2$. 

As in case 1, the two new rows can be inserted so that columns are be lined 
up in such a way that the three factors belong to $G_{3,1}(0,1;\#)$ (and not 
just to $G_{3,1}$). 

\smallskip

The Lemma now follows by induction on the table-size. Elements of table-size 
$< 9$ are then used as generators. Since over an alphabet of size 3, 
maximal prefix codes have size $1+2i$, it follows that the generators of 
table-size $< 9$ actually have table-size $\leq 7$.

\medskip

\noindent {\bf (2)} \ The proof that $G_{3,1}^{\rm mod \, 3}(0,1;\#)$ is 
finitely generated follows the same outline as the proof for 
$G_{3,1}(0,1;\#)$. The 
only difference is that now we have to check that the factors are in 
$G_{3,1}^{\rm mod \, 3}(0,1;\#)$, not just in $G_{3,1}(0,1;\#)$.
Let $\varphi \in G_{3,1}^{\rm mod \, 3}(0,1;\#)$.

\smallskip

\noindent {\sc Case 1} --- The column index sets $\{1, 2, 3\}$ and
$\{i_1, i_2, i_3\}$ are disjoint:

\smallskip

\noindent Again, $\varphi$ has the form
\[  \left[ \begin{array}{ccc ccc ccc}
x0 & x1 & x\# & x_4 & x_5 & x_6 & x_7 & \ldots & x_n \\
y_1 & y_2 & y_3 & y0 & y1 & y\# & y_7 & \ldots & y_n
 \end{array}        \right]
\]
where $x, y, x_4, x_5, y_1, y_2 \in \{0,1\}^*$, \ 
$|x| + 1 \equiv |y_1| \equiv |y_2| \equiv j$ mod 3, and 
$|y|+1 \equiv |x_4| \equiv |x_5| \equiv i$ mod 3. 

Let $Q = {\rm domC}(\varphi)$ and let
$Q_1 = {\rm domC}(\varphi) \cap \{0,1\}^*$. Then $x$ labels a leaf of
$T_{\rm in}(Q_1)$. Moreover, $x_4, x_5$ are either the children of a leaf of 
$T_{\rm in}(Q_1)$ or they are the children of two one-child vertices, both 
having equivalent depths modulo 3. So we can apply Lemma \ref{mod_3_code} in 
order to obtain a maximal prefix code $P_1 \subset \{0,1\}^*$ with the same 
{\rm mod 3} cardinality as $Q_1$, such that $T_{\rm in}(P_1)$ has a leaf at 
depth $\equiv j-1$ and a leaf at depth $\equiv i-1$ mod 3. So, $P_1$ has the 
form $P_1 = \{u0, u1, v0, v1, \ldots \}$, with $|u| + 1 \equiv i$ and 
$|v| + 1 \equiv j$ mod 3. By Lemma \ref{endmarker_code}, this uniquely 
determines a maximal prefix code $P = P_1 \cup P_2 \#$ over $\{0,1,\#\}$ 
(where $P_2 \subset \{0,1\}^*$ consists of the strict prefixes of elements 
of $P_1$). 

Now, as in proof (1) for $G_{3,1}(0,1;\#)$, we insert the code $P$ as
a row into the table of $\varphi$. We line up  the colums as in case 1 of (1).  
The columns of the table can be lined up so that the factors of $\varphi$ are 
in $G_{3,1}^{\rm mod \, 3}(0,1;\#)$; indeed, \ $P_1$, ${\rm domC}(\varphi)$,
and ${\rm imC}(\varphi)$ have the same mod 3 cardinality, and \ 
$|u| + 1 \equiv |x| + 1 \equiv |y_1| \equiv |y_2|$, \ 
$|v| + 1 \equiv |y|+1 \equiv |x_4| \equiv |x_5|$ mod 3.
\medskip

\noindent {\sc Case 2} --- The column index sets $\{1, 2, 3\}$ and
$\{i_1, i_2, i_3\}$ have a non-empty intersection:

\smallskip

Again, let $Q = {\rm domC}(\varphi)$, 
$Q_1 = {\rm domC}(\varphi) \cap \{0,1\}^*$; also, $x$ labels a leaf of
$T_{\rm in}(Q_1)$. We assume that $n$ is large enough in order to make sure 
that $T_{\rm in}(Q_1)$ has either another leaf or two one-child vertices,
both having equivalent depths modulo 3. If the depth of $T_{\rm in}(Q_1)$ is
at least 4 then this will be the case, and we can apply Lemma \ref{mod_3_code}.
In order to make sure that $T_{\rm in}(Q_1)$ has depth $\geq 4$ we assume
that $|Q_1| \geq 2^5$, and this is equivalent to assuming \ 
$n = |{\rm domC}(\varphi)| = |Q_1 \cup Q_2 \#| = 2 \, |Q_1| - 1$
$\geq 2 \, 2^5 - 1$.\ (Recall the $Q_2 \subset \{0,1\}^*$ consists of all 
strict prefixes of elements of $Q_1$, hence $|Q_2| = |Q_1| - 1$.) 
Thus, we assume $n \geq 2^6 -1 = 63$.
 
Now we insert the new code $P$ twice into the table, in the same way as in 
case 2 of (1).  Elements of odd table size $< 63$ can thus be used as 
generators.
  \ \ \ $\Box$

\bigskip

Next, we want to prove that $G_{3,1}(0,1;\#)$ and 
$G_{3,1}^{\rm mod \, 3}(0,1;\#)$ are finitely {\it presented}. 
Following Higman \cite{Hig74} (p.~25), we associate a table with a relation 
in $G_{3,1}(0,1;\#)$ or $G_{3,1}^{\rm mod \, 3}(0,1;\#)$. 
Let us fix a finite generating set for $G_{3,1}(0,1;\#)$, and let \   
$\varphi_1, \ldots, \varphi_n$ \ be a sequence of generators.
By restriction of the generators, we can choose a table for 
each generator in such a way that the image code of $\varphi_i$ is equal to 
the domain code of $\varphi_{i+1}$ $(1 \leq i < n)$.   
Then all these domain and image codes have the same cardinality, say $m$.
Putting these $n$ tables together in an $(n+1) \times m$ table yields the 
{\it table of the sequence} $\varphi_1 \ldots \varphi_n$:
\[  \left[ \begin{array}{ccc}
x_{1,1} & \ldots & x_{1,m} \\  
\ldots &  \ldots & \ldots \\  
x_{n,1} & \ldots & x_{n,m} \\   
x_{n+1,1} & \ldots & x_{n+1,m}  \end{array} \right],
\]
where the following is a table for $\varphi_i$ ($1 \leq i \leq n$): 
\[  \left[ \begin{array}{ccc}
x_{i,1} & \ldots & x_{i,m} \\  
x_{i+1,1} & \ldots & x_{i+1,m} \end{array} \right].\]
Note that $\varphi_1 \ldots \varphi_n$ is a relator of $G_{3,1}(0,1;\#)$ 
 \ iff \ $x_{1,j} = x_{n+1,j}$ \ for all $j = 1, \ldots, m$ (i.e., the first 
and the last rows are equal).

The smallest $m$ for which a sequence (or, in particular, a relator) \ 
 $\varphi_1, \ldots, \varphi_n$ \ has a table, is called the {\it table-size}
of the sequence (or the relator).

The concepts of  ``table of a relator'', and ``table-size of a relator''
make sense for any sequence \ $\varphi_1, \ldots, \varphi_n$ \ of elements 
of $G_{3,1}$, or in particular of $G_{3,1}^{\rm mod \, 3}(0,1;\#)$.

Thanks to this concept we can formulate the previous Lemma 
\ref{V_endm_fin_gen}
in a slightly stronger way (similar to Higman's Lemma 4.3).

\begin{lem} \label{V_endm_fin_gen_rels} \ 
Every element $\varphi \in G_{3,1}(0,1;\#)$ of table-size $\|\varphi\| > 7$ 
can 
be represented by a word $w_{\varphi}$ over the set of elements of table-size 
$\leq 7$, such that the sequence $w_{\varphi}$ has table-size 
$\leq \|\varphi\|$. 

Similarly, every element $\varphi \in G_{3,1}^{\rm mod \, 3}(0,1;\#)$ of 
table-size $\|\varphi\| > 61$ can be represented by a word $w_{\varphi}$ over 
the set of elements of table-size $\leq 61$, and such that the sequence 
$w_{\varphi}$ has table-size $\leq \|\varphi\|$.
\end{lem}
{\bf Proof.} \  This follows from the proof of Lemma \ref{V_endm_fin_gen}.
In that proof, we started out with a table of $\varphi$ (of table-size 
$\|\varphi\|$), and repeatedly inserted rows. No {\em columns} are ever 
added, hence the table-size doesn't increase. 
See also the proof of Higman's Lemma 4.3 in \cite{Hig74}.
 \ \ \ $\Box$

\begin{pro} \label{V_endm_fin_pres} \
The group $G_{3,1}(0,1;\#)$ is presented by relators of table-size
$\leq 9$, in terms of generators of table-size $\leq 7$.

The group $G_{3,1}^{\rm mod \, 3}(0,1;\#)$ is presented by relators of 
table-size $\leq 125$ in terms of generators of table-size $\leq 61$. 
Hence $G_{3,1}(0,1;\#)$ and $G_{3,1}^{\rm mod \, 3}(0,1;\#)$ are finitely 
presented.
\end{pro}
{\bf Proof.} \  Higman's method for proving that $G_{N,r}$ is finitely 
presented can be applied directly (see \cite{Hig74}, pp.~29-33). 
Now we use part (3) of Lemma \ref{V_endm_two_leaves} for $G_{3,1}(0,1;\#)$,
and Lemma \ref{mod_3_code} for $G_{3,1}^{\rm mod \, 3}(0,1;\#)$.

For the same reason as in Lemma \ref{V_endm_fin_gen}, the new rows that are
inserted can be lined up (column-wise), in such a way that all pairs of 
adjacent rows represent elements of $G_{3,1}(0,1;\#)$ or 
$G_{3,1}^{\rm mod \, 3}(0,1;\#)$ (and not just of $G_{3,1}$). 

The number 9 for $G_{3,1}(0,1;\#)$ comes from the fact that, in order to do 
the row insertions the table-size $n = 1 + 2i$ has to be at least 
$4 \times |A| -2 =$ $4 \times 3 -2 = 10$  \ (where $A = \{0,1,\#\}$).
Hence $i \geq 5$, hence $n \geq 11$. So, for the generators we can pick 
table-size $< 11$  (which implies table-size $i \leq 9$, since over a
three-letter alphabet table-sizes are odd). Refer to p.~31 of \cite{Hig74} 
(the ``linkages between them'' occupy at most $4 \, |A|-2$ columns). 
 
For $G_{3,1}^{\rm mod \, 3}(0,1;\#)$, Higman's ``type III reductions'' 
require that
we insert a row corresponding to a prefix code with 3 leaves in the inner
tree. Since one of the pre-existing rows in the table already has two leaves,
we need the table-size to be large enough so that the Observation before Lemma
\ref{mod_3_code} applies. If $T_{\rm in}$ (over $\{0,1\}$) has depth at least 
5, and $T_{\rm in}$ has at least two leaves, then it either has 3 (or more 
leaves) or it has at least two one-child vertices such that the depths of 
these two vertices are equivalent mod 3. In the latter case we apply Lemma 
\ref{mod_3_code} to obtain the desired code. For $T_{\rm in}$ to have depth 5, 
it is sufficient for the code (over $\{0,1\}$) to have cardinality $2^6$. 
Hence (by Lemma \ref{mod_3_code}), the code over $\{0,1,\#\}$ has cardinality 
$2 \times 2^6 -1 = 127$.
So the presentation of $G_{3,1}^{\rm mod \, 3}(0,1;\#)$ uses tables of size 
$< 127$,
hence of size $\leq 125$ (since code sizes over a 3-letter alphabet are odd).
 \ \ \ $\Box$

\bigskip

\noindent The group $G_{3,1}(0,1;\#)$ maps onto $G_{3,1}$ by the homomorphism 

\medskip 

\hspace{1in}
$ \left[ \begin{array}{ccc ccc}
x_1 & \ldots & x_n & x_1'\# & \ldots & x_m'\#  \\
y_1 & \ldots & y_n & y_1'\# & \ldots & y_m'\# 
\end{array}        \right] $
 $ \ \ \ \ \longmapsto \ \ \ \ $
$ \left[ \begin{array}{ccc}
x_1 & \ldots & x_n \\
y_1 & \ldots & y_n \end{array}        \right] $

\medskip 

\noindent whose kernel is the normal subgroup ${\rm pFix}_G(\{0,1\}^*)$ 
of $G = G_{3,1}(0,1;\#)$ (by Lemma \ref{pStabpFixGroups} this partial fixator
is indeed a group). The group ${\rm pFix}_G(\{0,1\}^*)$ consists of the 
elements that have a table of the form

\bigskip 

\hspace{1in}
$ \left[ \begin{array}{ccc ccc}
x_1 & \ldots & x_n & x_1'\# & \ldots & x_m'\#  \\
x_1 & \ldots & x_n & y_1'\# & \ldots & y_m'\# 
\end{array}        \right]. $

\bigskip 

\noindent Hence, $G_{3,1}(0,1;\#)$ is not a simple group. 

In a similar way, $G_{3,1}^{\rm mod \, 3}(0,1;\#)$ maps onto 
$G = G_{3,1}^{\rm mod \, 3}$ with kernel ${\rm pFix}_G(\{0,1\}^*)$; hence
$G_{3,1}^{\rm mod \, 3}(0,1;\#)$ is not a simple group. 

\bigskip

\noindent In summary, we proved:

\begin{thm} \
The group $G_{3,1}^{\rm mod \, 3}(0,1;\#)$ is finitely presented, and not 
simple.

The word problem of $G_{3,1}^{\rm mod \, 3}(0,1;\#)$, over the infinite 
generating set \ 
$\Delta \cup \{\tau_{i,i+1} : 0 \leq i\}$, is coNP-hard with respect to 
constant-arity conjunctive polynomial-time reduction. (Here $\Delta$ is a 
finite generating set of $G_{3,1}^{\rm mod \, 3}(0,1;\#)$.)
\end{thm}
  
\bigskip

%%%%%%%%%%%%%%%%%%%%%%%%%%%%%%%%%%%%%%%%%%%%%%%%%%%%%%%%%%%%%%%%%%%
% Section 7
%%%%%%%%%%%%%%%%%%%%%%%%%%%%%%%%%%%%%%%%%%%%%%%%%%%%%%%%

\section{Reduction to the word problem of a finitely presented group}

So far, the word problems that we have focused on were over infinite 
generating sets, although the groups used also admit finite generating sets.
This is a crucial point, because the groups $G_{3,1}$, etc., have
their word problem in {\bf P} over a finite generating set; but their 
word problem over certain infinite generating sets, as seen here, is 
coNP-hard.
 
In this section we obtain different Thompson groups with finite generating 
sets. These groups are obtained by expressing the transpositions
$\tau_{n,n+1}$ over a finite set of generators, according to tbe methods 
of Section 2; we saw that $\tau_{n,n+1}$ has polynomial word length (in $n$)
over those generators. Thus Section 2 now gives us a finitely generated 
Thompson group with coNP-hard word problem.   
We will work next at obtaining a finitely presented group.

\begin{pro} \label{V_endm_closed} \ 
Conjugation by $\kappa_i$ or $\kappa_i^{-1}$ ($i = 0, 1, 2, 3$) is an 
automorphism of $G_{3,1}^{\rm mod \, 3}(0,1;\#)$, and also an automorphism of 
$G_{3,1}^{\rm mod \, 3}(0,1)$.
\end{pro}
{\bf Proof.} \ It is enough to prove that $G_{3,1}^{\rm mod \, 3}(0,1;\#)$ and 
$G_{3,1}^{\rm mod \, 3}(0,1)$ are  closed under conjugation by $\kappa_i$ and 
by $\kappa_i^{-1}$.   

The definition of $\kappa_i$ directly shows that $\kappa_i$ and 
$\kappa_i^{-1}$ stabilize $\{0,1\}^*$ and $\{0,1\}^*\#$. Hence for every 
$\varphi \in G_{3,1}^{\rm mod \, 3}(0,1)$, 
 \ $\kappa_i \varphi \kappa_i^{-1}$ and $\kappa_i^{-1} \varphi \kappa_i$
stabilize $\{0,1\}^*$; and for every $\varphi \in $
$G_{3,1}^{\rm mod  \, 3}(0,1;\#)$,
 \ $\kappa_i \varphi \kappa_i^{-1}$ and $\kappa_i^{-1} \varphi \kappa_i$ 
stabilize $\{0,1\}^*$ and $\{0,1\}^*\#$. 
 
The definition of $\kappa_i$ also directly shows that $\kappa_i$ is 
length-preserving. Hence, or every $\varphi \in G_{3,1}^{\rm mod \, 3}$, \ 
$\kappa_i \varphi \kappa_i^{-1}$ and $\kappa_i^{-1} \varphi \kappa_i$ 
preserve length of strings in $\{0,1\}^*$ modulo 3.

Thus, all we still need to show is that if 
$\varphi \in G_{3,1}^{\rm mod \, 3}(0,1)$, then
$\kappa_i \varphi \kappa_i^{-1}$ and $\kappa_i^{-1} \varphi \kappa_i$
belong to $G_{3,1}$. We will do this by showing that they have 
``finite depth''.  An element $\psi \in {\mathcal G}_{3,1}$ is said to have
depth $\leq d$ iff for all $w \in \{0,1,\# \}^*$ with $|w| > d$, there is a
prefix $v$ of $w = vs$ (for some $s \in \{0,1,\#\}^*$, with $|v| \leq d$ and 
$\psi(w) = \psi(v) \ s$.  Obviously, $\psi$ belongs to $G_{3,1}$ iff $\psi$ 
has finite depth.

Let $\varphi \in G_{3,1}^{\rm mod \, 3}(0,1)$.
Since $\kappa_i \varphi \kappa_i^{-1}$ stabilizes $\{0,1\}^*$, it has
domain and image codes of the form  \ 
$P_1 \cup \bigcup_{v \in P_2} v \, \# \, P(v)$, \ where 
$P_1$ ($\subset \{0,1\}^*$) is a maximal prefix code over $\{0,1\}$,
$P_2$ ($\subset \{0,1\}^*$) consists of all strict prefixes of elements of
$P_1$, and each $P(v)$ is a maximal prefix code over $\{0,1,\#\}$. 
In order to show that $\kappa_i \varphi \kappa_i^{-1}$ and 
$\kappa_i^{-1} \varphi \kappa_i$ belong to $G_{3,1}$, we have to show
that $P_1$ is finite (hence $P_2$ is finite), and that each $P(v)$ is finite,
as $v$ ranges over $P_2$).  Let 

\smallskip

$m \ = \ {\rm max}\{ |w| \ : \ w \in $
  ${\rm domC}(\varphi) \ \cup \ {\rm imC}(\varphi) \}$,

\smallskip

$d \ = \ i \ + \ 3 \cdot \lceil m/3 \rceil$ \ \ 
  (i.e., $3 \cdot \lceil m/3 \rceil$ \ is ``$m$ rounded up to the next 
      multiple of 3''). 

\smallskip

\noindent We claim:

\smallskip

$\kappa_i \varphi \kappa_i^{-1}$ \ and \ $\kappa_i^{-1} \varphi \kappa_i$
 \ \ have depth $\leq d$.

\smallskip

\noindent For $w \in \{0,1\}^*$, if $|w| > d$ we can write 
$w = xs \in P_1 \{0,1\}^* \subset \{0,1\}^*$, with $|x| = d$. Then 

\smallskip

$\kappa_i(xs\#) = \kappa_i(x) \, \kappa_0(s) \ \#$,

\smallskip

\noindent by the choice of $d$, and since $|x| = d$. Next,
applying $\varphi$ yields \

\smallskip

$\varphi(\kappa_i(x)) \ \kappa_0(s) \ \#$,

\smallskip

\noindent since \ $|\kappa_i(x)| = |x| \geq m$.
Now, applying $\kappa_i^{-1}$ yields \

\smallskip

$\kappa_i^{-1}(\varphi(\kappa_i(x))) \ \kappa_0^{-1}(\kappa_0(s)) \ \#$,

\smallskip

\noindent since
$|\varphi(\kappa_i(x))| \equiv |\kappa_i(x)| \ ({\rm mod} \ 3) $
$\equiv d \ ({\rm mod} \ 3)$.
 \ Important remark: Here we used the fact that $\varphi$ 
preserves length modulo 3 (on $\{0,1\}^*$). 

\smallskip

\noindent Thus we have:

\smallskip

$\kappa_i^{-1} \varphi \kappa_i(xs\#) \ = \ $
$\kappa_i^{-1}(\varphi(\kappa_i(x))) \ s \ \#$,

\smallskip

\noindent for all $xs \in \{0,1\}^*$ with $|xs| \geq d$, $|x| = d$.

\medskip

\noindent For $v\#s \, t$ with $v \in P_2$, \ $st \in P(v) \{0,1\}^*$,  
and $|v\#s| \leq d$ we have 

\smallskip

$\varphi(\kappa_i(v\#s \, t)) = \varphi(\kappa_i(v) \ \#s \, t)$
$= \varphi(\kappa_i(v) \ \#s) \ t$;

\smallskip

\noindent the last equality holds because $|v\#s| \leq d$.
Note that $\varphi(\kappa_i(v) \ \#s)$ contains at least one copy of 
$\#$, since $\varphi \in G_{3,1}^{\rm mod \, 3}(0,1)$, 
i.e., $\varphi(\kappa_i(v) \# s) = y \# z$ for some $y \in \{0,1\}^*$,
 $z \in \{0,1,\# \}^*$.

\smallskip

\noindent Therefore, when we apply $\kappa_i^{-1}$ we obtain 

\smallskip

$\kappa_i^{-1}(\varphi(\kappa_i(v) \ \#s)) \ t$ 
 $ \ = \ \kappa_i^{-1}(y) \# zt$.

\smallskip

This shows that $\kappa_i^{-1} \varphi \kappa_i$ has depth $\leq d$.
Hence, $P_1$, $P_2$, and all $P(v)$ are finite.
For $\kappa_i \varphi \kappa_i^{-1}$ the proof is the same.
  \ \ \  $\Box$

\bigskip

\noindent 
As a consequence of Proposition \ref{V_endm_closed} we can consider the 
following HNN-extension:

\medskip

$H(0,1;\#) \ = \ $
$\langle G_{3,1}^{\rm mod \, 3}(0,1;\#) \cup \{t\} \ : \ $
  $\{t \, g \, t^{-1} = g^{\kappa_{321}} : g \in $
  $G_{3,1}^{\rm mod \, 3}(0,1;\#) \} \rangle$.

\medskip

\noindent Since $G_{3,1}^{\rm mod \, 3}(0,1;\#)$ is finitely generated, the 
HNN-relations form a finite set; moreover, since 
$G_{3,1}^{\rm mod \, 3}(0,1;\#)$ 
is finitely presented, the whole HNN-extension is a finitely presented group. 

This HNN-extension is rather special, since the group being extended
is the same as the group being conjugated. Therefore, the normal form of 
elements of the HNN-extension $H(0,1;\#)$ is  

\smallskip 

\hspace{1.5in}  $g t^n$, \ \ where $n \in {\mathbb Z}$ and 
                        $g \in G_{3,1}^{\rm mod \, 3}(0,1;\#)$.

\smallskip 

\noindent 
It follows that this HNN-extension is a {\it semidirect product}:

\medskip

 \ \ \ \ \ 
$H(0,1;\#) \ \cong \ G_{3,1}^{\rm mod \, 3}(0,1;\#) \rtimes {\mathbb Z}$.

\smallskip 
 
\begin{lem} \label{H_t_kappa} \ 
The HNN-extension $H(0,1;\#)$ is isomorphic to the subgroup \
$\langle G_{3,1}^{\rm mod \, 3}(0,1;\#) \cup \{\kappa_{321}\} \rangle$ \ of 
the group ${\mathcal G}_{3,1}$.
\end{lem}
{\bf Proof:} \ 
By the normal form theorem for HNN extensions, the mapping defined by 
 \ $t \mapsto \kappa$ \ determines a surjective homomorphism from $H(0,1;\#)$ 
onto $\langle G_{3,1}^{\rm mod \, 3}(0,1;\#) \cup \{\kappa\} \rangle$.

Here we abbreviate $\kappa_{321}$ to $\kappa$. 

In order to show that the map 
 \ $g t^n \longmapsto g \kappa^n$ \ 
has trivial kernel, suppose by contradiction that for some $n \neq 0$, an 
element \ $\varphi \ = \ g \kappa^n$ \ is the identity. 

Since $g \in G_{3,1}$, it has finite domain and finite image codes. Let 
$\ell$ be an upper bound on the longest length of any element in the 
domain code and the image code of $g$. Let $B$ be an integer such that \  
$B > 6 \, |n|$, and $B > 2 \, \ell$.   

Let $x \in \{0,1\}^*$ be of length $> 3B$, and let us apply $\varphi$ to
the argument $x\#$. The map $g \in G_{3,1}^{\rm mod \, 3}(0,1;\#)$ can change 
at most $\ell$ bits of the argument, 
or shorten or lengthen the argument by $< \ell$ bits. The map $\kappa^n$ 
moves bits over a distance $\leq 6 \, |n|$.
Therefore, the effect of $g$ on the argument $x\#$ is only felt on the 
leftmost \ $B = 2 \, \ell + 6 \, |n|$ \ bits of the argument. Further to the
right inside $x\#$, only $\kappa^n$ has an effect.
So we can write $x$ as $x = ps$ with $p, s \in \{0,1\}^*$, \ 
$|p| = B$; then $\varphi(ps\#)$ has the form 

\smallskip

 \ \ \ \ \   $\varphi(ps\#) = p's'\#$

\smallskip

\noindent for some $p', s' \in \{0,1\}^*$, with $|p'| \leq B$.
Most importantly, $s'$ is changed (according to $\kappa^n$) at every bit 
position, except perhaps in the rightmost $6 \,|n|$ ($< B$) bits. Since 
we chose $|x| > 3B$, we conclude: $\varphi$ changes $x\#$. So $\varphi$ is 
not the identity map.  The completes the proof by contradiction.
  \ \ \  $\Box$

\bigskip

\noindent  In summary, so far we have proved the following.

\begin{thm} \label{finprescoNPhard} \  
There exists a finitely presented Thompson group $G$ 
($\subset {\mathcal G}_{3,1}$), with the following properties:

\smallskip

\noindent $\bullet$ \ The word problem of $G$ (over a fixed finite 
generating set) is coNP-hard, with respect to polynomial-time 
constant-arity conjunctive reduction. 

\smallskip

\noindent $\bullet$ \ $G$ is an HNN extension (by one stable letter) of
some finitely presented subgroup {\rm Th} of $G_{3,1}$. In fact, $G$ is 
isomorphic to the semidirect product ${\rm Th} \rtimes {\mathbb Z}$.

\smallskip

An example of such a group $G$ is the subgroup  \ 
$\langle G_{3,1}^{\rm mod \, 3}(0,1;\#) \cup \{\kappa_{321}\} \rangle$ \ of \  
${\mathcal G}_{3,1}$, where {\rm Th} is $G_{3,1}^{\rm mod \, 3}(0,1;\#)$.
\end{thm}
{\bf Proof.} \ We use $\kappa_{321}$ to replace the transpositions 
$\tau_{n,n+1}$ by words over $\Delta \cup \{\kappa_{321}\}$ of linear
length (according to Lemma \ref{tau_generated}); here $\Delta$ is a finite 
generating set of $G_{3,1}^{\rm mod \, 3}(0,1;\#)$. 
Now the previously seen reductions reduce the circuit equivalence problem
to the word problem of \ 
$\langle G_{3,1}^{\rm mod \, 3}(0,1;\#) \cup \{\kappa_{321}\} \rangle$.
 \ \ \ $\Box$

\bigskip

In the next Section we will show that the word problem of \ 
$\langle G_{3,1}^{\rm mod \, 3}(0,1;\#) \cup \{\kappa_{321}\} \rangle$  \  
is in coNP, thus showing that this word problem is coNP-complete.

%%%%%%%%%%%%%%%%%%%%%%%%%%%%%%%%%%%%%%%%%%%%%%%%%%%%%%%%%%%%%%%%%%%
% Section 8 
%%%%%%%%%%%%%%%%%%%%%%%%%%%%%%%%%%%%%%%%%%%%%%%%%%%%%%%%

\section{Complexity of some word problems}

\noindent Consider the following subgroups of the Thompson-Higman group 
${\mathcal G}_{3,1}$: 

\medskip

$H(0,1) \ = \ $
$\langle G_{3,1}^{\rm mod \, 3}(0,1) \cup \{\kappa_{321}\} \rangle$,

\smallskip

$H(0,1;\#) \ = \ $
$\langle G_{3,1}^{\rm mod \, 3}(0,1;\#) \cup \{\kappa_{321}\} \rangle$, 

\smallskip

$\langle G_{3,1} \cup \{\kappa_{321}\} \rangle$, \ \ and

\smallskip

$\langle G_{3,1} \cup \{\kappa_0, \kappa_1, \kappa_2\} \rangle$ .

\medskip

\noindent We see from the definition of $\kappa_0$ and $\kappa_3$ that they
differ only by a finite permutation; hence \ 
$\langle G_{3,1} \cup \{\kappa_0, \kappa_1, \kappa_2\} \rangle \ = \ $
$\langle G_{3,1} \cup \{\kappa_0, \kappa_1, \kappa_2, \kappa_3 \} \rangle$.

Before we analyze the word problem of these groups we need a 
result about the permutation group  
$\langle \gamma_0, \gamma_1, \gamma_2 \rangle$ 
(of permutations of $\mathbb{N}$),  and about the subgroup 
$\langle \kappa_0, \kappa_1, \kappa_2 \rangle$ of
${\mathcal G}_{3,1}$.
For $\pi \in \langle \gamma_0, \gamma_1, \gamma_2 \rangle$
we denote the word-length of $\pi$ over 
$\{\gamma_0, \gamma_1, \gamma_2 \}^{\pm 1}$ by $|\pi|$; 
similarly, for $K \in \langle \kappa_0, \kappa_1, \kappa_2 \rangle$,
the word-length of $K$ over 
$\{\kappa_0, \kappa_1\, \kappa_2 \}^{\pm 1}$ is denoted by $|K|$.

\begin{lem} \label{subgroupKappa0123} \  
Let $\pi \in $ $\langle \gamma_0, \gamma_1, \gamma_2  \rangle$.
Then for all $n \in \mathbb{N}$ with $n \geq 2 \, |\pi| + 1:$ 
 \ \ $\pi(n + 3) = \pi(n) + 3$.
Hence the displacement function $n \mapsto \pi(n) - n$  is ultimately 
periodic, with period $3$, when $n \geq 2 \, |\pi| + 1$. 

As a consequence, $\pi \neq {\bf 1}$ \ iff \ $\pi(n) \neq n$ for some 
$n \leq  2 \, |\pi| + 3$.  Similarly, for $K \in $
$\langle \kappa_0, \kappa_1, \kappa_2 \rangle$ we have:
$K \neq {\bf 1}$ \ iff \ there exists $x \in \{0,1\}^*$ with 
$|x| \leq 6 \, |K| + 3$, such that \ $K(x\#) \ \neq \ x\#$.  

The word problems of the groups 
$\langle \gamma_0, \gamma_1, \gamma_2 \rangle$ and 
$\langle \kappa_0, \kappa_1, \kappa_2 \rangle$ 
can be decided deterministically in quadratic time. 
\end{lem}
{\bf Proof.} \ From the definition of $\gamma_0$, $\gamma_1$, and $\gamma_2$,
one sees immediately that $\gamma_i(n+3) = \gamma_i(n) + 3$, for all 
$n \geq 3$, $i = 0, 1, 2$.
For $\pi \in \langle \gamma_0, \gamma_1, \gamma_2 \rangle$, the 
relation $\pi(n + 3) = \pi(n) + 3$ (when $n \geq 2 \, |\pi| + 1$)    
follows by a straightforward induction on $|\pi|$. Indeed, 
$\gamma_i \pi(n+3) = \gamma_i(\pi(n)+3) = \gamma_i\pi(n)+3$, if $n \geq 3$ 
and $\pi(n) \geq 3$. Moreover, since each $\gamma_i$ can decrement its 
argument by at most 2, we have $\pi(n) \geq 3$ if $n \geq 2 \, |\pi| +3$
$ = 2 \, |\gamma_i \pi| + 1$.

Let $\pi \in \langle \gamma_0, \gamma_1, \gamma_2 \rangle$.
To check whether $\pi = {\bf 1}$, we compute the $2 \, |\pi| + 4$ numbers
$\pi(n)$ with $0 \leq n \leq 2 \, |\pi| + 3$, and check whether $\pi(n) = n$.
Let $\pi = \pi_k \ldots \pi_1$, with $\pi_k, \ldots, \pi_1 \in$
$\{\gamma_0, \gamma_1, \gamma_2 \}^{\pm 1}$. To compute $\pi(n)$
we successively compute $\pi_1(n)$, $\pi_2\pi_1(n)$, $\ldots \ldots$,
$\pi_j \ldots \pi_1(n)$, $\ldots \ldots$, $\pi_k \ldots \pi_j \ldots \pi_1(n)$.
For this, all we need is a deterministic push-down automaton, whose input tape 
contains the word $(\pi_k, \ldots, \pi_1)$; inputs are read from right to left.
After reading $(\pi_j, \ldots, \pi_1)$ with $k \geq j \geq 1$, the machine's
stack contains the number $\pi_j \ldots \pi_1(n)$ in unary, and the machine's 
internal state remembers $\pi_j \ldots \pi_1(n)$ mod 3. To apply $\pi_{j+1}$ 
to $\pi_j \ldots \pi_1(n)$, the machine only needs to know
$\pi_j \ldots \pi_1(n)$ mod 3, and it needs to know whether
$\pi_j \ldots \pi_1(n)$ is equal to 0, 1, 2, or $> 2$. Since a push-down 
automaton has linear running time, $\pi(n)$ can thus be computed in time 
$O(n)$ $(\leq O(|\pi|)$.  Since
$0 \leq n \leq 2 \, |\pi| + 3$, the total time to compute $\pi(0)$, $\pi(1)$,
$\ldots,$ $\pi(2 \, |\pi| +3)$ is $O(|\pi|^2)$.

For $K \in \langle \kappa_0, \kappa_1, \kappa_2 \rangle$ and 
$x\# \in \{0,1\}^*\#$, the action of $K$ on $x\#$ permutes the bits of the 
bitstring $x$.
Note that $\kappa_i$ permutes the bits of $x\#$ in the same way as 
$\gamma_i^{-1}$ permutes the bit positions, except near $\#$.
More generally, when $|x| \geq 2 \, |K|$, the action of $K$ on $x\#$ permutes
the bits of $x$ in the same way as 
$\pi_K \in \langle \gamma_0, \gamma_1, \gamma_2 \rangle$, except perhaps for
the right-most $2 \, |K|$ bits of $x$ (near $\#$); here $\pi_K$ is obtained 
from $K$ by replacing every $\kappa_i$ by $\gamma_i^{-1}$ $(i = 0,1,2)$. 
Indeed, every $\kappa_i$ in $K$ differs from the corresponding $\gamma_i^{-1}$
at most on the 2 bits near $\#$; this effect propagates $|K|$ times,
to a distance $\leq 2 \, |K|$ from $\#$. 

If $K \neq {\bf 1}$, then either $\pi_K \neq {\bf 1}$, or $K$ is a 
non-identity permutation on the right-most $2 \, |K|$ positions of some words 
$x\# \in \{0,1\}^*\#$. Note that $|\pi_K| \leq |K|$.
When $|x| \geq 6 \, |K| + 3$, the action of $K$ on $x\#$ consists of applying
$\pi_K$  on $x$, except for the right-most $2 \, |K|$ bits. 
Thus, if $\pi_K \neq {\bf 1}$, we can check this on the left-most 
$4 \, |\pi_K| +3 \ ( \, \leq 4 \, |K| +3)$ bits of $x\#$; if 
$\pi_K = {\bf 1}$, we can check that $K$ is a non-identity permutation 
on the right-most $2 \, |K|$ positions by inspecting these $2 \, |K|$ 
positions. Therefore, if $K \neq {\bf 1}$, there is a position 
$n \leq 6 \, |K| + 3$ which is permuted non-identically by $K$.
Therefore, to decide the word problem for 
$K \in \langle \kappa_0, \kappa_1, \kappa_2 \rangle$ we can proceed as for 
$\langle \gamma_0, \gamma_1, \gamma_2 \rangle$, above, but we check how 
$K$ permutes all $n$ with $n \leq 6 \, |K| + 3$ (instead of 
$\leq 4 \, |\pi| + 3$).
   \ \ \ $\Box$

\begin{thm} \label{coNP} \  
The word problem of 
$\langle G_{3,1} \cup \{\kappa_0,\kappa_1, \kappa_2 \} \rangle$, 
and hence of $H(0,1)$, $H(0,1;\#)$, and
$\langle G_{3,1} \cup \{\kappa_{321}\} \rangle$, are in {\rm coNP}.
\end{thm}
{\bf Proof.} \ Since $G_{3,1}^{\rm mod \, 3}(0,1)$ and 
$G_{3,1}^{\rm mod \, 3}(0,1;\#)$
are finitely generated subgroups of $G_{3,1}$, and 
$\langle G_{3,1} \cup \{\kappa_{321}\} \rangle$ is a finitely generated
subgroup of 
$\langle G_{3,1} \cup \{\kappa_0,\kappa_1, \kappa_2 \} \rangle$, 
it is sufficient to show that the word problem of 
$\langle G_{3,1} \cup \{\kappa_0,\kappa_1, \kappa_2 \} \rangle$
is in coNP. Indeed, it is a general fact that if a group's word problem has 
a complexity $\leq f(n)$ (regarding time of space, deterministic,
nondeterministic, or co-nondeterministic), then every finitely generated 
subgroup has a word problem of complexity $\leq f(cn)$, for some positive
constant $c$ (see \cite{MO}, and \cite{Bi}).

Let $\Delta_{3,1}$ be a finite generating set of $G_{3,1}$. We will prove 
(in the Claim below) that if a word $w$ over the generating set \ 
$\Delta_{3,1}^{\pm 1} \cup \{\kappa_0,\kappa_1, \kappa_2 \}^{\pm 1}$
 \ is not the identity then there exists a word $x \in \{0,1,\#\}^*$ of 
length  $|x| \leq c \, |w|$ (for some constant $c$), such that $w(x)$ is
defined and $w(x) \neq x$.

Therefore, a nondeterministic algorithm for the negated word problem of 
$\langle G_{3,1} \cup \{\kappa_0,\kappa_1, \kappa_2 \} \rangle$ 
simply needs to guess $x$, then compute $w(x)$, then check that $x \neq w(x)$. 
Guessing $x$ takes linear time (since $|x| \leq c \, |w|$). Applying an 
element $\delta \in \Delta_{3,1}^{\pm 1}$ to a word $z \in \{0,1,\#\}^*$ takes 
constant time (since $\delta$ just changes a bounded-length prefix of $z$), 
and changes the length of $z$ by an additive constant: 
$|\delta(z)| \leq |z| + c$.
Applying $\kappa_i^{\pm 1}$ ($i = 0,1,2$) to $z$ will not change the length 
of $z$ and takes linear time ($\leq c \, |z|$). Finally, since 
$|w(x)| \leq c\, |w|$ for some constant $c$), one can check in linear time 
whether $x \neq w(x)$. So the Theorem will follow from the following Claim.

\medskip

\noindent {\sc Claim:} \ Let \ $w \in$
$(\Delta_{3,1}^{\pm 1} \cup \{\kappa_0,\kappa_1, \kappa_2 \}^{\pm 1})^*$ 
 \ be such that as an element of ${\mathcal G}_{3,1}$, $w$ is not the 
identity. Then there exists $x \in \{0,1,\#\}^*$ such that $w(x)$ is defined, 
$x \neq w(x)$, and $|x| \leq c \, |w|$.

\smallskip

\noindent {\sc Proof} of the Claim: \   
Let $\ell$ be the length of the longest word in the domain and image codes 
of the elements of $\Delta_{3,1}$.

The word $w$ is of the form \  
$w \ = \ g_n K_n g_{n-1} K_{n-1} \ \cdots \ g_1 K_1 g_0$, \ 
where $g_n, \ldots, g_1, g_0 \in (\Delta_{3,1}^{\pm 1})^*$,  and 
$K_n, \ldots, K_1 \in $ 
$(\{\kappa_0,\kappa_1, \kappa_2 \}^{\pm 1})^*$.
Since $w$ does not represent the identity, there exists a word 
$z \in \{0,1,\#\}^*$ such that $z \neq w(z)$. We can assume that $z$ is long
enough (indeed, $w(zZ) = w(z) \, Z \neq zZ$ for any word $Z \in \{0,1,\#\}^*$; 
so we could replace $z$ by $zZ$ and thus make $z$ as long as we wish).
So we can assume that $|z| > 3N$, where \ $N \ = \ $
$\ell \, \sum_{j=0}^n |g_j|  + 6 \, \sum_{j=1}^n |K_j|$ \  
($\leq (\ell+6) \, |w|$). 
Let $pqr$ be the prefix of length $3N$ of $z$, where $|p| = |q| = |r|$.
We will show that $x = pqr\#$ satisfies $w(x) \neq x$. 

The first (i.e., the right-most) generator in $g_0$ affects only the 
left-most $\ell$ letters of $z$. Since the right-most letter in $g_0$ could
shorten $z$ by up to $\ell - 1$, the right-most two letters of $g_0$ could 
affect at most the first $2 \ell$ letters of $z$. In total, $g_0$ can 
affect the left-most $\ell \, |g_0|$ (or fewer) letters of $z$.

Next, $K_1$ moves each bit of $g_0(z)$ over a distance $\leq 6 \, |K_1|$.  
So, $K_1 g_0$ changes the left-most $6 \, |K_1| + \ell \, |g_0|$ (or fewer) 
letters of $z$ (in ways that we will not try to specify). The letters further 
to the right in $g_0(z)$ (at positions $> 6 \, |K_1| + \ell \, |g_0|$) are 
permuted by $K_1$ iff $\#$  does not appear within the left-most 
$6 \, |K_1| + \ell \, |g_0|$ positions of $g_0(z)$.  Note that since 
$w(z)$ is defined, $g_0(z)$ must contain some $\#$ (otherwise, $K_1$ would 
not be defined on $g_0(z)$). 

For the same reason,  $w$ changes the left-most \ 
$N \ = \ 6 \sum_{j=1}^n |K_j| + \ell \sum_{j=0}^n |g_j|$ \ (or fewer) letters 
of $z$ in fairly arbitrary ways; those are the positions in the prefix $p$ of
$z$. The letters further to the right in $z$ (at positions \ $> N$) are 
only permuted according to some of the $K_m$'s ($n \geq m \geq 1$),
namely for those $m$ for which \, 
$g_{m-1} K_{m-1} \ldots g_1 K_1 g_0 (z)$ \, does not contain $\#$ within
the $N$ leftmost positions. Let 
$K \in \langle \kappa_0,\kappa_1, \kappa_2 \rangle$ be the 
concatenation of those $K_m$ ($m = n, \ldots, 1$) for which there is no $\#$ 
in \, $g_{m-1} K_{m-1} \ldots g_1 K_1 g_0 (z)$ \, within the $N$ leftmost 
positions. 

Since $w$ changes $z$, it either changes the prefix $p$ of $z$, and in that
case, $w$ will of course also change $x = pqr\#$. Or $w$ does not change the 
prefix $p$, but $K$ permutes bits at positions $> N$ in $z$, 
non-identically. Moreover, by Lemma \ref{subgroupKappa0123}, if $K$ acts 
non-identically at a position $i+3$ of $z$, with $N < i$, and 
$|z| \geq N + 4 \, |K|$, then $K$ also acts non-identically on position $i$ 
of $z$. Thus, acts non-identically on a position $i$ of $z$, with 
$N+3 \leq i > N$. Then $w$ changes $p$, hence $x$.
This proves the Claim, and hence the Theorem.
 \ \ \ $\Box$

\bigskip

The main theorem of the previous section now becomes:

%% Main Theorem %%  
\begin{thm}  \label{finPrescoNPcompl} \   
There exists a finitely presented Thompson group $G$
($\subset {\mathcal G}_{3,1}$), with the following properties:

\smallskip

\noindent $\bullet$ \ The word problem of $G$ (over a fixed finite
generating set) is coNP-complete (with respect to polynomial-time
constant-arity conjunctive reduction).

\smallskip

\noindent $\bullet$ \   
$G$ is an HNN extension (by one stable letter) of some finitely presented 
subgroup {\rm Th} of $G_{3,1}$. In fact, $G$ is
isomorphic to the semidirect product ${\rm Th} \rtimes {\mathbb Z}$.

\smallskip

An example of such group $G$ is the subgroup  \
$\langle G_{3,1}^{\rm mod \, 3}(0,1;\#) \cup \{\kappa_{321}\} \rangle$ \ of \
${\mathcal G}_{3,1}$, where {\rm Th} is $G_{3,1}^{\rm mod \, 3}(0,1;\#)$.
\end{thm}
{\bf Proof.} \ This follows directly by combining Theorems 
\ref{finprescoNPhard} and \ref{coNP}. 
 \ \ \ $\Box$ 

\bigskip

Next we give a coNP-completeness result about finitely generated 
{\it simple} groups. First, recall the following: 
{\it If $G_{N,1} \subseteq G \subseteq {\mathcal G}_{N,1}$ then the 
commutator subgroup $G'$ is a simple group} (see R.~Thompson's comment
before Corollary 1.11 in \cite{Th}; an actual proof of this claim and a 
generalization to the Thompson-Higman groups $G_{N,1}$ was given by 
E.~Scott, Lemma 20 in \cite{ESc}). Note the symbols ``$\subseteq$'' in the
result; it is not sufficient that $G$ contains a copy of $G_{N,1}$ and 
${\mathcal G}_{N,1}$ contains a copy of $G$, but the copy of $G_{N,1}$ inside
$G$ must be identical with the subgroup $G_{N,1}$ of ${\mathcal G}_{N,1}$. 

When $H$ is a subgroup of a group of $G$, recall the {\it 
Reidemeister-Schreier rewrite process} (see e.g., \cite{MagnusKaSo}
pp.~90-93, \cite{LyndonSchupp} pp.~102-104, \cite{Reidemeister} pp.~69-78). 
The graphical form of the process is quite intuitive. 
One first takes the {\it Schreier graph}, whose vertex set is the set of 
cosets $Hg_i$, where  $g_i$ $(i = 1, \ldots, k)$ are coset representatives (we only use the case when $k$ is finite). The set of (labeled)
edges of the Schreier graph is \  
$\{ Hg_i \stackrel{a}{\longrightarrow} Hg_ia \ : \ $
$a \in A, \ i = 1, \ldots, k \}$,     
where $A$ is a generating set of $G$. We will only consider the case when
$A$ is finite. 

Hence, when $G$ is finitely generated and $H$ has finite index in 
$G$ then the Schreier graph is a finite automaton (if we pick the coset $H$ 
as both start and accept state), which decides the generalized word problem 
of $H$ in $G$ (deterministically in linear time). We also have the following 
interesting properties, assuming $H$ has finite index in $G$: If $G$ is 
finitely generated then $H$ is finitely generated; if $G$ is finitely
presented then $H$ is finitely presented (see the above references).
Moreover, the Reidemeister-Schreier rewrite process shows that when $G$ is
generated and $H$ has finite index in $G$ then the distortion of $H$ in $G$ 
is linear.

Let us pick a spanning tree in the Schreier graph, with root $H$; this is 
the graphical way of choosing a Schreier transversal: for every vertex 
$Hg_i$ let $t_i \in (A^{\pm 1})^*$ be the label of the path in the spanning 
tree from the root $H$ to $Hg_i$; then the word $t_i$ represents an element 
of $Hg_i$, so we can write $Ht_i$ for $Hg_i$; let 
$T = \{t_i : i = 1, \ldots, k \}$.  For any word $w \in (A^{\pm 1})^*$, we
denote the coset representative of $w$ by $\overline w$ ($\in T$). 
The following set, called the {\it Reidemeister-Schreier generators}, 
generates $H$: \ \  
 $R = \{ t_i \, a \, (\overline{t_i a})^{-1} \ : \ t_ia \not\in T, $
$a \in A^{\pm 1}, \ i = 1, \ldots, k \}$. 

We need an auxiliary result:

\begin{pro} \label{finIndex} \ 
Suppose $G$ is a finitely generated group, and $H$ is a subgroup of $G$ of
finite index. Then the word problems of $G$ and $H$ are reducible to each
other by linear-time many-to-one reductions. 
\end{pro}
{\bf Proof.} \ Recall that by the Reidemeister-Schreier rewrite process,
$H$ is finitely generated; let $R$ be a finite generating set of $H$. 
Hence, the identity embedding of $H$ into $G$ is a one-to-one reduction of 
the word problem of $H$ to the word problem of $G$. The reduction just consists
of expressing each generator in $R$ by a string over $A$ in a fixed way, so 
this reduction has linear time complexity.

Conversely, let us reduce the word problem of $G$ to the word problem of
$H$. Let $h_0$ be some fixed word over the  such that 
$h_0 \neq {\bf 1}$ in $H$. A function that reduces the word problem of $G$ 
to the word problem of $H$ can be defined by

\smallskip

$w \in (A^{\pm 1})^* \ \longmapsto \ $
$\left\{\begin{array}{lll}
h_0 \ & \ \ \ {\rm if}  & w \not\in H,  \\
(w)_R \ & \ \ \ {\rm if}  &  w \in H.
\end{array}        \right. $ 

\smallskip

\noindent Here, $(w)_R$ denotes the expression of $w$ over the 
Reidemeister-Schreier generating set $R$ of $H$ (when $w \in H$). 
By the Reidemeister-Schreier rewrite process, $(w)_R$ can be obtained from
$w$ in linear time. Since the generalized word problem of $H$ in $G$ is
decidable in linear time (using the Schreier graph automaton), it follows 
that the above reduction function is computable in linear time. 
Finally, $w = {\bf 1}$ in $G$ iff $f(w) = {\bf 1}$ in $H$. 
 \ \ \ $\Box$

\begin{thm} \label{finGenSimple} \ 
There exists a finitely generated {\em simple} Thompson group whose word 
problem is coNP-complete (with respect to polynomial-time constant-arity 
conjunctive reduction).

An example of such a group is \   
$\langle G_{3,1}\cup \{\kappa_0, \kappa_1, \kappa_2 \} \rangle'$, \  
i.e., the commutator subgroup of 
$\langle G_{3,1}\cup \{\kappa_0, \kappa_1, \kappa_2 \} \rangle$.
\end{thm}
{\bf Proof.} \ The group  
$G = \langle G_{3,1}\cup \{\kappa_0, \kappa_1, \kappa_2 \} \rangle$
satisfies $G_{3,1} \subseteq G \subseteq {\mathcal G}_{3,1}$. We immediately
conclude that the commutator subgroup $G'$ is a {\em simple} group, by the
earlier remarks on R.~Thompson's comments.
Also, $G'$ has finite index in $G$. Indeed, 
$\kappa_0^3 = \kappa_1^3 = \kappa_2^3 = 1$, and $G_{3,1}'$ has
index 2 in $G_{3,1}$ (by \cite{Hig74}).
Clearly, $G$ is finitely generated (since the Thompson-Higman group $G_{3,1}$ 
is finitely generated). It follows that $G'$ is finitely generated, by our
remarks above on the Reidemeister-Schreier rewrite process. 

By Proposition \ref{finIndex}, the word problems of $G$ and $G'$ are 
reducible to each other. Hence, since the word problems of $G$ is 
coNP-complete, the word problems of $G'$ is also coNP-complete. 
 \ \ \ $\Box$

\bigskip

We have now completed the proofs of the main theorems, which give us 
finitely presented Thompson groups, and finitely generated simple Thompson
groups with coNP-complete word problems. To finish, let us give some more
explanations of the fact that the finitely presented group $G_{3,1}$ (over 
a finite set of generators) has a word problem in {\bf P} (deterministic 
polynomial time), but over an infinite set of generators (obtained by 
including all letter transpositions) the word problem of $G_{3,1}$ is 
coNP-complete.
This is related to the concept of {\it distortion}. See \cite{Gromov} 
for the original definition by Gromov, and \cite{Olsh}, \cite{OlSap}
for a slightly more natural definition and some interesting results;
results on distortion in Thompson groups appear in \cite{BiThomps}; 
the complexity version of the Higman embedding theorem (in \cite{Bi}
for semigroups and \cite{SBR}, \cite{BORS} for groups) show that the 
embeddings given there have linear distortion. 

Originally, Gromov only defined distortion to characterize the relation
between a group and a subgroup. In the present context it is useful to
also consider the self-distortion of a group, relative to different 
generating sets. In the case of $G_{3,1}$ we consider a finite generating
set $\Delta_{3,1}$ and the infinite generating set
$\Delta_{3,1} \cup \{ \tau_{i,i+1} : i \geq 0\}$ (with infinitely many
transpositions included in the generating set). The {\it self-distortion}
of $G_{3,1}$ relative to  
$B = \Delta_{3,1} \cup \{ \tau_{i,i+1} : i \geq 0\}$ is said to have 
upper bound $f$ iff $f: {\mathbb N} \to {\mathbb N}$ is a non-decreasing
function such that for every $g \in G_{3,1}$ we have: \   
$|g|_A \ \leq \ f(|g|_B)$.
Here, $|g|_A$ denotes the word-length of $g$ over the generating set $A$, 
i.e., the length of the shortest word in $(A^{\pm 1})^*$ representing $g$
(and similarly for $|g|_B$).  
The next Lemma shows that the self-distortion of $G_{3,1}$ for the above 
generating set is at least exponential. The self-distortion of $G_{3,1}$ 
for the above generating set is closely related to the Gromov distortion of
$G_{3,1}$ within $\langle G_{3,1} \cup \{\kappa_{321}\} \rangle$.

\noindent {\bf Definition.} \ {\it
Two functions \  $f_1, f_2 : {\mathbb N} \to {\mathbb N}$ \ are said to be 
{\em linearly equivalent} iff there exist positive constants 
$c_0, c_1, c_2, c_3, c_4$ such that for all $n \geq c_0$: \ \   
$f_1(n) \leq c_1 \, f_2(c_2n)$ \ and \ $f_2(n) \leq c_3 \, f_1(c_4n)$.

A function $f : {\mathbb N} \to {\mathbb N}$ is {\em at least exponential}
iff there is a constant $c > 1$ such that for infinitely many $n$: \ 
$f(n) > c^n$.
}

\begin{lem} \label{distortionVinTboolLoweBound} \
The self-distortion of $G_{3,1}$ relative to the generating set 
$\Delta_{3,1} \cup \{ \tau_{i,i+1} : i \geq 0\}$ is at least exponential.
The distortion of $G_{3,1}$ in 
$\langle G_{3,1} \cup \{\kappa_{321}\} \rangle$ is at least exponential. 
Similarly, the distortion of $G_{3,1}^{\rm mod \, 3}(0,1;\#)$ in 
$H(0,1;\#)$ is at least exponential.
\end{lem}
{\bf Proof.} \ We saw already in Lemma \ref{tau_generated} that the 
transposition $\tau_{n-1,n}$ (where $n \geq 0$) has word length \   
$|\tau_{n,n+1}| \leq 2n -1$ \ in 
$\langle G_{3,1} \cup \{\kappa_{321}\} \rangle$, 
over any generating set containing $\tau_{1,2}$ and $\kappa_{321}$.
Moreover, we will prove next that the transpositions have exponential 
table size; this implies an exponential word length over any fixed finite
generating set of $G_{3,1}$, as we will see. 

\smallskip

\noindent {\sc Claim.} \ The table-size of the maximum extension of
$\tau_{n-1,n}$ is  \ \ $\|\tau_{n-1,n}\| \ = \ 2^{n+2} - 1$.

\smallskip

\noindent Proof of the Claim: \
The domain and image code of $\tau_{n-1,n}$, as originally defined, are both
equal to \ $\{0,1\}^{n+1} \ \cup \ \{0,1\}^{\leq n} \# $.
However, in order to find $\|\tau_{n-1,n}\|$ we must maximally extend 
$\tau_{n-1,n}$. Recall that a bijection $\varphi$ between finite maximal 
prefix codes is extendable iff there exist $u, v \in \{0,1,\#\}^*$ such
that domC($\varphi$) contains the triple $u0, u1, u\#$ and imC($\varphi$)
contains the triple $v0,v1, v\#$ with $\varphi(u0) = v0$, 
$\varphi(u1) = v1$, $\varphi(u\#) = v\#$.

A triple of arguments in the domain code of $\tau_{n-1,n}$ that could possibly
lead to an extension is of the form 

\smallskip

$(x_0 \ldots x_{n-2} x_{n-1} 0, \ x_0 \ldots x_{n-2} x_{n-1} 1,$
$ \ x_0 \ldots x_{n-2} x_{n-1} \#)$,

\smallskip

\noindent where $x_0, \ldots, x_{n-2}, x_{n-1} \in \{0,1\}$.
The transposition $\tau_{n-1,n}$ maps this triple to the triple \

\smallskip

$(x_0 \ldots x_{n-2} 0 x_{n-1}, \ x_0 \ldots x_{n-2} 1 x_{n-1}, $
$ \ x_0 \ldots x_{n-2} x_{n-1} \#)$.

\smallskip

\noindent So, whether $x_{n-1} = 1$ or $x_{n-1} = 0$, no extension is 
possible. The set \ $\{0,1\}^{n+1} \cup \{0,1\}^{\leq n} \# $ \ has 
cardinality \ $2^{n+1} + 2^{n+1} -1$.
This proves the Claim.

\bigskip

\noindent Now, by the relation \
$c_{_{\Delta}} \cdot \|\tau_{n,n+1}\| \leq |\tau_{n,n+1}|_{_{\Delta}}$ \
(Corollary 4.7 in \cite{BiThomps}) for some constant $c_{_{\Delta}} > 0$, 
depending on the choice of a finite generating set $\Delta_{3,1}$ chosen 
for $G_{3,1}$:

\smallskip

$|\tau_{n-1,n}|_{_{\Delta}} \ \geq \ c_{_{\Delta}} \cdot \|\tau_{n-1,n}\|$
$  \ \geq \ c_{_{\Delta}} \cdot 2^{n+2} - c_{_{\Delta}}$

\smallskip

\noindent
Hence, since \ $2n-1 \geq |\tau_{n-1,n}|_{_{\Delta \cup \kappa}}$ \
(as we already saw at the beginning of this proof),

\smallskip

$|\tau_{n-1,n}|_{_{\Delta}} \ \geq \ c_{_{\Delta}} \cdot $
$2^{\frac{1}{2} \, |\tau_{n-1,n}|_{_{\Delta \cup \kappa}} \ + \ \frac{5}{2}}$
$ - c_{_{\Delta}}$.  
 \ \ \ \ \ \ $\Box$

\bigskip

\begin{thm} \label{distortionVinTbool} \
The distortion of $G_{3,1}$ in 
$\langle G_{3,1} \cup \{\kappa_{321}\} \rangle$
is exponential (i.e., it is linearly equivalent to $2^n$). 

Similarly, the distortion of $G_{3,1}^{\rm mod \, 3}(0,1)$ in $H(0,1)$ or in 
$\langle G_{3,1} \cup \{\kappa_{321}\} \rangle$ is exponential.
And the distortion of $G_{3,1}^{\rm mod \, 3}(0,1;\#)$ in $H(0,1;\#)$ or in 
$H(0,1)$ or in $\langle G_{3,1} \cup \{\kappa_{321}\} \rangle$ is 
exponential.
\end{thm}
{\bf Proof.} \ We already saw an exponential lower bound, in Lemma
\ref{distortionVinTboolLoweBound}. We will now prove an exponential upper 
bound, of the form $c^n$ (for some constant $c > 1$). 

Let $\Delta_{3,1}$ be a fixed finite generating set for $G_{3,1}$. For 
$\varphi \in G_{3,1}$, let $|\varphi|_{\Delta_{3,1}}$ denote the word-length 
of $\varphi$ over the generating set $\Delta$ (i.e., the length of a shortest
word over $\Delta_{3,1}^{\pm 1}$ that represents $\varphi$). 
Similarly, $|\varphi|_{\Delta_{3,1}, \kappa}$ denotes the word-length of
$\varphi$ over $\Delta_{3,1} \cup \{\kappa\}$ (i.e., the length of a shortest
word over $\Delta_{3,1}^{\pm 1} \cup \{\kappa^{\pm 1}\}$ that represents 
$\varphi$).

\smallskip

\noindent {\sc Claim:} \ Let $w$ be a word over 
$\Delta_{3,1}^{\pm 1} \cup \{\kappa^{\pm 1}\}$ that represents an element 
$\varphi$ of $G_{3,1}$, and assume that $w$ is in shortest form (i.e., there 
is no shorter word over $\Delta_{3,1}^{\pm 1} \cup \{\kappa^{\pm 1}\}$, 
representing the same group element). Then the longest entry in the table of 
$\varphi$ has length \, $\leq (6 + \ell) \, |w|$. 

\smallskip

\noindent Proof of the Claim: \ Let \ 
$w \ = \  \ g_n \kappa^{i_n} g_{n-1} \kappa^{i_{n-1}} \ \cdots \ $
       $ g_1 \kappa^{i_1} g_0$, \
where $i_n, \cdots, i_1 \in {\mathbb Z} - \{0\}$, and
$g_n, \cdots, g_1, g_0 \in (\Delta_{3,1}^{\pm 1})^*$.

As in the proof of Theorem \ref{coNP}, let $x \in \{0,1,\#\}^*$ be any 
word of length at least $3N$, where  \  
$N = \ell \, \sum_{j=0}^n |g_j|  + 6 \, \sum_{j=1}^n |i_j|$ \
($\leq (\ell+6) \, |w|$), and where $\ell$ is the length of the longest word 
in the domain and image codes of the elements of $\Delta_{3,1}$. 
In the proof of Theorem \ref{coNP} we saw that the action of $w$ on $x$
changes the left-most $N$ (or fewer) letters of $x$ in fairly arbitrary ways.
The letters of $x$ at positions further to the right (i.e., at positions
$> N$) are only permuted according to $\kappa^{i_{\rm sum}(x)}$.
We have $i_{\rm sum}(x) = 0$, otherwise $w$ would change bits at 
arbitrarily remote positions on $x$ (for arbitrarily long words $x$; this 
would imply that $w$ has an infinite table (contradicting the assumption that
$w$ represents an element of $G_{3,1}$). 

Now, since $i_{\rm sum}(x) = 0$, $w$ only changes letters at positions 
$\leq N$ ($\leq (6 + \ell) \, |w|$) in $x$. Therefore, the longest word in
the domain code of $w$ has length $\leq (6 + \ell) \, |w|$.
This proves the Claim.

\medskip

\noindent It follows immediately from the Claim that the table size of 
$\varphi$ satisfies \ 
$\|\varphi\| \ \leq \ 3^{(6 + \ell) \, |w|}$. 
 \ Note that here, $w$ is the word length of $\varphi$ over the
generating set $\Delta_{3,1} \cup \{\kappa\}$; i.e., \ 
$|w| = |\varphi|_{\Delta_{3,1}, \kappa}$.
By Theorem 4.8 in \cite{BiThomps}, \ \ $|\varphi|_{_{\Delta_{3,1}}} \ \leq \ $
$c_{_{\Delta}} \cdot \|\varphi\| \cdot \log_2 \|\varphi\|$, \ where 
$c_{_{\Delta}} > 0$ is a constant. Hence, 

\smallskip

$|\varphi|_{_{\Delta}} \ \leq \ $
$c \ 3^{c \, |\varphi|_{\Delta, \kappa}} \ c \ $
$|\varphi|_{\Delta, \kappa}$
$ \ \leq \ C^{|\varphi|_{\Delta, \kappa} }$. 

\smallskip

\noindent for some constants $c, C > 1$. This proves the Theorem.
  \ \ \ $\Box$ 

\bigskip

%%%%%%%%%%%%%%%%%%%%%%%%%%%%%%%%%%%%%%%%%%%%%%%%%%%%%%%%
% Section 9 
%%%%%%%%%%%%%%%%%%%%%%%%%%%%%%%%%%%%%%%%%%%%%%%%%%%%%%%%

\newpage

\section{Appendix}

%%%%%%%%%%%%%%%%%%%%%%%%%%%%%%%%%%%%%%%%%%%%%%%%%%%%%%%%

\subsection{Properties of prefix codes}  % A1

In this appendix we prove various properties of prefix codes that are 
used in the paper. Recall that $\varepsilon$ denotes the empty word.

\begin{lem} {\bf (Lemma \ref{endmarker_code})} \\  
{\bf (1)} \
If \ $P \subset \{0,1\}^* \cup \{0,1\}^* \#$ \ is a maximal prefix code
over $\{0,1,\# \}$ then \ $P = P_1 \, \cup \, P_2 \#$ \ for some
$P_1, P_2 \subset \{0,1\}^*$, with the following properties:

\smallskip

$\bullet$  \ $P_1$ is a maximal prefix code over $\{0,1\}$;

\smallskip

$\bullet$ \ $P_2 \ = \ \{ p \in \{0,1\}^* \ : \ $
       $p$ is a strict prefix of some element of $P_1  \}$.

\smallskip

\noindent When $P_1$ is finite, this last property implies: \
 $|P_2| = |P_1| - 1$.

\smallskip

\noindent {\bf (2)} \
 Conversely, if \ $P \ = \ P_1 \ \cup \ P_2 \#$ \ for some
$P_1, P_2 \subset \{0,1\}^*$ with the above two properties, then $P$ is a
maximal prefix code over $\{0,1,\# \}$.
\end{lem}
{\bf Proof.} \ If $P \subset \{0,1\}^* \ \cup \ \{0,1\}^* \# $ \ is
a maximal prefix code then $P$ has the form
 \ $P = P_1 \ \cup \ P_2 \#$, with $P_1, P_2 \subset \{0,1\}^*$.
Since $P$ is a maximal prefix code,
$P_1$ is a maximal prefix code over $\{0,1\}$. Also, the set
$P_2 \#$ \, is a prefix code for any subset $P_2 \subset \{0,1\}^*$ (since
any two elements $p_2\# \neq p_3\#$ with
$p_2, p_3 \in \{0,1\}^*$ are prefix incomparable).

Let us prove that $P_2$ is as in the Lemma. Since $P_1$ is a maximal prefix
code, every $p_2 \in P_2$ (and in fact every string in $\{0,1\}^*$) is
prefix comparable with some element of $P_1$. Let's say, $p_1 \in P_1$
is prefix comparable with $p_2 \in P_2$. If $p_1$ were a prefix of $p_2$
then $p_1$ would also be a strict prefix of $p_2\#$, which would contradict
the fact that $P$ is a prefix code. This shows that every element of $P_2$
is strict prefix of an element of $P_1$. Since $P$ is a maximal prefix code,
$P_2$ consists of all strict prefixes of elements of $P_1$.

It is straightforward to prove the converse, namely that every set
$P = P_1 \cup P_2 \#$, with $P_1, P_2$ as above, is a maximal prefix code.
 \ \ \ $\Box$

\begin{lem} \label{endmarker_code2} \
Every maximal prefix code over the alphabet $\{0,1,\#\}$ can be written in the
form \ $P_1 \ \cup \ \bigcup_{v \in P_2} v  \# \, P(v)$, for some
$P_1, P_2 \subset \{0,1\}^*$ and $P(v) \subset \{0,1,\# \}^*$,
with the following properties:

\smallskip

$\bullet$  \ $P_1$ is a maximal prefix code over $\{0,1\}$;

\smallskip

$\bullet$ \ $P_2 \ = \ \{ p \in \{0,1\}^* : \ p$ is a strict prefix of some
                      element of $P_1  \}$.

\smallskip

 \ \ \ \ Hence, when $P_1$ is finite, this last property implies: \
 $|P_2| = |P_1| - 1$.

\smallskip

$\bullet$  \ For every $v \in P_2$, the set $P(v)$ is a maximal prefix code
over $\{0,1,\#\}$.

\smallskip

 Conversely, if \ $P \ = \ P_1 \ \cup \ \bigcup_{v \in P_2} v \# \, P(v)$ \
for some $P_1, P_2$ and $P(v)$ with the above three properties, then $P$ is a
maximal prefix code over $\{0,1,\# \}$.
\end{lem}
{\bf Proof.} \ The proof is straightforward, and similar to the proof of
Lemma \ref{endmarker_code}. \ \ \  $\Box$

\begin{lem} \label{pref_comp} \
Let $x, y, u, v \in A^*$. If $xu$ and $yv$ are prefix-comparable
then $x$ and $y$ are prefix-comparable. Contrapositively,
if $x$ and $y$ are prefix incomparable then $xu$ and $yv$ are prefix
incomparable.
\end{lem}
{\bf Proof.} \ The proof is straightforward. \ \ \ $\Box$

\begin{lem} \label{x_inv_pref_code} \
If $P \subset A^*$ is a maximal prefix code and if $x \in A^*$,
then \ $x^{-1} P = \{w \in A^* : xw \in P \}$ \
is either empty or a maximal prefix code. Equivalently, if $P A^*$ is an
essential right ideal then $x^{-1} P A^*$ is also an essential right
ideal.
\end{lem}
{\bf Proof.} \ This is a classical property of maximal prefix codes
(see e.g. Lemma 8.4 in \cite{BiThomps}).    \ \ \ $\Box$

\begin{lem} \label{pref_code_length} \
If a maximal prefix code over an alphabet $A$, $|A| \geq 2$, contains a word 
of length $n$ then this prefix code has cardinality at least $n+1$. 
\end{lem}
{\bf Proof.} \ This is a classical property of maximal prefix codes
(see e.g. \cite{BiThomps}).    \ \ \ $\Box$

\begin{lem} {\bf (Lemma \ref{non_id})} \\  
Let $P \subset A^*$ be any prefix code, where $|A| = n \geq 2$.  Assume
$\varphi \in {\rm pStab}_{{\mathcal G}_{n,1}}(P A^*)$, but
$\varphi \not\in {\rm pFix}_{{\mathcal G}_{n,1}}(P A^*)$.
Then there exists $x \in P A^*$ such that $x$ and
$\varphi(x)$ are not prefix-comparable.

In particular, if $\varphi \in {\mathcal G}_{n,1}$ is not the identity element
then there exists $x \in {\rm domC}(\varphi)$ such that $x$ and
$\varphi(x)$ are not prefix-comparable.
\end{lem}
{\bf Proof.} \
Let $P'$ be another prefix code such that $P$ and $P'$ complementary prefix
codes \ ($P' = \emptyset$ if $P$ is a maximal prefix code).
Let $\psi$ be the restriction of $\varphi$ to the essential
right ideal $(P \cup P') A^*$. So $\psi$ is a right-ideal isomorphism
that represents $\varphi$.

We will prove the contrapositive of the Lemma:
{\it  Assume that $x$ and $\psi(x)$ are prefix-comparable for all
$x \in {\rm  domC}(\psi) \cap P A^*$, and that $\psi$ and $\psi^{-1}$
stabilize $P A^*$ where $\psi$ and $\psi^{-1}$ are defined. Then the 
restriction of $\psi$ to ${\rm  domC}(\psi) \cap P A^*$ is the identity 
map. }

\smallskip

\noindent Case 1: \ $x >_{\rm pref} \psi(x)$ \ (i.e., $x$ is a strict prefix
of $\psi(x)$), for some $x \in {\rm domC}(\psi) \cap P A^*$.

Then $\psi(x) = xv$ for some $v \in AA^*$, and
$\psi(x) = xv \in {\rm imC}(\psi) \cap xAA^*$.
By Lemma \ref{x_inv_pref_code}, $x^{-1}({\rm imC}(\psi) \cap xAA^*)$ \ is a
maximal prefix code. Since it contains the non-empty word $v$, it contains at
least two elements (by Lemma \ref{pref_code_length}).
Hence ${\rm imC}(\psi) \cap xAA^*$ contains at least two elements.  
So there exists $x' \in {\rm  domC}(\psi)$ such that
$\psi(x') \neq \psi(x)$ and
$\psi(x') = xw \in {\rm imC}(\psi) \cap xAA^*$ (for some $w \in AA^*$).
Since $\psi^{\pm 1}$ stabilizes $P A^*$ we also have $x' \in P A^*$.
Since ${\rm imC}(\psi)$ is a prefix code, the inequality
$\psi(x') \neq \psi(x)$ implies that $\psi(x')$ and $\psi(x)$ are not
prefix-comparable.

By the assumption of the Lemma (or its contrapositive): $x'$ and
$\psi(x')$ are prefix-comparable. Hence we have two possibilities:

\smallskip

\noindent (1) \ $x' \leq_{\rm pref} \psi(x')$: \ then
$x' \leq_{\rm pref} \psi(x') <_{\rm pref} x$, so $x' \leq_{\rm pref} x$,
which contradicts the fact that ${\rm domC}(\psi)$ is a prefix code.

\smallskip

\noindent (2) \ $x' >_{\rm pref} \psi(x')$: \ then
$x' >_{\rm pref} \psi(x') = x'z$ for some $z \in AA^*$, and
$\psi(x') = xw$ (as seen above). Hence $x'z = xw$, which implies that
$x'$ and $x$ are prefix-comparable; again, this contradicts the fact that
${\rm domC}(\psi)$ is a prefix code.

\smallskip

We conclude that case 1 is impossible.

\smallskip

\noindent Case 2: \ $x <_{\rm pref} \psi(x)$ \, for some
$x \in {\rm domC}(\psi) \cap P A^*$.

Then $x = \psi(x) \ u$ \ for some $u \in AA^*$; so,
$x \in {\rm domC}(\psi) \cap \psi(x) \, AA^*)$.
Moreover, by Lemma \ref{x_inv_pref_code}, \
$\psi(x)^{-1}({\rm domC}(\psi) \cap \psi(x) \, AA^*)$ \
is a maximal prefix code. Since it contains the non-empty word $u$,
it contains at least two elements (by Lemma \ref{pref_code_length}); hence,
${\rm domC}(\psi) \cap \psi(x) \, AA^*$ contains at least two
elements. Therefore, there is $x' \neq x$ with
$x' \in {\rm domC}(\psi) \cap \psi(x) \, AA^*$. Moreover, since
$\psi^{\pm 1}$ stabilizes $P A^*$, we have $x' \in P A^*$.

By the assumption of the Lemma (or its contrapositive): $x'$ and
$\psi(x')$ are prefix-comparable. Hence we have two possibilities:

\smallskip

\noindent (1) \ $x' \geq_{\rm pref} \psi(x')$: then
$\psi(x') \leq_{\rm pref}  x'$ and $x' <_{\rm pref} \psi(x)$ (since
$x' \in \psi(x) \, AA^*$).
Therefore, $\psi(x') <_{\rm pref} \psi(x)$, which contradicts the fact
that ${\rm imC}(\psi)$ is a prefix code.

\smallskip

\noindent (2) \ $x' <_{\rm pref} \psi(x')$: then
$\psi(x') >_{\rm pref} x' = \psi(x') \ z$ \ for some $z \in AA^*$;
and $x' = \psi(x) \ w$ for some $w \in AA^*$
(since $x' \in \psi(x) \, AA^*$).
Therefore, $x' = \psi(x') \ z = \psi(x) \ w$, which implies that
$\psi(x')$ and $\psi(x)$ are prefix-comparable; again, this contradicts
the fact that ${\rm imC}(\psi)$ is a prefix code.

\smallskip

We conclude that case 2 is impossible.  Now, having ruled out cases 1
and 2, the only remaining possibility is that $x = \psi(x)$.
 \ \ \ $\Box$

\bigskip

It is well known (and easy to prove) that every finite maximal prefix code
$P$ over an alphabet $A$ has cardinality \ $|P| = 1 + (|A| -1) \, i_P$,
where $i_P$ is the number of internal vertices of the prefix tree of $P$.
Also, for every integer $i \geq 0$,
there exists a maximal prefix code $P$ over an alphabet $A$ of cardinality \
$1 + (|A| -1) \, i$.

\begin{lem} {\bf (Lemma \ref{extending_to_max_pref_code})} \\
Suppose $P, P' \subset A^*$ are complementary finite prefix codes.
Let $x_1, \ldots, x_k \in P A^*$, and assume $x_1$, $\ldots$, $x_k$ are
two-by-two prefix-incomparable.
Then for all $n$ of the form \ $n = 1 + i \, (|A|-1)$, with \
$n \geq \ |P| - k + (|A|-1) \, (|x_1| + \ldots + |x_k|)$,
there exists a prefix code $Q$  such that

\smallskip

$\bullet$ \ \ $Q \cup \{x_1, \ldots, x_k\}$ and $P'$ are complementary
prefix codes, with \ $Q \cup \{x_1, \ldots, x_k \} \subset P A^*$;

\smallskip

$\bullet$ \ \ $|Q| = n$.

\smallskip

$\bullet$ \ \ The set of prefixes of $P$ is a subset of the set of prefixes of
$Q \cup \{x_1, \ldots, x_k\}$.
\end{lem}
{\bf Proof.} \ Since $x_i \in P A^*$ we can write $x_i = p_i u_i$ for
some $p_i \in P$ and $u_i \in A^*$ ($i = 1, \ldots, k$).
Moreover, $p_i$ and $u_i$ are uniquely determined by $x_i$ (since $P$ is a
prefix code).

We will prove the Lemma only when $k = 2$, for clarity; the proof for general
$k$ is very similar.

Let $u_1 = c_1 c_2 \ldots c_m$ and $u_2 = d_1 d_2 \ldots d_n$.
In this Lemma we define \ ${\overline a} = A - \{a\}$, for all $a \in A$.
We define

\smallskip

$Q_0 \ = \ (P - \{p_1, p_2\}) \ \ \cup \ $
$p_1 \, ({\overline c_1} \cup c_1 {\overline c_2} \cup c_1c_2 {\overline c_3}$
  $ \cup \ldots \ \ldots \cup c_1c_2 \ldots c_{m-1} {\overline c_m} )$

\smallskip

 \ \ \ \ \  \ \ \ \ \  \ \ \ \ \    \ \ \ \ \  \ \ \ \ \   \ \ \ \ \
 $\cup \ $
$p_2 \, ({\overline d_1} \cup d_1 {\overline d_2} \cup d_1d_2 {\overline d_3}$
  $ \cup \ldots \ \ldots \cup d_1d_2 \ldots d_{n-1} {\overline d_n} )$.

\smallskip

\noindent In the special case where $x_1 = p_1$, we let

\smallskip

$Q_0 \ = \ (P - \{p_2\}) \ \ \cup \ $
$p_2 \, ({\overline d_1} \cup d_1 {\overline d_2} \cup d_1d_2 {\overline d_3}$
  $ \cup \ldots \ \ldots \cup d_1d_2 \ldots d_{n-1} {\overline d_n} )$.

\smallskip

\noindent and similarly if $x_2 = p_2$. Moreover, if both $x_1 = p_1$
and $x_2 = p_2$, we simply let $Q_0  = P$.

The formula for $Q_0$ implies directly that \
$|Q_0| = |P| - 2 + (|u_1| + |u_2|)(|A|-1) \leq $
  $|P| - 2 + (|x_1| + |x_2|)(|A|-1)$.

Intuitively, we can picture the set \
$p_1 \, ({\overline c_1} \cup c_1 {\overline c_2} \cup c_1c_2 {\overline c_3}$
$\cup \ \ldots \ \cup c_1c_2 \ldots c_{m-1}{\overline c_m}) \ \cup \ \{x_1\}$
 \ on the prefix tree of $A^*$: \
Consider the path labeled by $p_1 u_1 = x_1$, and consider all the vertices
{\it attached} to this path, but not {\it on} the path. The set above
consists of these attached
vertices, but excluding the prefixes of $p_1$ (i.e., we start after $p_1$),
and excluding the leaves of $x_1$. E.g., if $A = \{a,b\}$, and
$x_1 = p aaba$, then  the set is
$\{pb, pab, paaa, paabb, paaba \}$.
For $x_2$, the intuition is similar.

From the prefix tree picture it is obvious that
$Q_0 \cup \{x_1, x_2\} \cup P'$ is a maximal prefix code, since it consists
of the set of leaves of a subtree of $A^*$.

Finally, if we want a prefix code $Q$ as in the Lemma, with cardinality
exactly $n$, we simply take $Q_0$ and repeatedly replace some leaf in the
set \
$p_1 \, \{ {\overline c_1}, c_1 {\overline c_2}, c_1c_2 {\overline c_3}, $
       $ \ldots , c_1c_2 \ldots c_{m-1} {\overline c_m}\}$ \
by {\it its} $|A|$ children; each such step increases $|Q|$ by  $|A|-1$, while
preserving the fact that $Q \cup \{x_1, x_2\} \cup P'$ is a maximal prefix
code, and $Q \subset P A^*$.

\medskip

The fact that the set of prefixes of $P$ is a subset of the set of prefixes 
of $Q \cup \{x_1, \ldots, x_k\}$ follows immediately from the fact that
$Q \cup \{x_1, \ldots, x_k \} \subset P A^*$
 \ \ \ $\Box$

\begin{lem} {\bf (Lemma \ref{making_pref_codes}).} \\  
{\bf (1)} \ For all $x, y \in \{0,1\}^*$ there exist letters 
$\ell_1, \ell_2 \in \{0,1\}$
such that $x \ell_1$, and $y \ell_2$ are prefix incomparable.

\noindent {\bf (2)} \
For all $x, y, z \in \{0,1\}^*$ there exist letters $\ell_1, \ldots, \ell_6$
$\in \{0,1\}$ such that $x \ell_1 \ell_3$, $y \ell_2 \ell_4$, and
$z \ell_5 \ell_6$, are prefix incomparable.
\end{lem}
{\bf Proof.} \ (1) \ If $x, y$ are prefix incomparable, then any
$\ell_1, \ell_2 \in \{0,1\}$ will work, by Lemma \ref{pref_comp}.

If $x = y$,then $x0$ and $y1$ are prefix incomparable.

So let's suppose $x$ is a strict prefix of $y$. Then either $x0$  is a
prefix of $y$ (and then $x1$ is prefix incomparable with $y$ and $y\ell_2$),
or $x1$ is a prefix of $y$ (and then $x0$ is prefix incomparable with $y$
and $y\ell_2$).

The case where $y$ is a strict prefix of $x$ is very similar to the
previous case.

\smallskip

\noindent (2) \ By (1) there are letters $\ell_1, \ell_2$
$\in \{0,1\}$ such that $x \ell_1$ and $y \ell_2$ are prefix incomparable.
Then by Lemma \ref{pref_comp}, $x \ell_1 \ell_3$, $y \ell_2 \ell_4$
are prefix incomparable, for any $\ell_3, \ell_4$.

If $z$ is prefix incomparable with both $x \ell_1$ and $y \ell_2$
then any choice of $\ell_3$, $\ldots$, $\ell_6$ will work, by
Lemma \ref{pref_comp}.

If $z$ is a common prefix of both $x \ell_1$ and $y \ell_2$ then
either $z0$ or $z1$ will be prefix incomparable with both $x \ell_1$
and $y \ell_2$; then any choice of $\ell_3$, $\ell_4$, and $ell_6$
will work.

The only remaining cases are when $z$ is prefix incomparable with
exactly one of $x \ell_1$ and $y \ell_2$. Let's say, $z$ is prefix
incomparable with $y \ell_2$. Then by (1), $z\ell_5$ and $x \ell_1 \ell_3$
are prefix incomparable for some $ell_5$, $\ell_3$ $\in \{0,1\}$.
Then by Lemma \ref{pref_comp}, $y \ell_2 \ell_4$ and $z\ell_5$ are also
prefix incomparable, as well as $y \ell_2 \ell_4$ and $z\ell_5\ell_6$.
       \ \ \ $\Box$

\begin{lem} {\bf (Lemma \ref{V_endm_two_leaves}).}  \\   
{\bf (0)} \ Every finite maximal prefix code $P$ over an alphabet $A$
(e.g., $A = \{0,1,\#\}$) has cardinality \ $|P| = 1 + (|A| -1) \, i_P$,
where $i_P$ is the number of inner vertices of the prefix tree of $P$.

If $|P| > 1$ then $P$ contains a subset
of the form \ $uA$ (for some word $u \in A^*$).

Also, for every integer $i \geq 0$,
there exists a maximal prefix code $P$ over an alphabet $A$ of cardinality \
$1 + (|A| -1) \, i$.

\medskip

\noindent {\bf (1)} \ If \ $P \subset \{0,1\}^* \, \{\varepsilon,\# \}$,
and $|P| > 1$, then $P$ contains a subset of the form \

\smallskip

 \ \ \ \ \  $u \ \{0,1,\# \}$, \ \ \  for some $u \in \{0,1\}^*$

\medskip

\noindent {\bf (2)} \ For every integer $i \geq 3$ there is a maximal
prefix code \ $P \subset \{0,1\}^* \, \{\varepsilon, \# \}$, with
\ $|P| = 1 + 2i$, and with the following property:

\smallskip

 $P$ contains a subset of the form \ $\{u, v\} \, \{0,1,\# \}$,
  \ \  for some $u, v \in \{0,1\}^*$, \ $u \neq v$.

\medskip

\noindent {\bf (3)} \ For every integer $i \geq 5$, there is a maximal
prefix code \ $P \subset \{0,1\}^* \, \{\varepsilon, \#\}$, with
 \ $|P| = 1 + 2i$, and with the following property:

\smallskip

 $P$ contains a subset of the form \ $\{u, v, w\} \, \{0,1,\# \}$,
  \ \  for some $u, v, w \in \{0,1\}^*$,

  with $u, v, w$ distinct two-by-two.
\end{lem}
{\bf Proof.} \ Property {\bf (0)} is well known (see e.g., \cite{Hig74},
\cite{Bi}).

For property {\bf (1)}, recall Lemma \ref{endmarker_code}
about maximal prefix codes  \ $\subset \{0,1\}^* \, \{\varepsilon, \# \}$.
If $|P| > 1$ then $|P_1| > 1$. Any maximal prefix code $P_1$ over $\{0,1\}$
with $|P_1| > 1$ contains a subset of the form $u \ \{0,1\}$ (for some
$u \in \{0,1\}^*$). Hence $u \in P_2$, and $P$ will also contain \
$u \#$.

For property {\bf (2)},
let $P_1$ be a maximal prefix code over $\{0,1\}$, such that \
$\{u, v\} \ \{0,1\} \ \subset \ P_1$, for some words $u, v \in \{0,1\}^*$
 \ with $u \neq v$. For any $n \geq 4$, such a $P_1$ exists with cardinality
 \  $|P_1| = n$. \
This is folklore knowledge on prefix codes. One can prove it, e.g., by
looking at the inner tree of the tree of a maximal prefix code. One takes
an inner tree with 2 leaves (it suffices for the inner tree to have 3
vertices, arranged in the shape \verb+ /\ +). Then among the vertices
of the tree of the maximal prefix code there will be 2 vertices, each of
which has 2 leaves.

Recall that $P_1$ determines $P_2$ (and $|P_2| = |P_1| -1 = n-1$), thus for
any $n \geq 4$ we obtain a maximal prefix code \, $P = P_1 \cup P_2 \#$ \,
of cardinality \ $|P| = n + n-1 = 1 + 2i$, \  with $i = n-1 \geq 3$.

Property {\bf (3)} is proved in a similar way as property (2). We take an
inner tree with 3 leaves (it suffices for the inner tree to have 5
vertices), in the shape
\begin{verbatim}
                           /\
                          /\
\end{verbatim}
This proves the Lemma.    \ \ \ $\Box$

\begin{lem} {\bf (Lemma \ref{mod_3_code}).} \\   
$\bullet$ \ Suppose that there exists a maximal prefix
code $Q$ over the alphabet $\{0,1\}$, whose inner tree $T_{\rm in}(Q)$ has
two one-child vertices at depths $\equiv i$ {\rm mod 3} (for some
$i \in \{0,1,2\}$).
Then there exists a maximal prefix code $P \subset \{0,1\}^*$ with the same
{\rm mod 3} cardinality as $Q$, and with the following property:

\smallskip

 there is a word $u \in \{0,1\}^*$ such that \
 $u \cdot \{0,1\} \subseteq P$ \ and \ $|u| \equiv i$ {\rm mod 3}.

\smallskip

\noindent
Equivalently, the inner tree of the prefix code $P$ has a leaf at
depth $\equiv i$ {\rm mod 3}.

\medskip

\noindent  $\bullet$ \ More generally, let $k \geq 2$, let $i_1, \ldots, i_k$
$\in \{0,1,2\}$, and suppose that $T_{\rm in}(Q)$ has the following property:
For every $\lambda$ ($1 \leq \lambda \leq k$), \ $T_{\rm in}(Q)$ has a leaf of
depth $\equiv i_{\lambda}$ or it has two one-child
vertices at depths $\equiv i_{\lambda}$ {\rm mod 3}.

Then there exists a maximal prefix code $P \subset \{0,1\}^*$ with the same
{\rm mod 3} cardinality as $Q$, and with the following property:

\smallskip

there are $k$ different words $u_1, \ldots, u_k$ $\in \{0,1\}^*$ such that \
$\{u_1, \ldots, u_k \} \cdot \{0,1\} \subseteq P$ \ and \

\smallskip

$|u_1| \equiv i_1$, \ \ldots \ , \  $|u_k| \equiv i_k$,  \  {\rm mod 3}.

\smallskip

\noindent Equivalently, the inner tree of the prefix code $P$ has at
least $k$ leaves that have depths respectively  $\equiv i_1$, $\ldots$,
$\equiv i_k$ {\rm mod 3}.
\end{lem}
{\bf Proof.} \ We start with the maximal prefix code $Q$ and we transform it
into a maximal prefix code $P$ that has the required properties. The
transformations consist of rearrangements of the existing vertices of
$T_{\rm in}(Q)$.

\smallskip

Let us look at two one-child vertices $A, B$ of $T_{\rm in}(Q)$, at depths
$\equiv i_1$ mod 3. If $A$ and $B$ are not on a common path from the root,
we transform $T_{\rm in}(Q)$ as follows:

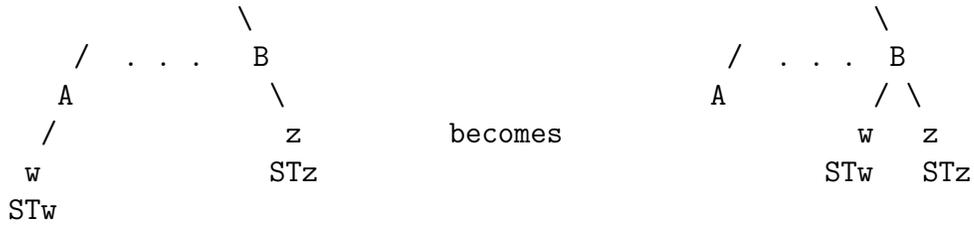
\begin{verbatim}
                     \                                      \
           /  . . .   B                            /  . . .  B
          A            \                          A         / \
         /              z         becomes                  w   z
        w              STz                               STw   STz
       STw
\end{verbatim}
Here, $w$ is the child of $A$ and $z$ is the child of $B$. Moreover, {\tt STw}
is the subtree with root $w$ and {\tt STz} is the subtree with root $z$.
We moved the subtree {\tt STw} with its root $w$ to the unoccupied child
position of $B$.
After the transformation, vertex $A$ is a leaf of depth $\equiv i_1$ mod 3.
Since the depths of $A$ and $B$ are equivalent mod 3, the above
transformation preserves the depth mod 3 of all nodes.

If $A$ and $B$ are on a common path from the root,
we transform $T_{\rm in}(Q)$ as follows:

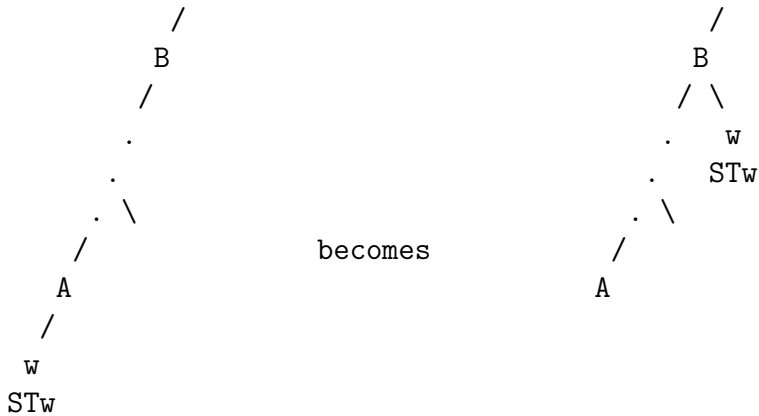
\begin{verbatim}
                /                                /
               B                                B
              /                                / \
             .                                .   w
            .                                .   STw
           . \                              . \
          /              becomes           /
         A                                A
        /
       w
      STw
\end{verbatim}
Again, after the transformation, vertex $A$ is a leaf of depth $\equiv i_1$
mod 3. The depths (mod 3) of all nodes are unchanged.

\smallskip

The proof in the general case is very similar, and can be done by induction.
If $T_{\rm in}$ already has some leaves (either present initially or obtained
by transformations as above) the additional transformations don't remove these
leaves and don't change their depths mod 3.
 \ \ \ $\Box$

%%%%%%%%%%%%%%%%%%%%%%%%%%%%%%%%%%%%%%%%%%%%%%%%%%%%%%%%

\newpage

% A2

\subsection{Commutation test for $G_{n,1}$ and $G_{n,1}^{\rm mod \, 3}$} 

We show in this subsection that the commutation test works for some fixators 
in $G_{n,1}$ and $G_{n,1}^{\rm mod \, 3}$.
The definitions and proofs are simpler that for the case of the group
$G_{3,1}^{\rm mod \, 3}(0,1;\#)$ studied in Section 5. 
The rest of this paper does not depend on this subsection.

For the finite alphabet $A$ below we assume $|A| \geq 2$. 

\begin{defn} \
Let $G \subseteq {\mathcal G}_{n,1}$ (i.e., $G$ is a subgroup with a
particular embedding), and let $P$, $P'$ be complementary prefix codes over
an alphabet $A$, with $|A| = n$.

The fixator ${\rm pFix}_G(P'A^*)$ is {\em separating} on $P A^*$ iff
for any ordered pair $(x,y)$ of prefix-incomparable words $x, y \in PA^*$ 
there exists $h \in {\rm pFix}_G(P' A^*)$ and there exists $u \in A^*$ 
such that

 \ \ \ \ \ \ $h(xu) = xu$ \ and \ $h(yu) \neq yu$.
\end{defn}

\begin{lem} \label{separ_maxA2} \
Let $G \subseteq {\mathcal G}_{n,1}$ and let $P$, $P'$ be complementary
prefix codes over an alphabet $A$.  If a fixator ${\rm pFix}_G(P'A^*)$ is 
separating on $P A^*$ then it is a maximal fixator.
\end{lem}
{\bf Proof.} \ Recall the definition of maximal fixator from Section 5.
Suppose by contradiction that for some
$y_0 \in PA^*$ we have for all $h \in {\rm pFix}_G(P'A^*)$: \ $h(y_0) = y_0$.
However, $x = y_0 a$ and $y = y_0 b$ are prefix-incomparable (for two
letters $a \neq b \in A$). Hence, by the separating property,
$h_0(yu_0) \neq yu_0$ for some $h_0 \in {\rm pFix}_G(P'A^*)$ and some
$u_0 \in A^*$.
On the other hand, $h_0(y_0) = y_0$ implies $h_0(yu_0) = h_0(y_0bu_0) = $
$h_0(y_0) \, bu_0$ $ = y_0bu_0 = yu_0$.
Now we have both $h_0(yu_0) \neq yu_0$ and $h_0(yu_0) = yu_0$.
  \ \ \ $\Box$

\medskip

\begin{pro} \label{general_fix_commA2} \
Let $G \subset {\mathcal G}_{n,1}$ be a subgroup
(with an embedding), and let $P$, $P'$ be complementary prefix codes over
$A$, with $n = |A|$.

\smallskip

\noindent {\bf (1)} \ For all $g \in {\rm pFix}_G(PA^*)$ and all
$h \in {\rm pFix}_G(P'A^*)$: \ $gh = hg$.

\smallskip

\noindent {\bf (2)} \ Suppose ${\rm pFix}_G(P'A^*)$ is separating on
$PA^*$.  Then we have for every $g \in G$:

\smallskip

  \ \ \ \ \  if \ $gh = hg$ \ for all \ $h \in {\rm pFix}_G(P'A^*)$, then
   \ $g \in {\rm pFix}_G(PA^*)$.
\end{pro}

\noindent {\bf Proof.} \ (1) is
straightforward. To prove (2), suppose $g \in G$
commutes with all $h \in {\rm pFix}_G(P'A^*)$. Then $g^{-1}$ also commutes
with all $h \in {\rm pFix}_G(P'A^*)$.

\smallskip

\noindent {\sc Claim:} \ $g$ stabilizes $PA^*$ and $P'A^*$, where defined.

Proof of the Claim: Assume, by contradiction, that $g(x') = y$ for some
$x' \in P'A^*$, $y \in PA^*$. Then we have for all
$h \in {\rm pFix}_G(P'A^*)$: \ $hg(x') = gh(x') = g(x')$. So, all of
${\rm pFix}_G( P'A^*)$ fixes $y \in PA^*$. This contradicts the maximality
of ${\rm pFix}_G( P'A^*)$, and hence it contradicts the separation property,
by Lemma \ref{separ_maxA2}. Therefore, $g$ maps $P'A^*$ into $P'A^*$.
Similarly, $g^{-1}$ maps $P'A^*$ into $P'A^*$.

If we had $g(x) = y'$ for some $x \in PA^*$, $y' \in P'A^*$, then
$g^{-1}(y') = x$, contradicting the fact that $g^{-1}$ maps $P'A^*$ into
$P'A^*$. Thus, $g$ maps $PA^*$ into $PA^*$. This proves the Claim.

\smallskip

To prove that $g \in {\rm pFix}_G(PA^*)$ , assume by contradiction that
$g(x_1) = y_1 \neq x_1$, for some $x_1 \in PA^*$; by the Claim, $y_1 \in PA^*$.
By Lemma \ref{non_id}, there exist therefore $x, y \in PA^*$ such that
$g(x) = y$ and $x$ and $y$ are prefix-incomparable.
Now we have for all $h \in {\rm pFix}_G( P'A^*)$: \ $gh(x) = h(y)$.
Hence for all $u \in A^*$, \ $gh(xu) = h(yu)$.

But by the separating assumption, there exists $h_0 \in {\rm pFix}_G(P'A^*)$
and there exists $u_0 \in A^*$ such that
$h_0(yu_0) \neq yu_0$ and $h_0(xu_0) = xu_0$; the latter, together with
$gh_0(xu_0) = h_0(yu_0)$, proved above, implies $(g(xu_0) =)$
$yu_0 = h_0(yu_0)$. Now we have both $h_0(yu_0) \neq yu_0$ and 
$yu_0 = h_0(yu_0)$, a contradiction.
 \ \ \ $\Box$

\begin{pro} \label{G_n1_sepA2} \
Let $P, P' \subset A^*$ be finite non-empty complementary prefix codes with
$|A| = n \geq 2$. Then ${\rm pFix}_G(P'A^*)$ is separating on $P A^*$
for the following groups taken for $G$:

\smallskip

{\bf (1)} \hspace{1in} $G \ = \ G_{n,1}$,

\smallskip

{\bf (2)} \hspace{1in} $G \ = \ G_{n,1}^{\rm mod \, 3}$.
\end{pro}
{\bf Proof.} \
{\bf (1)}  Let $x, y \in PA^*$ be two prefix-incomparable words. Let
$a \neq b \in A$ (any two different letters); note that this makes the words
$x, ya, yb$ prefix-incomparable two-by-two (for $x$ and $ya$, use Lemma
\ref{pref_comp}, and similarly for $x$ and $yb$). Now use Lemma
\ref{extending_to_max_pref_code} to construct a maximal prefix code \   
$Q \cup \{x, ya, yb\} \cup P'$, with $Q \subset PA^*$.
Define $h_0 \in G_{n,1}$ by \

\smallskip

$h_0(ya) = yb$, \  $h_0(yb) = ya$, \ $h_0(x) = x$, and 

$h$ is the identity on $Q \cup P'$. 

\smallskip

\noindent 
So, $Q \cup \{x, ya, yb\} \cup P'$ is the domain code and image code of $h_0$.
Then $h_0 \in {\rm pFix}_{G_{n,1}}(P'A^*)$, \ $h_0(ya) \neq ya$, and
$h(x_0a) = xa$ (since $h_0(x) = x$). So here, $a$ plays the role of $u_0$
in the separation property.

\medskip

\noindent {\bf (2)} \
Note that $h_0$ preserves lengths, so $h_0 \in G_{n,1}^{\rm mod \, 3}$, which
proves that ${\rm pFix}_{G_{n,1}^{\rm mod \, 3}}(P'A^*)$ is separating too.
  \ \ \ $\Box$

\begin{cor} {\bf (Commutation test).} \\      
Let $A$ be an alphabet with $|A| = n \geq 2$. 
Let $G = G_{n,1}$ or $G = G_{n,1}^{\rm mod \, 3}$, and let $P, P'$ be
complementary prefix codes over $A$. Then for any $g \in G$ 
we have:

\smallskip

$g \in {\rm pFix}_G(P A^*)$ \ \ iff \ \
$gh = hg$ \ for all $h \in {\rm pFix}_G(P' A^*)$.
\end{cor}

%%%%%%%%%%%%%%%%%%%%%%%%%%%%%%%%%%%%%%%%%%%%%%%%%%%%%%%%%%%%%%%

\newpage

\subsection{Commutation test, finite presentation, and word problem of 
$G_{3,1}^{\rm mod \, 3}(0,1)$}  % A3

We prove that the commutation test works for $G_{3,1}^{\rm mod \, 3}(0,1)$,
thus reducing the circuit equivalence problem to the word problem of
$G_{3,1}^{\rm mod \, 3}(0,1)$ (over an infinite generating set).

Then we show that $G_{3,1}^{\rm mod \, 3}(0,1)$ is finitely presented.

Finally, we embed $G_{3,1}^{\rm mod \, 3}(0,1)$ into a finitely presented 
Thompson group $H(0,1)$, thus showing that $H(0,1)$ is a finitely presented
group with coNP-hard word problem. 

So, $G_{3,1}^{\rm mod \, 3}(0,1)$ and the corresponding group $H(0,1)$
have similar properties as $G_{3,1}^{\rm mod \, 3}(0,1;\#)$ and 
$H(0,1;\#)$. The proofs are similar too, but a little more complicated in
the case of $G_{3,1}^{\rm mod \, 3}(0,1)$. The other sections of this paper
do not depend on this subsection.

\begin{defn} \label{def_max_sep_endmA3} \
Let \ $G \subset {\mathcal G}_{3,1}^{\rm mod \, 3}(0,1)$.
Let $P, P'$ be complementary prefix codes over $\{0,1,\#\}$,  with
$P \cap \{0,1\}^* \neq \emptyset$, $P' \cap \{0,1\}^* \neq \emptyset$.
The fixator ${\rm pFix}_G(P'\{0,1,\#\}^*)$ is {\em separating} on 
$P\{0,1,\#\}^*$ iff the following hold for any ordered pair $(x,y)$ of
prefix-incomparable words \    
$x, y \in \{0,1\}^* \ \cup \ \{0,1\}^*\# \{0,1,\#\}^*$:

\noindent $\bullet$ \  If \ $x, y \in (P \cap \{0,1\}^*) \, \{0,1\}^*$, 
then there exists $h \in {\rm pFix}_G(P' \{0,1,\#\}^*)$ and there exists 
$u \in \{0,1\}^*$ such that

\smallskip

 \ \ \ \ \    $h(xu) = xu$ \ and \ $h(yu) \neq yu$.

\smallskip

\noindent $\bullet$ \ If \ $x, y \not\in \{0,1\}^*$, and \  
$x, y \in P \{0,1,\#\}^*$ \ there exists
$h \in {\rm pFix}_G(P' \{0,1,\#\}^*)$ such that

\smallskip

 \ \ \ \ \   $h(x) = x$ \ and \ $h(y) \neq y$.
\end{defn}
(Note that we don't have any requirements in the case where $x \in \{0,1\}^*$
and $y \notin \{0,1\}^*$, or the case where $x \notin \{0,1\}^*$ and
$y \in \{0,1\}^*$.)

\begin{lem} \label{sep_implies_max_endmA3} \
If ${\rm pFix}_G(P'\{0,1,\#\}^*)$ is separating on $P \{0,1,\#\}^*$ then it 
is a maximal fixator.
\end{lem}
{\bf Proof.} \ Recall the definition of a maximal fixator from Section 5. 
Suppose by contradiction that there exists $x_0 \in P\{0,1,\#\}^*$ such that 
$h(x_0) = x_0$ for all $h \in {\rm pFix}_G(P'\{0,1,\#\}^*)$. By Lemma 
\ref{endmarker_code2}, the prefix code $P$ is of the form \ 
$P \ = \ P_1 \cup \bigcup_{v \in P_2} v \, \# \, P(v)$, where 
$P_1 = P \cap \{0,1\}^*$.

\smallskip

\noindent Case 1: \ $x_0 \in P_1 \{0,1\}^*$.

Choose $x = x_00$ and $y = x_01$. Then $x$ and $y$ are prefix incomparable,
hence by the separation property of the fixator, there exists
$h_0 \in {\rm pFix}_G(P'\{0,1,\#\}^*)$ and $u_0 \in \{0,1\}^*$ with

\smallskip

 \ \ \ \ \   $h_0(xu_0) = xu_0$, \ $h_0(yu_0) \neq yu_0$.

\smallskip

\noindent However, $h_0(yu_0) \neq yu_0$ contradicts the fact that
$h_0(x_0) = x_0$.

\smallskip

\noindent Case 2: \ $x_0 \in P_1 \{0,1\}^* \# \{0,1,\#\}^*$, or \
$x_0 \in \bigcup_{v in P_2} v \, \# \, P(v)$ \
  with $|P_2| \geq 2$.

Let $x_0 = v_0 \# s$.
Let $w_0 \in P_2$ with $w_0 \neq v_0$, and choose $x = w_0\# t$
(for some $t \in P(w_0)$) and $y = v_0\# s$.
Then $x$ and $y$ are prefix incomparable, and both are in
$\{0,1\}^* \# \{0,1,\#\}^*$; so there exists
$h_0 \in {\rm pFix}_G(P'\{0,1,\#\}^*)$ with

\smallskip

 \ \ \ \ \   $h_0(x) = x$, \ $h_0(y) \neq y$.

\smallskip

\noindent However, $h_0(y) \neq y$ contradicts the fact that
$h_0(x_0) = x_0$.

\smallskip

\noindent Case 3: \ $x_0 \in \bigcup_{v in P_2} v \, \# \, P(v)$ \
and $|P_2| = 1$. (Obviously the case $|P_2| = 0$ cannot occur when 
$x_0 \in P_2 \#$.)

Then $P_2 = \{v_0 \}$, so we have \ $x = v_0\# s  \in v_0\#P(v_0)$, for 
some $s \in P(v_0)$. Let $z_0 \in P_1$ (recall that we assume 
$P_1 \neq \emptyset$).  Let $x = z_0\#$ and $y = x_0 = v_0\#s$. 
Since $z_0 \neq v_0$, \
$x$ and $y$ are prefix incomparable, and both are in
$\{0,1\}^* \# \{0,1,\#\}^*$; so there exists
$h_0 \in {\rm pFix}_G(P'\{0,1,\#\}^*)$ with

\smallskip

 \ \ \ \ \   $h_0(x) = x$, \ $h_0(y) \neq y$.

\smallskip

\noindent Again, $h_0(y) \neq y$ contradicts the fact that
$h_0(x_0) = x_0$.
  \ \ \ $\Box$

\begin{pro} \label{max_sep_end_commA3} \
Let $P, P'$ be complementary prefix codes over $\{0,1,\#\}$ with
$P \cap \{0,1\}^*$ and $P'\cap \{0,1\}^*$ non-empty. Suppose that \    
$G \subset {\mathcal G}_{3,1}^{\rm mod \, 3}(0,1)$ is a group, and that
${\rm pFix}_G(P'\{0,1,\#\}^*)$ is separating on $P \{0,1,\#\}^*$.
Then for all $g \in G$ we have:

\smallskip

 \ \ \  If $g$ commutes with all elements of 
${\rm pFix}_G(P'\{0,1,\#\}^*)$ then $g \in {\rm pFix}_G(P\{0,1,\#\}^*)$.
\end{pro}
{\bf Proof.} \  \ Let $g \in G$ and assume $g$ commutes with all elements 
of ${\rm pFix}_G(P'\{0,1,\#\}^*)$. We want to show that 
$g \in {\rm pFix}_G(P\{0,1,\#\}^*)$. We first prove:

\smallskip

\noindent {\sc Claim:} \ $g$ stabilizes $P' \{0,1,\#\}^*$ and
$P \{0,1,\#\}^*$.

\smallskip

\noindent Proof of the Claim: \ Assume by contradiction that $g(x') = y$ for
some $x' \in P' \{0,1,\#\}^*$ and $y \in P \{0,1,\#\}^*$. Since $g$ commutes
with all elements of the fixator we have for all
$h \in {\rm pFix}_G(P'\{0,1,\#\}^*)$:
 \ $gh(x') = hg(x') = g(x') = y$, i.e., $h(y) = y$. This contradicts the
maximality of the fixator ${\rm pFix}_G(P'\{0,1,\#\}^*)$; so $g$ maps
$P'\{0,1,\#\}^*$ into itself.

Similarly, $g^{-1}$ maps $P'\{0,1,\#\}^*$ into itself. From this it follows
(as in the proof of Proposition \ref{general_fix_commA2}) that
$g$ also maps $P \{0,1,\#\}^*$ into itself, and similarly for $g^{-1}$.
This proves the Claim.

\medskip

\noindent Assume now by contradiction that $g$ does not fix some element
$x_1 \in P \{0,1,\#\}^*$: \ $g(x_1) = y_1 \neq x_1$.
By the Claim, $y_1 \in P \{0,1,\#\}^*$.

By Lemma \ref{non_id} there exist $x, y \in P \{0,1,\#\}^*$ such that $x$
and $y$ are prefix incomparable and $g(x) = y$.
And since $g$ commutes with the fixator,
we have for all $h \in {\rm pFix}_G(P'\{0,1,\#\}^*)$: \ $gh(x) = hg(x) = h(y)$.

On the other hand, the separation property of the fixator implies that there
exists $h_0 \in {\rm pFix}_G(P'\{0,1,\#\}^*)$ and $u_0 \in \{0,1,\#\}^*$
such that \ $h_0(yu_0) \neq yu_0$ and $h_0(xu_0) = xu_0$.

The equality $gh(x) = h(y)$ implies $gh_0(xu_0) = h(yu_0)$; this, together
with $h_0(xu_0) = xu_0$, implies $yu_0 = gh_0(xu_0) = h(yu_0)$. But this
contradicts $h_0(yu_0) \neq yu_0$.
       \ \ \ $\Box$

\medskip

In the next two Lemmas we will check that the group
$G = G_{3,1}^{\rm mod \, 3}(0,1)$ 
satisfies the conditions of Proposition \ref{max_sep_end_commA3}, i.e.,
that the fixator ${\rm pFix}_G(P'\{0,1,\#\}^*)$ is separating.

\begin{pro} \label{G_31_01_sepA3} \
Let $P, P'$ be complementary prefix codes over $\{0,1,\#\}$ with
$P \cap \{0,1\}^*$ and $P'\cap \{0,1\}^*$ non-empty. Let
$G = G_{3,1}^{\rm mod \, 3}(0,1)$.
Then the fixator \ ${\rm pFix}_G(P'\{0,1,\#\}^*)$ \ is separating on
$P \{0,1,\#\}^*$.
\end{pro}
{\bf Proof.} \ Let $x, y \in P_1 \{0,1\}^*$ or
$x, y \in P_1 \{0,1\}^*\# \cup \bigcup_{v \in P_2} v \, \# \, P(v)$,
and assume $x$ and $y$ are prefix incomparable.
We want to find $h_0 \in {\rm pFix}_G(P'\{0,1,\#\}^*)$ and
$u_0 \in \{0,1\}^*$ such that $h_0(xu_0) = xu_0$ and
$h_0(yu_0) \neq yu_0$, or $h_0(xu_0) \neq xu_0$ and $h_0(yu_0) = yu_0$;
 \ if $x, y \notin \{0,1\}^*$ then $u_0$ is empty.

\medskip

\noindent {\sc Case 1:} \ $x, y \in P_1 \{0,1\}^*$.

\smallskip

In this case we can apply the same proof as for Proposition \ref{G_n1_sepA2},
with alphabet $A = \{0,1\}$.

\medskip

\noindent {\sc Case 2:} \
$x, y \in P_1 \{0,1\}^*\# \{0,1,\#\}^* \ \cup $
 $ \ \bigcup_{v \in P_2} v \, \# \, P(v) \ \{0,1,\#\}^*$.

\smallskip

Let $x = x_0\#s_0$ and $y = y_0\#t_0$

\smallskip

\noindent {\sc Case 2.1:} \ $y_0 \in P_1 \{0,1\}^*$.

\smallskip

\noindent $\bullet$ \ Assume $x_0 \in P_2$.

Now, $x_0 \neq y_0$, since $x_0 \in P_2$ and $y_0 \in P_1 \{0,1\}^*$;
since $P_2$ is closed under prefix (by Lemma \ref{endmarker_code2}), 
$y_0$ is not a prefix of $x_0$.
By Lemma \ref{extending_to_max_pref_code} over the
alphabet $A = \{0,1\}$, there is a finite prefix code
$Q_1 \subset P_1 \{0,1\}^*$ such that
$Q_1 \cup \{y_0 00 \}$ and $P_1'$ and complementary prefix codes
(over $\{0,1\}$). Therefore the following will be a finite
maximal prefix code over $\{0,1,\#\}$:

\smallskip

 \ \ \ \ \  $C \ = \ Q_1 \cup \{y_0 00\} \cup P'_1$
    $\cup \ \bigcup_{v \in Q_2} v \, \# \ \Pi(v)$,

\smallskip

\noindent where
 \ $Q_2 \ = \ \ >_{{\rm pref}}\!\!(Q_1 \cup \{y_0 00\} \cup P'_1)$;
moreover, $\Pi(v) = P(v)$ if $v \in P_2$, \ $\Pi(v) = P'(v)$ if 
$v \in P_2'$, and $\Pi(v)$ consists of just the empty word $\varepsilon$
otherwise.

\smallskip

\noindent Now we define $h_0$, with domain code and image code $C$, by

\smallskip

 \ \ \ \ \  $h_0(y_0\# t_0) = y_0 0 \#$,
          \ $h_0(y_0 0 \#) = y_0 \# t_0$,  \ and 

 \ \ \ \ \ $h_0$ is the identity everywhere else on $C$.

\smallskip

\noindent Thus, $h_0(y) \neq y$.
Moreover, $h_0 \in {\rm pFix}_G(P'\{0,1,\#\}^*)$, because $y_0 0 \# $
$\notin \ \bigcup_{v in P_2'} v \# P'(v)$; indeed,
$y_0 0 \# \in P_1 \{0,1\}^* \#$ $\subset P \{0,1,\# \}^*$.

And $h_0$ preserves the length of strings in $\{0,1\}^*$ (since $h_0$ is
the identity on $\{0,1\}^*$ wherever $h_0$ is defined).

\smallskip

We also claim that $h_0(x) = x$.
Indeed, $x_0$ belongs to $P_2$, which is contained in
$>_{{\rm pref}}\!\!(P_1 \cup P'_1)$; moreover,
$>_{{\rm pref}}\!\!(P_1) \ \subset \ >_{{\rm pref}}\!\!(Q_1)$, by the 3rd 
point of Lemma \ref{extending_to_max_pref_code}; and $s_0 \in P'(x_0)$.
Therefore, $x_0 \# s_0$ belongs to $C$.  On the other hand, $x_0$ is
different from $y_0$ and $y_0 0$.

\smallskip

\noindent $\bullet$ \ Assume $x_0 \in P_1 \{0,1\}^*$.

Then, by Lemma \ref{making_pref_codes}, there are $\ell_1, \ell_2$
$\in \{0,1\}$ such that $x_0 \ell_1$ and $y_0 \ell_2$ are prefix
incomparable. By applying Lemma \ref{extending_to_max_pref_code} over
the alphabet $A = \{0,1\}$ we obtain a finite prefix code
$Q_1 \subset P_1 \{0,1\}^*$ such that
$Q_1 \cup \{x_0 \ell_1, y_0 \ell_2 0\} $ and $P_1'$ and complementary
prefix codes (over $\{0,1\}$). Therefore the following set \
$C \subset \{0,1\}^* \cup \{0,1\}^*\#$ \ will be a finite
maximal prefix code over $\{0,1,\#\}$:

\smallskip

 \ \ \ \ \  $C \ = \ Q_1 \cup \{x_0 \ell_1, y_0 \ell_2 0\} \cup P'_1$
  $\cup \ \ \bigcup_{v \in Q_2} v \, \# \ \Pi(v)$,

\smallskip

\noindent where
 \ $Q_2 \ = \ $
$>_{{\rm pref}}\!\!(Q_1 \cup \{x_0 \ell_1, y_0 \ell_2 0\} \cup P'_1)$,
and where \
$\Pi(v) = P(v)$ if $v \in P_2$, \ $\Pi(v) = P'(v)$ if $v \in P_2'$, and
$\Pi(v)$ consists of just the empty word otherwise.

\smallskip

\noindent Now we define $h_0$, with domain code and image code $C$, by

\smallskip

 \ \ \ \ \ $h_0(y_0\# t_0) = y_0 \ell_2 \#$,
         \ $h_0(y_0 \ell_2 \#) = y_0 \# t_0$, \  and 

 \ \ \ \ \ $h_0$ is the identity everywhere else on $C$.

\smallskip

\noindent Thus, $h_0(y) \neq y$.  Moreover,
$h_0 \in {\rm pFix}_G(P'\{0,1,\#\}^*)$, because $y_0 \ell_2 \#$
$\notin \ \bigcup_{v \in P_2'} v \# P'(v)$;  indeed, $y_0 \ell_2 \#$
$\in P_1 \{0,1\}^*\#$ $\subset P \{0,1,\# \}^*$.

And $h_0$ preserves the
length of strings in $\{0,1\}^*$ (since $h_0$ is the identity on
$\{0,1\}^*$ wherever it is defined).  Also, $h_0(x) = x$, since
$x_0 \# s_0$  belongs to $C$ (since $x_0$ is a strict prefix of
$x_0 \ell_1$), and since $x_0$ is different from $y_0$ and
$y_0\ell_2$.

\medskip

\noindent {\sc Case 2.2:} \ $y_0 \in P_2$.

\smallskip

Since $P_1 \neq \emptyset$, there exists $w_0 \in P_1$; hence
$y_0$ is different from $w_0$, $w_0 0$, and $w_0 00$.    Also,
$x_0$ is different from $w_0 0$ or from $w_0 00$ (or from both).
Let $z_0 0$ be one of $w_0 0$ or $w_0 00$, so that $z_0 0\neq x_0$.
We still have $z_0 0 \neq y_0$ and $z_0 0 \in P_1 \{0,1\}^*$.

\smallskip

\noindent $\bullet$ \ Assume $x_0 \in P_2$.

By Lemma \ref{extending_to_max_pref_code} over the
alphabet $A = \{0,1\}$, there is a finite prefix code
$Q_1 \subset P_1 \{0,1\}^*$ such that
$Q_1 \cup \{z_00\}$ and $P_1'$ and complementary prefix codes
(over $\{0,1\}$). Therefore the following set \
$C \subset \{0,1\}^* \cup \{0,1\}^*\#$  \ will be a finite maximal
prefix code over $\{0,1,\#\}$:

\smallskip

 \ \ \ \ \  $C \ = \ Q_1 \cup \{z_0 0\} \cup P'_1 \cup \ $
  $\bigcup_{v \in Q_2} v \, \# \ \Pi(v)$,

\smallskip

\noindent where
 \ $Q_2 \ = \ >_{{\rm pref}}\!\!(Q_1 \cup \{z_0 0\} \cup P'_1)$,
and where \
$\Pi(v) = P(v)$ if $v \in P_2$, \ $\Pi(v) = P'(v)$ if $v \in P_2'$, and
$\Pi(v)$ consists of just the empty word otherwise.

\smallskip

\noindent Now we define $h_0$, with domain code and image code $C$,
by

\smallskip

 \ \ \ \ \  $h_0(y_0\# t_0) = z_0 \#$,
          \ $h_0(z_0\#) = y_0\# t_0$, \ and 

 \ \ \ \ \  $h_0$ is the identity everywhere else on $C$.

\smallskip

\noindent Thus, $h_0(y) \neq y$ and $h_0(x) = x$. Note that
$h_0(x_0\# s_0)$ and $h_0(y_0\# t_0)$ are defined since
$x_0, y_0 \in P_2 \ \subset \ $
$>_{{\rm pref}}\!\!(P_1 \cup P'_1)$; moreover,
$>_{{\rm pref}}\!\!(P_1) \ \subset \ >_{{\rm pref}}\!\!(Q_1)$, by the 3rd 
point of Lemma \ref{extending_to_max_pref_code}. Therefore, $x_0\# s_0$ and 
$y_0\# t_0$ belong to $C$.

Also, $h_0 \in {\rm pFix}_G(P'\{0,1,\#\}^*)$, because $z_0 \# \notin$
$ \ \bigcup_{v \in P_2'} v \# P'(v)$; indeed, $z_0 \# \in$
$P_1 \{0,1\}^* \#$ $\subset P \{0,1,\#\}^*$.

Also, $h_0$ preserves the length of strings in $\{0,1\}^*$ since $h_0$
is the identity on $\{0,1\}^*$ wherever it is defined.

\smallskip

\noindent $\bullet$ \ Assume $x_0 \in P_1 \{0,1\}^*$.

Then, by Lemma \ref{making_pref_codes}, there are $\ell_1, \ell_2$
$\in \{0,1\}$ such that $x_0 \ell_1$ and $z_0 \ell_2$ are prefix
incomparable. By applying Lemma \ref{extending_to_max_pref_code} over
the alphabet $A = \{0,1\}$ we obtain a finite prefix code
$Q_1 \subset P_1 \{0,1\}^*$ such that
$Q_1 \cup \{x_0 \ell_1, z_0 \ell_2 \} $ and $P_1'$ and complementary
prefix codes (over $\{0,1\}$). Therefore the following set \
$C \subset \{0,1\}^* \cup \{0,1\}^*\#$ \ will be a finite
maximal prefix code over $\{0,1,\#\}$:

\smallskip

 \ \ \ \ \  $C \ = \ Q_1 \cup \{x_0 \ell_1, z_0 \ell_2 \} \cup P'_1$
  $\cup \ \bigcup_{v \in Q_2} v \, \# \ \Pi(v)$,

\smallskip

\noindent where
 \ $Q_2 \ = \ $
  $>_{{\rm pref}}\!\!(Q_1 \cup \{x_0 \ell_1, z_0 \ell_2 \} \cup P'_1)$,
and where \
$\Pi(v) = P(v)$ if $v \in P_2$, \ $\Pi(v) = P'(v)$ if $v \in P_2'$, and
$\Pi(v)$ consists of just the empty word otherwise.

\smallskip

\noindent Now we define $h_0$, with domain code and image code $C$, by

\smallskip

 \ \ \ \ \ $h_0(y_0\# t_0) = z_0 \#$,
         \ $h_0(z_0\#) = y_0\# t_0$, \ and 

 \ \ \ \ \ $h_0$ is the identity everywhere else on $C$.

\smallskip

\noindent Thus, $h_0(y) \neq y$, and $y \in C$ (since
$y_0 \in P_2 \ \subset \ >_{{\rm pref}}\!\!(P_1 \cup P'_1)$ $ \ \subset \ $
$>_{{\rm pref}}\!\!(Q_1 \cup P'_1)$ ).

Moreover, $h_0 \in {\rm pFix}_G(P'\{0,1,\#\}^*)$, because $z_0 \#$
$\notin P' \{0,1,\#\}^*$ for the same reasons as in the previous subcase.

And $h_0$ preserves the
length of strings in $\{0,1\}^*$ (since $h_0$ is the identity on
$\{0,1\}^*$ wherever it is defined).  Also, $h_0(x) = x$, since
$x_0 \# s_0$  belongs to $C$ (since $x_0$ is a strict prefix of
$x_0 \ell_1$), and since $x \neq y$ and $x_0 \neq z_0$.
           \ \ \ $\Box$

\bigskip

\noindent Propositions \ref{max_sep_end_commA3} and \ref{G_31_01_sepA3} 
immediately imply:
\begin{cor} {\bf (Commutation test for $G_{3,1}^{\rm mod \, 3}(0,1)$).} \\      
Let \ $G = G_{3,1}^{\rm mod \, 3}(0,1)$. For any $g \in G$ we have:

\smallskip

$g \in {\rm Fix}_G(0 \, \{0,1,\#\}^*)$ \ \ iff \ \ $gh = hg$
 \ for all $h \in {\rm pFix}_G( \{1,\# \}\{0,1,\#\}^*)$.
\end{cor}

\smallskip

Now, by the same reasoning as in Section 5, the above Corollary reduces 
the circuit equivalence problem to the word problem of 
$G_{3,1}^{\rm mod \, 3}(0,1)$; the reduction is an unbounded conjunctive
reduction. The next Lemma implies that
${\rm pFix}_G( \{1,\# \}\{0,1,\#\}^*)$ is isomorphic to 
$G = G_{3,1}^{\rm mod \, 3}(0,1)$.  This and the fact (proved later in this
subsection) that $G_{3,1}^{\rm mod \, 3}(0,1)$ is finitely generated 
implies that 
only the finitely many generators of ${\rm pFix}_G( \{1,\# \}\{0,1,\#\}^*)$ 
need to be used in the role of ``$h$'' in the above Corollary. 
This then yields:

\begin{cor} \
 The circuit equivalence problem reduces to the word problem of
$G_{3,1}^{\rm mod \, 3}(0,1)$ (over an infinite generating set),
by a polynomial-time $k$-bounded conjunctive reduction. Here,
$k$ is the minimum number of generators of $G_{3,1}^{\rm mod \, 3}(0,1)$.
\end{cor}

\begin{lem} \label{Fix_iso_V_endmA3} \   
For $G = G_{3,1}^{\rm mod \, 3}(0,1)$,
the subgroup ${\rm pFix}_G(\{1,\# \}\{0,1,\#\}^*)$ is isomorphic to $G$.
\end{lem}
{\bf Proof.} \ An element
$\varphi \in G = G_{3,1}^{\rm mod \, 3}(0,1)$
$ =  {\rm pStab}_{G_{3,1}^{\rm mod \, 3}}(\{0,1\}^*)$ belongs to
${\rm Fix}_G(\{1,\# \}\{0,1,\#\}^*)$ iff $\varphi$ has
a table of the form

\[ \varphi \ = \   \left[ \begin{array}{cc ccc ccc}
1 & \# & 0x_1 & \ldots & 0x_n & 0x_1'\#s_1 & \ldots & 0x_m'\#s_m  \\
1 & \# & 0y_1 & \ldots & 0y_n & 0y_1'\#t_1 & \ldots & 0y_m'\#t_m
\end{array}        \right]
\]
where $x_i,y_i, x'_j, y'_j$ range over $\{0,1\}^*$, $|x_i| \equiv |y_i|$
mod 3 ($i = 1, \ldots, n$), and $s_j, t_j \in \{0,1, \#\}^*$.
The isomorphism to $G$, as above, just maps this table to
\[ \psi \ = \  \left[ \begin{array}{ccc ccc}
x_1 & \ldots & x_n & x_1'\#s_1 & \ldots & x_m'\#s_m  \\
y_1 & \ldots & y_n & y_1'\#t_1 & \ldots & y_m'\#t_m
\end{array}        \right]
\]
It is straightforward to see that this is an isomorphism.
 \ \ \ $\Box$

\bigskip

We will prove next that the group $G_{3,1}^{\rm mod \, 3}(0,1)$ is finitely
presented. As in Section 6, we will follow Higman's method and, accordingly, 
we will have
to prove appropriate facts about maximal prefix codes over $\{0,1,\#\}$.
We will use the following notation (as before):
For any maximal prefix codes $P, Q \subset \{0,1,\# \}^*$, we let
$Q_1 = Q \cap \{0,1\}^*$ and $P_1 = P \cap \{0,1\}^*$; these are maximal
prefix codes over $\{0,1\}$.

\begin{lem} \label{G(0,1)codesA3}  \
Let $Q$ be a finite maximal prefix code over $\{0,1,\# \}$. Suppose the
inner tree $T_{\rm in}(Q)$ has a leaf $\ell \in \{0,1,\# \}^*$.
Assume that there exist elements $q_1, q_2 \in Q$ such that,
$q_1, q_2$ are not children of $\ell$ (in the prefix tree of $Q$),
and such that either both $q_1, q_2 \notin \{0,1\}^*$, or both
$q_1, q_2 \in \{0,1\}^*$; in the latter case we also assume that
$|q_1| \equiv |q_2|$ {\rm mod 3}.

Then there exists a finite maximal prefix code $P$ over $\{0,1,\# \}$
such that $|P| = |Q|$, $P_1$ has the same {\rm mod 3} cardinality as $Q_1$,
and the inner tree $T_{\rm in}(P)$ has two leaves
$\ell_1$, $\ell_2 \in \{0,1,\# \}^*$ such that:

\smallskip

\noindent (1) \ If both $q_1, q_2 \in \{0,1\}^*$ then
 $\ell_2 \in \{0,1\}^*$, with
 $|\ell_2| + 1 \equiv  |q_1| \equiv |q_2|$ {\rm mod 3}.
 Moreover, if $\ell \notin \{0,1\}^*$ then $\ell_1 \notin \{0,1\}^*$;
 if $\ell \in \{0,1\}^*$ then $\ell_1 \in \{0,1\}^*$,
 and $|\ell_1| \equiv  |\ell|$ {\rm mod 3}.

\smallskip

\noindent (2) \ If both $q_1, q_2 \notin \{0,1\}^*$, and if $q_1$
 or $q_2 \notin \{0,1\}^*\#$, then $\ell_2 \notin \{0,1\}^*$.
 Moreover, if $\ell \notin \{0,1\}^*$ then $\ell_1 \notin \{0,1\}^*$;
 if $\ell \in \{0,1\}^*$ then $\ell_1 \in \{0,1\}^*$,
 and $|\ell_1| \equiv  |\ell|$ {\rm mod 3}.

\smallskip

\noindent (3) \ If both $q_1, q_2 \in \{0,1\}^*\#$, and if
 $\ell \notin \{0,1\}^*$, then $\ell_2 \notin \{0,1\}^*$.
\end{lem}
Note that the case where $q_1, q_2 \in \{0,1\}^*\#$ and $\ell \in \{0,1\}^*$
is not considered in the Lemma (and will not be needed).

\smallskip

\noindent
{\bf Proof.} \ Obviously, $T_{\rm in}(Q)$ has at least one leaf.

\smallskip

\noindent (1) \ If $q_1, q_2 \in \{0,1\}^*$ we do the following transformation
on $Q$, where $z_1$ and $z_2$ are the parent vertices of $q_1$, respectively
$q_2$, and where ST1,ST2, ST3, ST4 are subtrees (below $z_1$ or $z_2$) of the
prefix tree of $Q$. If $z_1$ and $z_2$ are on a common path from the root we
let $z_1$ be the deeper one of the two.
\begin{verbatim}
        |         \       becomes          |         \
        z1         z2                      z1         z2
      / | \#     / |  \#                 / | \#     / |  \#
       ST1 ST2    ST3 ST4                        ST1 ST3 ST4
                                                          |
                                                         ST2


            /                                             /
           z2                                            z2
         / | \#                                        / | \#
       .     ST4                                      . ST1 ST4
      .                   becomes                    .       |
     .                                              .       ST2
    /                                              /
   z1                                             z1
 / | \#                                         / | \#
  ST1 ST2
\end{verbatim}

Let $P$ be the prefix code described by the transformed tree.
Then $|Q| = |P|$  (since the number of vertices has not changed), and
$P_1$ has the same mod 3 cardinality as $Q_1$ (since the subtree ST1
is moved from $z_1$ to $z_2$ and $z_1$, $z_2$ have equivalent depths
mod 3). The subtree ST2 was under $\#$ and is still below a $\#$-edge.
Finally, $T_{\rm in}(P)$ has two leaves, namely $\ell$ and $z_1$
(both $\in \{0,1\}^*$), and $|z_1| + 1 \equiv |q_1|$ (actually the
two numbers are equal). The existing leaf $\ell$ is either unchanged, or
(in case it was in ST1) changed to a leaf that has an equivalent depth
modulo 3, or (in case it was in ST2, hence was $\notin \{0,1\}^*$)
changed to a leaf $\notin \{0,1\}^*$.

\medskip

\newpage

\noindent (2) \ If $q_1, q_2 \notin \{0,1\}^*$, and
$q_1 \notin \{0,1\}^*\#$ (or, similarly, if $q_2 \notin \{0,1\}^*\#$),
we do the following transformation on $Q$. As above, $z_1$ and $z_2$ are the
parent vertices of $q_1$, respectively $q_2$, and ST1, ST2, ST3, ST4, ST5
are subtrees (below $z_1$ or $z_2$) of the prefix tree of $Q$; one of ST3 or
ST5 is a single vertex (corresponding to the word $q_2$). If $z_1$ and $z_2$
are on a common path from the root we let $z_1$ be the deeper one of the two.
Since $q_1 \notin \{0,1\}^*\#$ and $q_1 \notin \{0,1\}^*$, we have
$z_1 \notin \{0,1\}^*$; i.e., $\#$ appears on the path root-to-$z_1$.

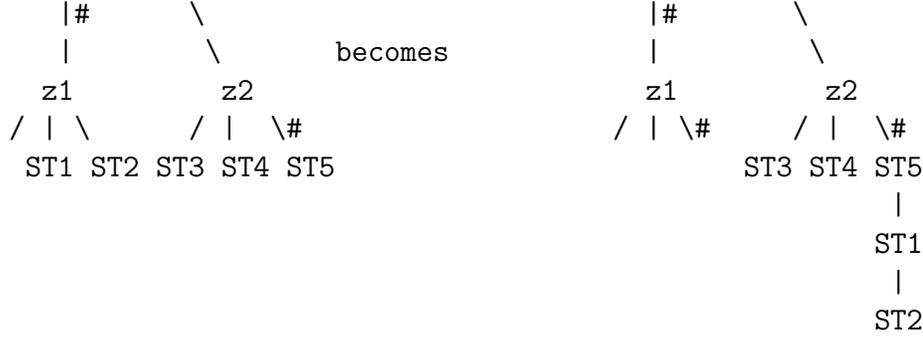
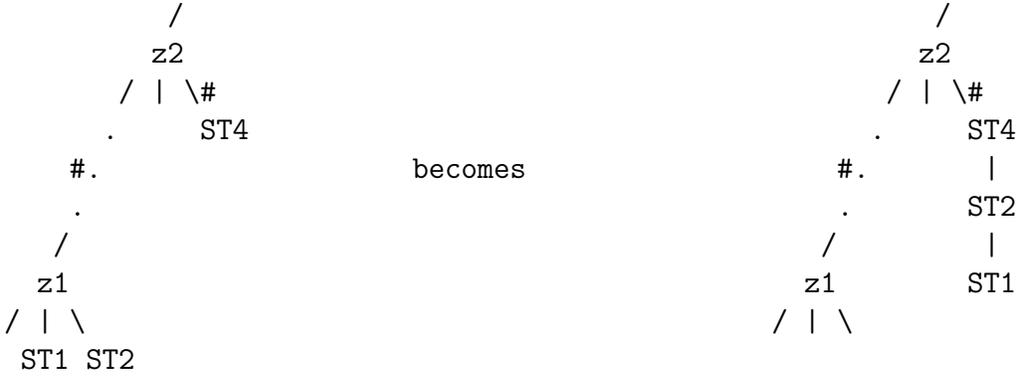
\begin{verbatim}
         |#      \                           |#       \
         |        \       becomes            |         \
        z1         z2                        z1         z2
      / | \      / |  \#                   / | \#     / |  \#
       ST1 ST2 ST3 ST4 ST5                         ST3 ST4 ST5
                                                            |
                                                           ST1
                                                            |
                                                           ST2


           /                                              /
          z2                                             z2
        / | \#                                         / | \#
       .     ST4                                      .     ST4
     #.                   becomes                   #.       |
     .                                              .       ST2
    /                                              /         |
   z1                                             z1        ST1
 / | \                                          / | \
  ST1 ST2

\end{verbatim}
Then the transformed tree describes the desired prefix code $P$.
In particular, $P_1 = Q_1$, since all the changes happen below $\#$-edges.
The existing leaf $\ell$ is either unchanged, or (in case it was in ST1 or 
ST2, in which case $\ell \notin \{0,1\}^*$) changed to a leaf that is below 
a $\#$-edge.

\medskip

\noindent (3) \ Suppose $q_1, q_2 \in \{0,1\}^*\#$ and
$\ell \notin \{0,1\}^*$.

\smallskip

\noindent $\bullet$ \ If $T_{\rm in}(Q)$ has any subtrees other than the
path root-to-$\ell$, then $T_{\rm in}(Q)$ has another leaf besides $\ell$.
In this case we have nothing to prove.

\smallskip

\noindent $\bullet$ \ If $T_{\rm in}(Q)$ consists only of the path
root-to-$\ell$, let $p$ be the parent vertex of $\ell$ in $T_{\rm in}(Q)$.
Since $q_1, q_2 \in \{0,1\}^*\#$, the path root-to-$\ell$ of $T_{\rm in}(Q)$
has some edges labeled over $\{0,1\}$, and at least one $z_1$ or $z_2$ (the
parent vertices of $q_1$ and $q_2$ in the prefix tree of $Q$) is at least
2 depth levels above $\ell$; assume $z_1$ is the deeper one.

\newpage 

\noindent We do the following transformation on $Q$.

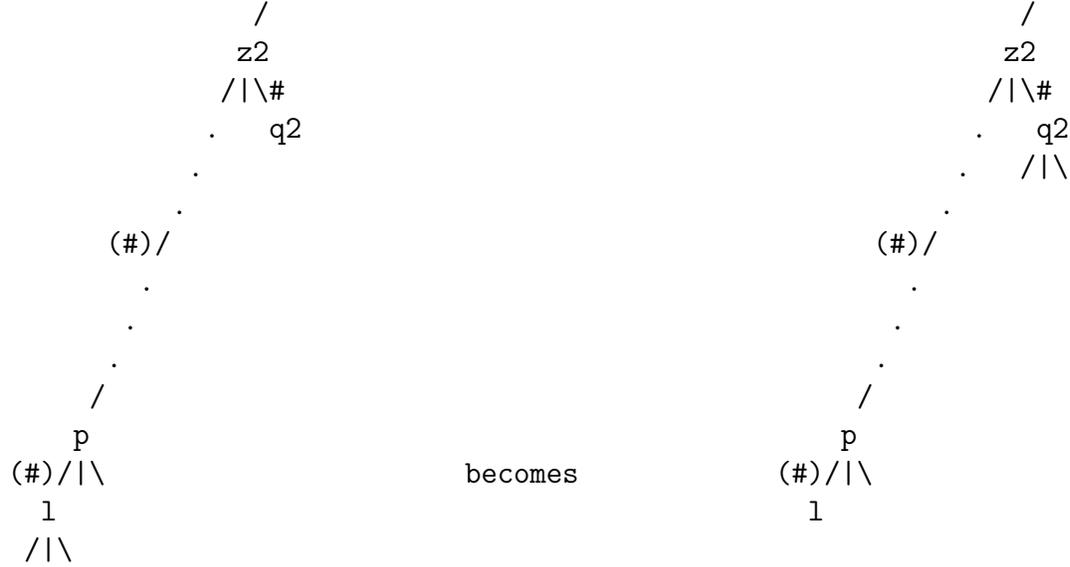
\begin{verbatim}
               /                                              /
              z2                                             z2
             /|\#                                           /|\#
            .   q2                                         .   q2
           .                                              .   /|\
          .                                              .
      (#)/                                           (#)/
        .                                              .
       .                                              .
      .                                              .
     /                                              /
    p                                              p
(#)/|\                      becomes            (#)/|\
  l                                              l
 /|\
\end{verbatim}
In other words, the three children of $\ell$ are moved to $q_2$.
Now $|P| = |Q|$ (since no additional vertices are added), and
$P_1 = Q_1$, hence $P_1$ and $Q_1$ have the same mod 3 cardinality.
Also, $T_{\rm in}(P)$  has two leaves, namely $q_2 \notin \{0,1\}^*$,
and $p$.
 \ \ \ $\Box$

\bigskip

\noindent Again, for maximal prefix codes $P, Q$ over $\{0,1,\# \}$ we will
use the notation $P_1 = P \cap \{0,1\}^*$, $Q_1 = Q \cap \{0,1\}^*$.

\begin{lem} \label{G(0,1)codes2A3} \
Let $Q$ be a finite maximal prefix code over $\{0,1,\# \}$. Suppose the
inner tree $T_{\rm in}(Q)$ has two leaves
$\ell_1, \ell_2 \in \{0,1,\#\}^*$.
Assume that there exist elements $q_1, q_2 \in Q$ such
that, $q_1, q_2$ are not children of $\ell_1$ or $\ell_2$ (in the prefix
tree of $Q$), and such that either both $q_1, q_2 \notin \{0,1\}^*$, or
both $q_1, q_2 \in \{0,1\}^*$; in the latter case we also assume that
$|q_1| \equiv |q_2|$ {\rm mod 3}.

Then there exists a finite maximal prefix code $P$ over $\{0,1,\# \}$
such that $|P| = |Q|$, $P_1$ has the same {\rm mod 3} cardinality as $Q_1$,
and the inner tree $T_{\rm in}(P)$ has three leaves
$\lambda_1$, $\lambda_2$, $\lambda_3 \in \{0,1,\# \}^*$ such that:

\smallskip

\noindent (1) \ If both $q_1, q_2 \in \{0,1\}^*$ then
 $\lambda_3 \in \{0,1\}^*$,
 with $|\lambda_3| + 1 \equiv  |q_1| \equiv |q_2|$ {\rm mod 3}.
 Moreover (for all $i = 1, 2$), both $\lambda_i, \ell_i \notin \{0,1\}^*$,
 or both $\lambda_i, \ell_i \in \{0,1\}^*$, and in the latter case
 $|\lambda_i| \equiv |\ell_i|$ {\rm mod 3}.

\smallskip

\noindent (2) \ If both $q_1, q_2 \notin \{0,1\}^*$, and if $q_1$
 or $q_2 \notin \{0,1\}^*\#$, then $\lambda_3 \notin \{0,1\}^*$.
 Moreover (for all $i = 1, 2$), both $\lambda_i, \ell_i \notin \{0,1\}^*$,
 or both $\lambda_i, \ell_i \in \{0,1\}^*$, and in the latter case
 $|\lambda_i| \equiv |\ell_i|$ {\rm mod 3}.

\smallskip

\noindent (3) \ If both $q_1, q_2 \in \{0,1\}^*\#$, and if
 $\ell_1 \notin \{0,1\}^*$, then $\lambda_2 = \ell_2$ and
 $\lambda_3 \notin \{0,1\}^*$.
 (However,  $\lambda_1$ could be $\in \{0,1\}^*$ or $\notin \{0,1\}^*$.)
\end{lem}
{\bf Proof.} \ The proof is similar to the proof of Lemma
\ref{G(0,1)codesA3}.

\smallskip

\noindent (1) \ If $q_1, q_2 \in \{0,1\}^*$ we do the following transformation
on $Q$, where $z_1$ and $z_2$ are the parent vertices of $q_1$, respectively
$q_2$, and where ST1, ST2, ST3, ST4 are subtrees (below $z_1$ or $z_2$) of the
prefix tree of $Q$. If $z_1$ and $z_2$ are on a common path from the root we
let $z_1$ be the deeper one of the two.
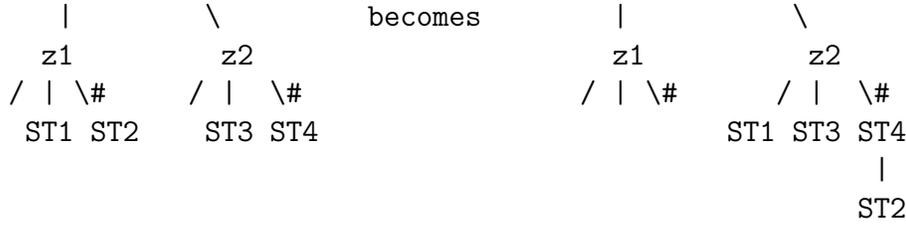
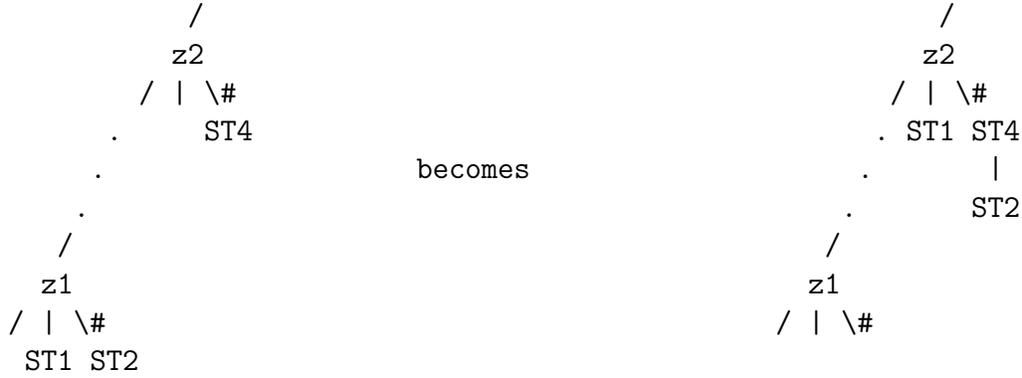
\begin{verbatim}
        |        \         becomes        |          \
       z1         z2                      z1          z2
     / | \#     / |  \#                 / | \#      / |  \#
      ST1 ST2    ST3 ST4                         ST1 ST3 ST4
                                                          |
                                                         ST2


            /                                             /
           z2                                            z2
         / | \#                                        / | \#
       .     ST4                                      . ST1 ST4
      .                   becomes                    .       |
     .                                              .       ST2
    /                                              /
   z1                                             z1
 / | \#                                         / | \#
  ST1 ST2
\end{verbatim}
Let $P$ be the prefix code described by the transformed tree.
Then $|Q| = |P|$  (since the number of vertices has not changed), and
$P_1$ has the same mod 3 cardinality as $Q_1$ (since the subtree ST1
is moved from $z_1$ to $z_2$ and $z_1$, $z_2$ have equivalent depths
mod 3). The subtree ST2 was under $\#$ and is still below a $\#$-edge.
Finally, $T_{\rm in}(P)$ has three leaves, namely $\ell_1$, $\ell_2$, and
$z_1$ (all $\in \{0,1\}^*$), and $|z_1| + 1 \equiv |q_1|$ (actually the
two numbers are equal). An existing leaf $\ell_1, \ell_2$ is either
unchanged, or (in case of a leaf in ST1) is changed to a leaf that has an
equivalent depth modulo 3, or (in case of a leaf in ST2, hence
$\notin \{0,1\}^*$) is changed to a leaf $\notin \{0,1\}^*$.

\smallskip

\noindent (2) \ If $q_1, q_2 \notin \{0,1\}^*$ and if
$q_1 \notin \{0,1\}^*\#$ (or, similarly, if $q_2 \notin \{0,1\}^*\#$),
we do the following transformation on $Q$. As above, $z_1$ and $z_2$ are the
parent vertices of $q_1$, respectively $q_2$, and ST1, ST2, ST3, ST4, ST5 are
subtrees (below $z_1$ or $z_2$) of the prefix tree of $Q$; one of ST3 and
ST5 is a single vertex (corresponding to the word $q_2$).
Again, if $z_1$ and $z_2$ are on a common path from the root we
let $z_1$ be the deeper one of the two.
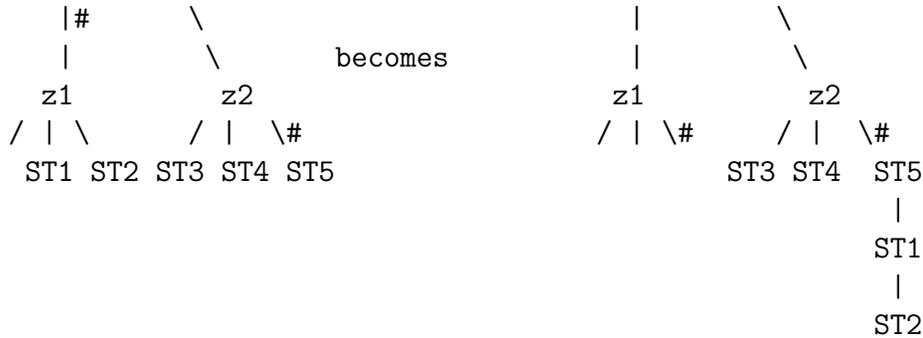
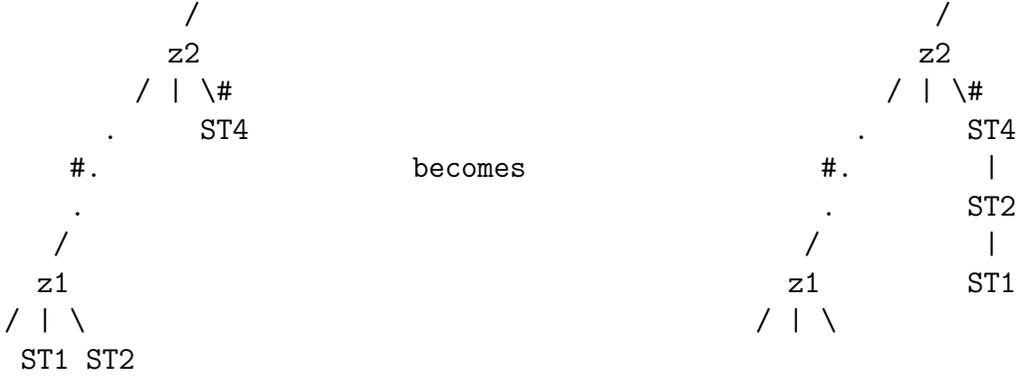
\begin{verbatim}
        |#      \                          |        \
        |        \       becomes           |         \
       z1         z2                      z1          z2
     / | \      / |  \#                  / | \#     / |  \#
      ST1 ST2 ST3 ST4 ST5                        ST3 ST4  ST5
                                                           |
                                                          ST1
                                                           |
                                                          ST2
\end{verbatim}

\newpage

\begin{verbatim}
            /                                             /
           z2                                            z2
         / | \#                                        / | \#
       .     ST4                                     .      ST4
     #.                   becomes                  #.        |
     .                                             .        ST2
    /                                             /          |
   z1                                            z1         ST1
 / | \                                         / | \
  ST1 ST2
\end{verbatim}
Then the transformed tree describes a maximal prefix code $P$ with the
desired properties. In particular, $P_1 = Q_1$, and $z_1 \notin \{0,1\}^*$
is now a leaf.

\smallskip

\noindent (3) \ Suppose $q_1, q_2 \in \{0,1\}^*\#$ and
$\ell_1 \notin \{0,1\}^*$.

\smallskip

\noindent $\bullet$ \ If $T_{\rm in}(Q)$ has other subtrees besides the
paths root-to-$\ell_1$ and root-to-$\ell_2$, then $T_{\rm in}(Q)$ has
another leaf besides $\ell_1$ and $\ell_2$. In this case we have nothing
to prove.

\smallskip

\noindent $\bullet$ \ If $T_{\rm in}(Q)$ consists only of the paths
root-to-$\ell_1$ and root-to-$\ell_2$, let $p_1$ be the parent vertex of
$\ell_1$ in $T_{\rm in}(Q)$.
Since $\ell_1 \notin \{0,1\}^*$, the path root-to-$\ell_1$ of
$T_{\rm in}(Q)$ has some edge(s) labeled by $\#$.

Since $q_1, q_2 \in \{0,1\}^*\#$, the paths root-to-$\ell_1$ or
root-to-$\ell_2$ of $T_{\rm in}(Q)$ have some edges labeled over $\{0,1\}$.
If $z_1$ and $z_2$ (the parent vertices of $q_1$ and $q_2$ in the prefix tree
of $Q$) are on a common root-to-leaf path, let $z_1$ be the name of
the deeper one of the two; then $z_2$ is at least 2 depth levels above
$\ell_1$ or $\ell_2$.
If $z_1$ and $z_2$ are on different paths root-to-leaf, let $z_2$ be on the
path root-to-$\ell_2$; then $z_2$ will be at least one depth level
above $\ell_2$ (since $q_1, q_2$ are not children of $z_1, z_2$).
We do the following transformation on $Q$.
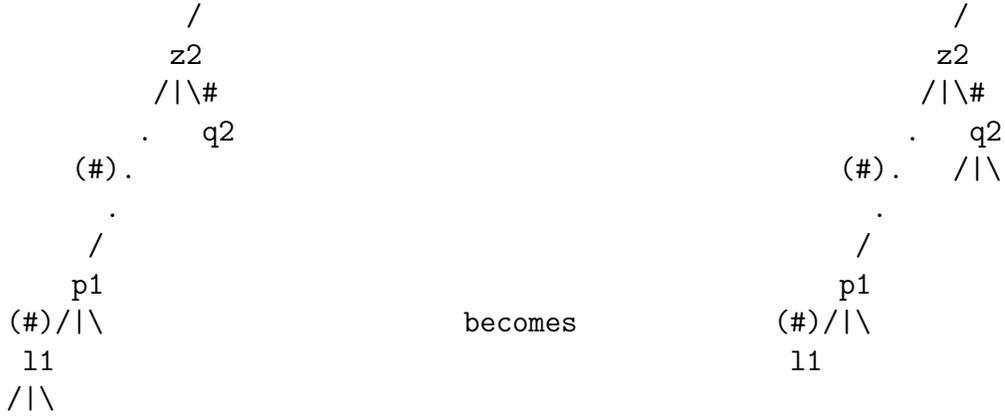
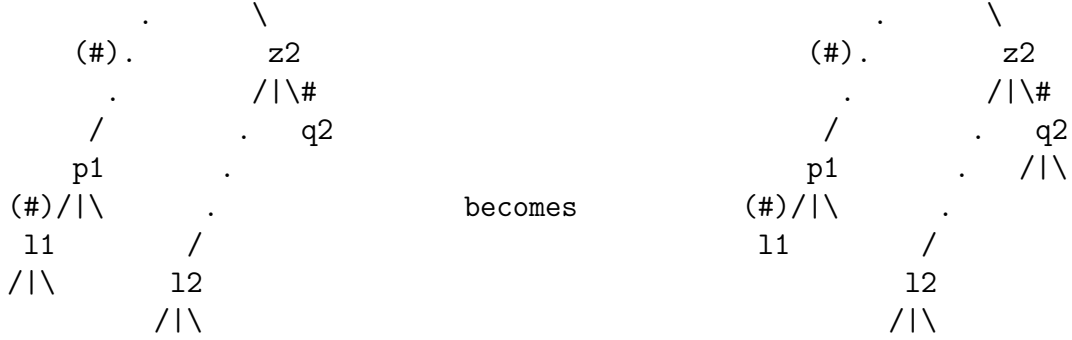
\begin{verbatim}
               /                                              /
              z2                                             z2
             /|\#                                           /|\#
            .   q2                                         .   q2
        (#).                                           (#).   /|\
          .                                              .
         /                                              /
        p1                                             p1
    (#)/|\                      becomes            (#)/|\
     l1                                             l1
    /|\
\end{verbatim}

\newpage 

\begin{verbatim}
           .      \                                     .      \
       (#).        z2                               (#).        z2
         .        /|\#                                .        /|\#
        /        .   q2                              /        .   q2
       p1       .                                   p1       .   /|\
   (#)/|\      .               becomes          (#)/|\      .
    l1        /                                  l1        /
   /|\       l2                                           l2
            /|\                                          /|\
\end{verbatim}
In other words, the children of $\ell_1$ are moved to $q_2$.
Now $|P| = |Q|$ (since no new vertices are created), and $P_1 = Q_1$,
hence $P_1$ and $Q_1$ have the same mod 3 cardinality.
Also, $T_{\rm in}(P)$  has three leaves, namely $q_2 \notin \{0,1\}^*$,
$\ell_2$, and $p_1$. If $\ell_1 \in \{0,1\}^*\#$ then $p_1 \in \{0,1\}^*$,
and if $\ell_1 \notin \{0,1\}^*\#$ (but still $\ell-1 \notin \{0,1\}^*$)
then $p_1 \notin \{0,1\}^*$.
    \ \ \ $\Box$

\begin{lem} \label{G_01_fin_genA3} \
The group $G_{3,1}^{\rm mod \, 3}(0,1)$ is generated by its elements of 
table-size $\leq c_{\rm gen}$, for some constant $c_{\rm gen}$.
\end{lem}
{\bf Proof.} \ We follow the same method as in the proof of Lemma
\ref{V_endm_fin_gen}. Let
\[ \varphi \ = \
\left[ \begin{array}{ccc ccc}
x 0 & x 1 & x \# & x_4 & \ldots & x_n \\
y_1 & y_2 & y_3  & y_4 & \ldots & y_n \end{array} \right]
 \ \in \ G_{3,1}(0,1).
\]
The image code's inner tree, $T_{\rm in}({\rm imC}(\varphi))$, has a leaf
$y$, so $\{y_1, \ldots, y_n\}$ also contains 3 words of the form  \
$y_{i_1} = y 0$, \ $y_{i_2} = y 1$, \ $y_{i_3} = y \#$, where
$y \in \{0,1,\#\}^*$. The three indices $i_1, i_2, i_3$ are in
$\{1, \ldots, n\}$, but any order relation between $i_1, i_2, i_3$ is
possible.

\medskip

\noindent {\bf Case 1:} \ The column index sets $\{1, 2, 3\}$ and
$\{i_1, i_2, i_3\}$ are disjoint.

\smallskip

\noindent
Then, after permuting columns (if necessary), the table of $\varphi$ has
the form
\[  \left[ \begin{array}{ccc ccc ccc}
x0 & x1 & x\# & x_4 & x_5 & x_6 & x_7 & \ldots & x_n \\
y_1 & y_2 & y_3 & y0 & y1 & y\# & y_7 & \ldots & y_n
 \end{array}        \right].
\]

\smallskip

\noindent {\sc Case 1.1:} \ $y \in \{0,1\}^*$. (The case where, instead,
$x \in \{0,1\}^*$ is very similar.)

Then $x_4, x_5 \in \{0,1\}^*$, and $|x_4| \equiv |x_5| \equiv |y| + 1$
mod 3, since $\varphi \in G_{3,1}^{\rm mod \, 3}(0,1)$. Then, applying Lemma
\ref{G(0,1)codesA3} (1) to the prefix code \ $Q = {\rm domC}(\varphi)$, we
obtain a maximal prefix code $P$ with the properties listed in that Lemma.
In particular, $T_{\rm in}(P)$ has two leaves, $\ell_1 \in \{0,1\}^*$,
$|\ell_1| \equiv |x|$ mod 3, and $\ell_2 \in \{0,1\}^*$ with
$|\ell_2| + 1$ $\equiv |x_4|$ $\equiv$ $|x_5| \equiv |y| + 1$ mod 3.
These properties imply that $P$ can be inserted into the table of $\varphi$
as an intermediary row, and that the columns can be lined up in such a way
that $\varphi$ is factored as two elements of 
$G_{3,1}^{\rm mod \, 3}(0,1)$:

\[  \left[ \begin{array}{ccc ccc ccc}
x0 & x1 & x\# & x_4 & x_5 & x_6 & x_7 & \ldots & x_n \\
\ell_1 0 & \ell_1 1 & \ell_1 \# & \ell_2 0 & \ell_2 1 & \ell_2 \# & z_7 &
\ldots
 & z_n \\
y_1 & y_2 & y_3 & y0 & y1 & y\# & y_7 & \ldots & y_n
 \end{array}        \right]
\]
Now, as in the proof of Lemma \ref{V_endm_fin_gen}, the two factors can
be extended, so as to get smaller tables.

\smallskip

\noindent {\sc Case 1.2:} \ Both $x, y \notin \{0,1\}^*$.
Then $x_4, x_5, y_1, y_2 \notin \{0,1\}^*$.

\smallskip

\noindent {\sc Case 1.2.1:} \ If $x_4$ or $x_5$ $\notin \{0,1\}^*\#$, then
we apply Lemma \ref{G(0,1)codesA3} (2) to the prefix code \ $Q = $
$ {\rm domC}(\varphi)$.
If $y_1$ or $y_2$ $\notin \{0,1\}^*\#$, then we apply Lemma 
\ref{G(0,1)codesA3}
(2) to the prefix code \ $Q = {\rm imC}(\varphi)$. Next, we insert $P$ into
the table of $\varphi$ in the same way as in case 1.1.

\smallskip

\noindent {\sc Case 1.2.2:} \ If $x_4, x_5, y_1, y_2 \in \{0,1\}^*\#$,
we can again apply Lemma \ref{G(0,1)codesA3} (3) to \  
$Q = {\rm domC}(\varphi)$.
If $P$ has both $\ell_1, \ell_2 \notin \{0,1\}^*$, we insert $P$ as a row,
as in case 1.1.

However, if $\ell_2 \notin \{0,1\}^*$ and $\ell_1 \in \{0,1\}^*$, we cannot
proceed as before, because both $x, y \notin \{0,1\}^*$; the resulting
factors of $\varphi$ would not stabilize $\{0,1\}^*$.
So this time we insert $P$ as two rows into the table of $\varphi$ (possibly
after permuting columns), as follows:

\[  \left[ \begin{array}{ccc ccc cc ccc}
x0 & x1 & x\# & x_4 & x_5 & x_6   & \ \ldots \ & \ \ldots \ &
                                        x_{n-2} & x_{n-1} & x_n \\
\ell_2 0 & \ell_2 1 & \ell_2 \# & \ldots & \ldots & \ldots &
   \ \ldots \ & \ \ldots \ & \ell_1 0 & \ell_1 1 & \ell_1 \#  \\
\ldots & \ldots & \ldots & \ell_2 0 & \ell_2 1 & \ell_2 \# & \ \ldots \ &
 \ \ldots \  &  \ell_1 0 & \ell_1 1 & \ell_1 \#   \\
y_1 & y_2 & y_3 & y0 & y1 & y\# & \ \ldots \  & \ \ldots \ &
  y_{n-2} & y_{n-1} & y_n
\end{array}        \right]
\]
The columns can be lined up in such a way that the three factors of
$\varphi$ belong to $G_{3,1}^{\rm mod \, 3}(0,1)$. Indeed, $|P| = |Q|$, and
$P_1$ has the same mod 3 cardinality as $Q_1$. Also, $x_4, x_5, y_1, y_2$,
$x,y, \ell_2$ $\notin \{0,1\}^*$.

\medskip

\noindent {\bf Case 2:} \ The column index sets $\{1, 2, 3\}$ and
$\{i_1, i_2, i_3\}$ overlap.

\smallskip

\noindent {\sc Case 2.1:} \ Suppose $\{1, 2\}$ overlaps with
$\{i_1,i_2,i_3\}$  {\bf and} $\{i_1,i_2\}$ overlaps with $\{1,2,3\}$.

Then both $x, y \in \{0,1\}^*$ or both $x,y \notin \{0,1\}^*$.
By Lemma \ref{G(0,1)codesA3} we find a prefix code $P$, with
$\ell_2 \in \{0,1\}^*$ if $x,y \in \{0,1\}^*$, and
$\ell_2 \notin \{0,1\}^*$ if $x,y \notin \{0,1\}^*$.
Then we insert $P$ as two rows:

\[  \left[ \begin{array}{ccc ccc cc ccc}
x0 & x1 & x\# & x_4 & x_5 & \ \ldots \ & \ \ldots \ &
                                        x_{n-2} & x_{n-1} & x_n \\
\ell_2 0 & \ell_2 1 & \ell_2 \# & \ldots & \ldots & \ \ldots \ & \ \ldots \
  & \ell_1 0 & \ell_1 1 & \ell_1 \#  \\
\ldots & \ell_2 a_1 & \ldots & \ell_2 a_2 & \ell_2 a_3 & \ \ldots \ &
 \ \ldots \  &  \ell_1 0 & \ell_1 1 & \ell_1 \#   \\
\ldots & y a_1 & \ldots & y a_2 & y a_3 & \ \ldots \  & \ \ldots \ &
  y_{n-2} & y_{n-1} & y_n
\end{array}        \right].
\]

\smallskip

\noindent {\sc Case 2.2:} \ Suppose \
$\{1, 2\} \cap \{i_1, i_2, i_3\} = \emptyset$ \ or \
$\{i_1,i_2\} \cap \{1,2,3\} = \emptyset$.

Then \ $\{1, 2, 3\} \cap \{i_1, i_2, i_3\} = \{ 3\}$ \ or \
$\{1, 2, 3\} \cap \{i_1, i_2, i_3\} = \{ i_3\}$. We only consider the case
where the intersection is $\{ 3\}$  (the case when it is $\{ i_3\}$ is
very similar).  Then the table of $\varphi$ is
\[  \left[ \begin{array}{ccccc cc }
x 0 & x 1 & x \#  & x_4   & x_5   & \ \ldots \ & \ \ldots \  \\
y_1 & y_2 & y a_1 & y a_2 & y a_3 & \ \ldots \  & \ \ldots \
\end{array}        \right]
\]

\smallskip

\noindent {\sc Case 2.2.1:} \ If $x, y \in \{0,1\}^*$, then
$ya_1  \notin \{0,1\}^*$ (since $\varphi \in G_{3,1}^{\rm mod \, 3}(0,1)$);
hence $a_1 = \#$.
Now we proceed as in case 1.1.

\smallskip

\noindent {\sc Case 2.2.2:} \ If $x, y \notin \{0,1\}^*$, then
$y_1, y_2, x_4, x_5 \notin \{0,1\}^*$. Now we proceed as in case 1.1.

\smallskip

\noindent {\sc Case 2.2.3:} \ If $x \in \{0,1\}^*$ and
$y \notin \{0,1\}^*$, then $y_1, y_2 \in \{0,1\}^*$.

We apply Lemma \ref{G(0,1)codes2A3} (1) to the maximal prefix code \
$Q = {\rm imC}(\varphi)$ \ with
existing leaf $\ell =$ $y \notin \{0,1\}^*$, and with $q_1, q_2$ equal to
$y_1, y_2 \in \{0,1\}^*$ respectively. Then we obtain a code $P$ with
$\ell_1 \notin \{0,1\}^*$, and with $\ell_2 \in \{0,1\}^*$,
 \ $|\ell_2| + 1 \equiv |y_1| \equiv |y_2|$ mod 3.
Now we insert $P$ into the table of $\varphi$ as two rows, to obtain
(after permuting columns, if necessary):
\[  \left[ \begin{array}{ccccc cc ccc ccc}
x 0 & x 1 & x \#  & x_4   & x_5   & \ \ldots \ & \ \ldots \ &
   \ x_{n-5} & x_{n-4} & x_{n-3} & x_{n-2} & x_{n-1} & x_n \\
\ell_2 0 & \ell_2 1 & \ell_2 \# & \ldots & \ldots & \ \ldots \ & \ \ldots \
  & \ell_1 0 & \ell_1 1 & \ell_1 \# & \ldots & \ldots & \ldots  \\
\ldots & \ldots & \ell_1 a_1 & \ell_1 a_2 & \ell_1 a_3 & \ \ldots \ &
 \ \ldots \ & \ldots & \ldots & \ldots & \ell_2 \# & \ell_2 1 & \ell_2 0   \\
y_1 & y_2 & y a_1 & y a_2 & y a_3 & \ \ldots \  & \ \ldots \ &
 y_{n-5} & y_{n-4} & y_{n-3} & y_{n-2} & y_{n-1} & y_n
\end{array}        \right]
\]
Here we assume that the columns $n-3$ and $n-2$ (that
contain $\ell_1 \#$, respectively $\ell_2 \#$) are disjoint.
This assumption can always be made if $n$ is large enough so that
$\{x_1, \ldots, x_n\}$ and $\{y_1, \ldots, y_n\}$ contain enough elements
$\in \{0,1\}^*$ and $\notin \{0,1\}^*$.
Then we can insert another copy of $P$ as follows:
\[  \left[ \begin{array}{ccccc cc ccc ccc}
x 0 & x 1 & x \#  & x_4   & x_5   & \ \ldots \ & \ \ldots \ &
   \ x_{n-5} & x_{n-4} & x_{n-3} & x_{n-2} & x_{n-1} & x_n \\
\ell_2 0 & \ell_2 1 & \ell_2 \# & \ldots & \ldots & \ \ldots \ & \ \ldots \
  & \ell_1 0 & \ell_1 1 & \ell_1 \# & \ldots & \ldots & \ldots  \\
\ldots & \ldots & \ldots & \ldots & \ldots & \ \ldots \ & \ \ldots \
  & \ell_1 0 & \ell_1 1 & \ell_1 \# & \ell_2 \# & \ell_2 1 & \ell_2 0 \\
\ldots & \ldots & \ell_1 a_1 & \ell_1 a_2 & \ell_1 a_3 & \ \ldots \ &
 \ \ldots \ & \ldots & \ldots & \ldots & \ell_2 \# & \ell_2 1 & \ell_2 0   \\
y_1 & y_2 & y a_1 & y a_2 & y a_3 & \ \ldots \  & \ \ldots \ &
 y_{n-5} & y_{n-4} & y_{n-3} & y_{n-2} & y_{n-1} & y_n
\end{array}        \right]
\]
This gives us a factorization of $\varphi$ as four elements of
$G_{3,1}^{\rm mod \, 3}(0,1)$, each of which can be reduced.

\smallskip

\noindent {\sc Case 2.2.4:} \ The case where $y \in \{0,1\}^*$ and
$x \notin \{0,1\}^*$ is similar to case 2.2.3, now using \
$Q = {\rm domC}(\varphi)$.
 \ \ \ $\Box$

\bigskip

In analogy with Lemma \ref{V_endm_fin_gen_rels}, Lemma \ref{G_01_fin_genA3}
can be strengthened as follows.

\begin{lem} \label{G_01_fin_gen_relsA3} \
Every element $\varphi \in G_{3,1}^{\rm mod \, 3}(0,1)$ of table-size
$> c_{\rm gen}$ can be represented by a word $w_{\varphi}$ over the set of
elements of table-size $\leq c_{\rm gen}$, and such that the sequence
$w_{\varphi}$ has table-size $\leq \|\varphi\|$. The constant $c_{\rm gen}$
is as in Lemma \ref{G_01_fin_genA3}.
\end{lem}
{\bf Proof.} \  This follows from the proof of Lemma \ref{G_01_fin_genA3}.
In that proof, we started out with a table of $\varphi$ (of table-size
$\|\varphi\|$), and repeatedly inserted rows. No columns are ever added,
hence the table-size doesn't increase.
See also the proof of Higman's Lemma 4.3 in \cite{Hig74}.
 \ \ \ $\Box$

\begin{lem} \label{G_01_fin_presA3} \
The group $G_{3,1}^{\rm mod \, 3}(0,1)$ is presented by relators of table-size
$\leq c_{\rm rel}$, in terms of generators of table-size $\leq c_{\rm gen}$,
where $c_{\rm gen}$ is the constant from Lemma \ref{G_01_fin_genA3}, and
$c_{\rm rel}$ is another constant.
Hence, $G_{3,1}^{\rm mod \, 3}(0,1)$ is finitely presented.
\end{lem}
{\bf Proof.} \ We use the same approach as in Proposition
\ref{V_endm_fin_pres} (based on Higman's proof  that $G_{N,r}$ is finitely
presented (see \cite{Hig74}, pp.~29-33).
We now use Lemma \ref{G(0,1)codes2A3}.

For the same reason as in Lemma \ref{G_01_fin_genA3}, the new rows that are
inserted have their columns lined up in such a way that all pairs of adjacent
rows represent elements of $G_{3,1}^{\rm mod \, 3}(0,1)$ (and not just 
of $G_{3,1}$).

Higman's ``type II'' reductions (described in the figure in the top
of p.~31 of \cite{Hig74}) can actually be replaced by his ``type III
reductions'' (described in the figure in the top of p.~32 of \cite{Hig74}).
Type II reductions are never needed (the reason why they were used by
Higman is probably that they are more efficient: they require a single
row insertion; on the other hand, a type III reduction consists of two
transformations).

Type III reductions require that we insert a row corresponding to a prefix
code with 3 leaves in the inner tree (see the figure at the top of p.~31 in
\cite{Hig74}). Since one of the pre-existing rows in the table already has
two leaves ($b$ and $c$ in Higman's notation), we want the table size to be
large enough so that the maximal prefix code
 \  $Q = \{b0, b1, b\#, c0, c1, c\#, \ldots \}$ \
(2nd row of figure at bottom of p.~31, and 2nd row of figure at top of p.~32
 in \cite{Hig74}) contains either another leaf in its inner tree or two words
that are not children of a leaf of the inner tree. In the latter case we
apply Lemma \ref{G(0,1)codes2A3} and obtain a maximal prefix code $P$ with
three leaves $x,y,z$ with $x$ equivalent to $b$ and
$y$ {\it equivalent} to $c$. (Here we define two words $u, v \in \{0,1,\#\}^*$
to be {\it equivalent} iff both $u,v \notin \{0,1\}^*$ or both
$u,v \in \{0,1\}^*$ and $|u| \equiv |v|$  mod 3.) Also, $P$ and $Q$ have the
same cardinality, and $P_1$ and $Q_1$ have the same mod 3 cardinality.
Therefore, we can insert a row corresponding to the prefix code $P$ in
exactly the same way as on p.~32 of \cite{Hig74}, taking care to line up the
columns so that the factors belong to $G_{3,1}^{\rm mod \, 3}(0,1)$.
 \ \ \ $\Box$

\bigskip

Just as for $G_{3,1}^{\rm mod \, 3}(0,1;\#)$, one can prove that
$G_{3,1}^{\rm mod \, 3}(0,1)$ is not a simple group; a very similar 
homomorphic image can be taken.
We can summarize the results for $G_{3,1}^{\rm mod \, 3}(0,1)$ as follows.

\begin{thm} \
The group $G_{3,1}^{\rm mod \, 3}(0,1)$ is finitely presented, and not simple.

The word problem of $G_{3,1}^{\rm mod \, 3}(0,1)$, over the generating set
 \ $\Delta_{0,1} \cup \{\tau_{i,i+1} : 0 \leq i\}$ is coNP-hard, with 
respect to constant-arity conjunctive polynomial-time reduction.
Here $\Delta_{0,1}$ is a finite generating set of 
$G_{3,1}^{\rm mod \, 3}(0,1)$.
\end{thm}
 
\noindent
As a consequence of Proposition \ref{V_endm_closed} we can consider the 
following HNN-extension:

\medskip

$H(0,1) \ = \ $
$\langle G_{3,1}^{\rm mod \, 3}(0,1) \cup \{t\} \ : \ $
  $\{t \, g \, t^{-1} = g^{\kappa_{321}} : g \in $
  $G_{3,1}^{\rm mod \, 3}(0,1) \} \rangle$.

\medskip

\noindent Since $G_{3,1}^{\rm mod \, 3}(\{0,1\}^*)$ is finitely generated, 
the HNN-relations form a finite set; moreover, since 
$G_{3,1}^{\rm mod \, 3}(0,1)$ is finitely presented (by teh above Theorem), 
the whole HNN-extension is a finitely presented group

For the same reason as for $G_{3,1}^{\rm mod \, 3}(0,1;\#)$ in Section 7, 
we obtain:

\begin{lem} \label{H_t_kappaA3}  \
The HNN-extension $H(0,1)$ is isomorphic to the subgroup \
$\langle G_{3,1}^{\rm mod \, 3}(0,1) \cup \{\kappa_{321}\} \rangle$ \ of the
Thompson group ${\mathcal G}_{3,1}$.
\end{lem}
In summary, we obtain Theorem \ref{finprescoNPhard}, as well as the other 
main theorems, for $G_{3,1}^{\rm mod \, 3}(0,1)$ and $H(0,1)$. 

Section 8 shows that the word problem of $H(0,1)$ (over a finite generating 
set) is in coNP.

%%%%%%%%%%%%%%%%%%%%%%%%%%%%%%%%%%%%%%%%%%%%%%%%%%%%%%%%%

\subsection{Miscellaneous}   % A4

The following is a converse of Proposition \ref{V_endm_closed}. This 
converse gives an interesting property of $G_{3,1}^{\rm mod \, 3}(0,1)$ (that 
$G_{3,1}^{\rm mod \, 3}(0,1;\#)$ does not have), but we make no use of it in 
this paper.

\begin{pro} \label{Converse_V_endm_closed} \
Let us abbreviate $\kappa_3 \kappa_2 \kappa_1(\cdot)$ to $\kappa$.
If $g \in G_{3,1}$ is such that the conjugates of $g$ under $\kappa$ or
$\kappa^{-1}$ belong to $G_{3,1}$, then $g \in G_{3,1}^{\rm mod \, 3}(0,1)$. 
In other words,

\smallskip

$G_{3,1}^{\rm mod \, 3}(0,1) = \{ g \in G_{3,1} : $
      $g^{\kappa}, g^{\kappa^{-1}} \in G_{3,1} \}$.
\end{pro}
{\bf Proof.} \  Suppose $g \in G_{3,1}$ and
$g^{\kappa}, g^{\kappa^{-1}} \in G_{3,1}$.

\smallskip

\noindent {\sc Claim 1:} \ $g \in {\rm Stab}(\{0,1\}^*)$.

\smallskip

\noindent Proof of Claim 1: \ By contraposition we assume that
$g \notin {\rm Stab}(\{0,1\}^*)$, hence $g^{-1} \notin$
${\rm Stab}(\{0,1\}^*)$, and we will prove that $g^{\kappa} \notin G_{3,1}$.

If $g, g^{-1} \notin {\rm Stab}(\{0,1\}^*)$
then (perhaps after replacing $g$ by $g^{-1}$), there is
$z \in \{0,1,\#\}^* - \{0,1\}^*$ such that $z \in {\rm Dom}(g)$ and
$g(z) \in \{0,1\}^*$; let $y = g(z)$.
Let $x = k^{-1}(z)$; note that $z \in {\rm Dom}(\kappa^{-1})$ since
$z$ contains the letter $\#$. Then $x$ contains $\#$ too, so
$x \in {\rm Dom}(\kappa)$; moreover, for all $v \in \{0,1,\#\}^*$, \
$\kappa(xv) = \kappa(x) \ v = zv$.
For any $w \in \{0,1\}^*$ we have

\smallskip

$x w \# \ \stackrel{\kappa}{\longmapsto} \ \kappa(x) \ w \# \ = \ zw\# $
$ \stackrel{g}{\longmapsto} \ g(z) \ w\# $
$ \ = \ yw \#  \ \ \ (\in \{0,1\}^* \#)$

\smallskip

$ \stackrel{\kappa^{-1}}{\longmapsto} \ \kappa^{-1}(yw \, \#)$.

\smallskip

\noindent We want to show now that ${\rm domC}(g^{\kappa})$ is infinite 
(when $g^{\kappa}$ is maximally extended). Assume by contradiction that
${\rm domC}(g^{\kappa})$ is finite; so the elements of 
${\rm domC}(g^{\kappa})$ have length $< b$ for some constant $b$.

Recall the definition of $\kappa$ and its relation with the permutation
$\gamma_3 \gamma_2 \gamma_1(\cdot)$ of $\mathbb{N}$, described in the
beginning of the paper:

\smallskip

$\gamma_3 \gamma_2 \gamma_1(\cdot) \ = \ $

\smallskip

 \ \ \ \ \
$( \ \ldots \ | \ 6(j+1) \ | \ 6j \ | \ \ldots \ | $
$\ 12 \ | \ 6 \ {\bf | \ 2 \ |} \ 5 \ | \ 8 \ | \ \ldots \ $
$ | \ 3i+2 \ | \ 3(i+1)+2 \ | \ \ldots \ ) \ \cdot \ $

\smallskip

 \ \ \ \ \
$( \ \ldots \ | \ 6(j+1)+3 \ | \ 6j+3 \ | \ \ldots \ | $
$ \  9 \ | \ 3 \ {\bf | \ 1 \ |} \ 4 \ | \ 7 \ | \ \ldots \ | $
$ \ 3i+1 \ | \ 3(i+1)+1 \ | \ \ldots \ )  (\cdot)$.

\smallskip

Therefore, the application of $\kappa^{-1}$ to $yw\#$ changes bit number
$3(i-1)+2$ to bit number $3i+2$ of $yw$,
for every $i$, \ $0 \leq i \leq (|yw|-2)/3$.
Let us pick a $w \in \{0,1\}^*$ which is much longer than $b$.
Then \  $g^{\kappa}(xw\#)$ $(= \kappa^{-1}(yw\#))$ \ cannot we written in
the form \ $g^{\kappa}(xw\#) = g^{\kappa}(uv\#) = g^{\kappa}(u) \, v\#$,
for any factorization of $xw\#$ as $xw\# = uv\#$ with $|u| < b$.
This contradicts the assumption that the elements of ${\rm domC}(g^{\kappa})$
have length $< b$.

Therefore, $g^{\kappa}$ does not belong to $G_{3,1}$. This proves Claim 1.

\smallskip

\noindent From here on we can assume that $g \in {\rm Stab}(\{0,1\}^*)$.

\medskip

\noindent {\sc Claim 2:} \ $g \in G_{3,1}^{\rm mod \, 3}$.

\smallskip

\noindent Proof of Claim 2: \ By contraposition we assume that
$g \notin G_{3,1}^{\rm mod \, 3}$, hence 
$g^{-1} \notin G_{3,1}^{\rm mod \, 3}$;
we will prove that $g^{\kappa} \notin G_{3,1}$.

If $g, g^{-1} \notin G_{3,1}^{\rm mod \, 3}$ then there is $z \in \{0,1\}^*$ 
such that $z \in {\rm Dom}(g)$, \, $g(z) \in \{0,1\}^*$, and
$|g(z)| \not\equiv |z|$ mod 3.

Since the action of $\kappa$ consists of permuting bits over a distance
$\leq 6$ we have the following:
There exist  $x, s \in \{0,1\}^*$, with $|s| \leq 6$ and $|x| = |z|$, such
that \ $\kappa(xs\#) = zt\#$ for some $t \in \{0,1\}^*$ with $|t| = |s|$.

For any $w \in \{0,1\}^*$ we have \ $\kappa(xsw\#) = zvw'\#$, for some
$v,w' \in \{0,1\}^*$ with $|v| = |s|$, $|w'| = |w|$, and where $w'$ depends
only on $s$ and $w$ (and not on $x$).
Let $y = g(z)$. Then we have:

\smallskip

$xsw\# \ \stackrel{\kappa}{\longmapsto} \ zvw'\# $
$ \stackrel{g}{\longmapsto} \ g(z) \ vw'\# $
$ \ = \ y vw' \#  \ \ \ (\in \{0,1\}^* \#)$

\smallskip

$ \stackrel{\kappa^{-1}}{\longmapsto} \ \kappa^{-1}(yvw' \, \#)$.

\smallskip

\noindent where $|y| \not\equiv |x|$ mod 3, and $|x| = |z|$.

Then $\kappa^{-1}(yvw'\#) = y'v'w''\#$, for some $y', v', w'' \in \{0,1\}^*$
with $|y'| = |y|$, $|v'| = |v|$ and $|w''| = |w'| = |w|$. However,
since $|y| \not\equiv |x| = |z|$, it follows that $w''$ differs from $w$
in every bit position.

We want to show now that ${\rm domC}(g^{\kappa})$ is infinite (when
$g^{\kappa}$ is maximally extended). Assume by contradiction that
${\rm domC}(g^{\kappa})$ is finite; so the elements of 
${\rm domC}(g^{\kappa})$ have length $< b$ for some constant $b$.
Let us pick a $w \in \{0,1\}^*$ which is much longer than $b$.

Since the application of $g^{\kappa}$ to $xsw\#$ changes all the bits
$w$, it follows that $g^{\kappa}(xsw\#)$ cannot we written in
the form \ $g^{\kappa}(xsw\#) = g^{\kappa}(uv\#) = g^{\kappa}(u) \, v\#$,
for any factorization of $xsw$ as $xsw = uv$ with $|u| < b$.
This contradicts the assumption that the elements of ${\rm domC}(g^{\kappa})$
have length $< b$.

Therefore, $g^{\kappa}$ does not belong to $G_{3,1}$. This proves Claim 2.
 \ \ \ $\Box$

\bigskip

%%%%%%%%%%%%%%%%%%%%%%%%%%%%%%%%%%%%%%%%%%%%%%%%%%%%%%%%%%%%%%%%%%%%%%%%%%

\bigskip

\bigskip

\noindent
{\bf Jean-Camille Birget} \\
Dept.\ of Computer Science \\
Rutgers University at Camden \\
Camden, NJ 08102, USA \\
{\tt birget@camden.rutgers.edu }


\begin{thebibliography}{99}

\bibitem{Ben73} C. Bennett, ``Logical reversibility of computation'', 
{\it IBM J.\ Research and Development} 17 (1973) 525-532.

\bibitem{Ben89} C. Bennett, ``Time/Space tradeoffs for reversible
computation'', {\it SIAM J. of Computing} {\bf 18} (1989) 766-776.

\bibitem{Bi} J.C.\ Birget, ``Time-complexity of the word problem for
  semigroups and the Higman Embedding Theorem",
  {\it International J. of Algebra and Computation} 8 (1998) 235-294.

\bibitem{BiRedFunct} J.C.\ Birget, ``Reductions and functors from problems 
  to word problems", {\it Theoretical Computer Science} 237 (2000) 81-104. 

\bibitem{BiFonG} J.C.\ Birget, ``Functions on groups and computational 
  complexity'', {\it International J. of Algebra and Computation}, to 
  appear. (Mathematics ArXiv: math.GR/0202124) 

\bibitem{BiThomps} J.C.\ Birget, ``The groups of Richard Thompson and 
  complexity'', {\it International J. of Algebra and Computation}, to appear.
  (Mathematics ArXiv: math.GR/0204292)

\bibitem{BORS} J.C.~Birget, A.~Ol'shanskii, E.~Rips, M.V.~Sapir,
  ``Isoperimetric functions of groups and computational complexity of the
  word problem'', {\it Annals of Mathematics} 156.2 (Sept.~2002) 467-518. 
  (Mathematics arXiv, math.GR/9811106,   http://front.math.ucdavis.edu)

\bibitem{Bridson} M.~Bridson, ``The Geometry of the Word Problem '', 
   in {\it Invitations to Geometry and Topology}, Oxford University Press, 
   2002.   

\bibitem{BradyBridson} N.~Brady, M.~Bridson, ``There is only one gap in 
   the isoperimetric spectrum'', {\it GAFA} 10 (2000) 1053-1070.

\bibitem{CFP} J.\ W.\ Cannon, W.\ J.\ Floyd, W.\ R.\ Parry,
``Introductory notes on Richard Thompson's groups'',
{\it L'Enseignement Math\'ematique} 42 (1996) 215-256.
 
\bibitem{FredToff} E.~Fredkin, T.~Toffoli, ``Conservative logic'', 
{\it International J.\ Theoretical Physics} 21 (1982) 219-253.

\bibitem{GZ} M.\ Garzon, Y.\ Zalcstein, ``The complexity of Grigorchuk
groups with application to cryptography'', {\it Theoretical Computer Science}
88 (1991) 83-98.

\bibitem{Gromov} M.\ Gromov, ``Asymptotic invariants of infinite groups'',
in {\it Geometric Group Theory} (G.\ Niblo, M.\ Roller, editors),
London Mathematical Society Lecture Notes Series 182, Cambridge Univ.\
Press (1993).

\bibitem{Hig74} G.\ Higman, ``Finitely presented infinite simple groups'',
Notes on Pure Mathematics 8, The Australian National University,
Canberra (1974).

\bibitem{Lec} Y. Lecerf, ``Machines de Turing r\'{e}versibles ...'',
{\em Comptes Rendus de l'Acad\'{e}mie des Sciences, Paris} 257 No. 18 (Oct.
1963) 2597 - 2600.

\bibitem{LiptZalc} R.\ Lipton, Y.\ Zalcstein, ``Word problems solvable in
log space'', {\it Journal of the Association for Computing Machinery} 24
(1977) 522-526.

\bibitem{LyndonSchupp} R.\ Lyndon, P.\ Schupp, {\it Combinatorial Group
Theory}, Springer-Verlag (1977).

\bibitem{MO} K.\ Madlener, F. Otto, ``Pseudo-natural algorithms for the
word problem for finitely presented monoids and groups'', {\it J.\ of
Symbolic Computation} 1 (1985) 383-418.

\bibitem{MagnusKaSo} W.~Magnus, A.~Karrass, D.~Solitar, {\it Combinatorial
Group Theory}, Dover 1976 (Interscience 1966).

\bibitem{McKTh} R.\ McKenzie, R.\ J.\ Thompson,
``An elementary construction of unsolvable word problems in group theory'',
in {\it Word Problems}, (W.\ W.\ Boone, F.\ B.\ Cannonito, R.\ C.\ Lyndon,
editors), North-Holland (1973) pp.\ 457-478.

\bibitem{Olsh} A.Y.\ Ol'shanskii, ``On subgroup distortion in finitely
presented groups'', {\it Matematicheskii Sbornik} 188 (1997) 51-98.

\bibitem{OlSap} A.Y.\ Ol'shanskii, M.V.\ Sapir, ``Length and area functions
on groups and quasi-metric Higman embedding'', {\it International J.\ of
Algebra and Computation} 11 (2001) 137-170.
 
\bibitem{Reidemeister} K.~Reidemeister, {\it Einf\"uhrung in die 
kombinatorische Topologie}, Chelsea, New York 1950 (Vieweg, Braunschweig 1932).

\bibitem{RoeverGR} C.\ R\"over, ``Constructing finitely presented simple groups
that contain Grigorchuk groups'', {\it  J.\ of Algebra} 220 (1999) 284-313.

\bibitem{SBR} M.V.~Sapir, J.C.~Birget, E.~Rips, ``Isoperimetric and
 isodiametric functions of groups", {\it Annals of Mathematics} 156.2
 (Sept.~2002) 345-466.  (Mathematics arXiv, math.GR/9811105,
  http://front.math.ucdavis.edu)

\bibitem{JESavage} J.E.\ Savage, {\it Models of Computation}, Addison-Wesley 
(1998).

\bibitem{ESc} Elizabeth A. Scott, ``A construction which can be used
to produce finitely presented infinite simple groups'',
{\it J. of Algebra} 90 (1984) 294-322.

\bibitem{EScConjug} Elizabeth A. Scott, ``A finitely presented simple group
with unsolvable conjugacy problem'', {\it J. of Algebra} 90 (1984) 333-353.

\bibitem{EScSurvey} Elizabeth A. Scott, ``A tour around finitely presented
simple groups'', in {\it Algorithms and Classification in Combinatorial
Group Theory} (G.\ Baumslag, Ch.F.\ Miller III, editors), MSRI Publications 23,
Springer-Verlag (1992).

\bibitem{Th0} Richard J. Thompson, Manuscript (1960s).

\bibitem{Th} Richard J. Thompson, ``Embeddings into finitely generated
simple groups which preserve the word problem'',
in {\it Word Problems II}, (S.\ Adian, W.\ Boone, G.\ Higman, editors),
North-Holland (1980) pp.\ 401-441.

\bibitem{Handb} J.\ van Leeuwen (editor), {\it Handbook of Theoretical
Computer Science}, volume {\bf A}, MIT Press and Elsevier (1990).

\bibitem{Wegener} I.\ Wegener, {\it The complexity of boolean functions},
Wiley/Teubner (1987).


\end{thebibliography}
\end{document}